\documentclass[12pt, a4paper]{amsart}

\usepackage[english]{babel}
\usepackage{amsfonts}
\usepackage{amsmath}
\usepackage{amssymb, latexsym}
\usepackage{graphicx}
\usepackage[usenames]{color}

\usepackage{fullpage}
\usepackage{amsmath}

\usepackage{mathrsfs}

\newtheorem{theorem}{Theorem}[section]
\newtheorem{proposition}[theorem]{Proposition}
\newtheorem{corollary}[theorem]{Corollary}
\newtheorem{lemma}[theorem]{Lemma}
\newtheorem{remark}[theorem]{Remark}
\newtheorem{definition}[theorem]{Definition}

\numberwithin{equation}{section} 

\numberwithin{equation}{section}

\newcommand{\R}{\mathbb{R}}

\newcommand{\ben}{\begin{eqnarray*}}
\newcommand{\enn}{\end{eqnarray*}}
\newcommand{\pa}{\partial}

\newcommand{\g}{\gamma}

\newcommand{\al}{\alpha}
\newcommand{\la}{\lambda}

\newcommand{\ol}{\overline}

\newcommand{\half}{\frac{1}{2}}

\newcommand{\na}{\nabla}
\newcommand{\be}{\begin{equation}}
\newcommand{\ee}{\end{equation}}
\newcommand{\ba}{\begin{aligned}}
\newcommand{\ea}{\end{aligned}}
\newcommand{\lf}{\left}
\newcommand{\rt}{\right}

\newcommand{\wt}{\tilde{w}}

\def\9{{\infty}}
\def\a{{\alpha}}
\def\b{{\beta}}
\def\g{{\gamma}}
\def\del{{\delta}}
\def\lbb{{\lambda}}
\def\t{{\theta}}

\def\calo{{\mathcal{O}}}
\def\calp{{\mathcal{P}}}

\def\calu{{\mathcal{U}}}
\def\calx{{\mathcal{X}}}
\def\bbp{{\mathbb{P}}}
\def\bbr{{\mathbb{R}}}
\def\ve{{\varepsilon}}
\def\vf{{\varphi}}

\def\wt{\widetilde}

\def\ol{\overline}
\def\({\left(}
\def\){\right)}
\def\<{\langle}
\def\>{\rangle}

\begin{document}

\title{Minimal mass blow-up solutions to rough nonlinear Schr\"odinger equations}

\author{Yiming Su}
\address[Y.M. Su]{School of science,
Zhejiang University of Technology, 310023 Hangzhou, China.}
\author{Deng Zhang}
\address[D. Zhang]{School of mathematical sciences,
Shanghai Jiao Tong University, 200240 Shanghai, China.}
\begin{abstract}
  We study the focusing mass-critical rough nonlinear Schr\"odinger equations,
  where the stochastic integration is taken in the sense of controlled rough path.
  We obtain the global well-posedness
  if the mass of initial data is below that of the ground state.
  Moreover, the existence of minimal mass blow-up solutions
  is also obtained in both dimensions one and two.
  In particular, these yield that the mass of ground state is
  exactly the threshold of global well-posedness and blow-up of solutions
  in the stochastic focusing mass-critical case.
  Similar results are also obtained for a class of
  nonlinear Schr\"odinger equations with lower order perturbations.
\end{abstract}


\maketitle
\section{Introduction} \label{Sec-Intro}
We are concerned with
the  focusing mass-critical rough nonlinear Schr\"{o}dinger equations
with linear multiplicative conservative noise
and the blow-up dynamics of their solutions.
Precisely,
consider the equation
\begin{align} \label{equa-X-rough}
     &  dX = i\Delta Xdt + i|X|^{\frac 4d}X dt - \mu X dt + XdW(t),  \\
     &  X(0)= X_0 \in H^1(\bbr^d). \nonumber
\end{align}
Here
$$W(t,x)=\sum_{k=1}^N i\phi_k(x)B_k(t),\ \ x\in \bbr^d,\ t\geq 0,$$
where
$\phi_k \in C_b^\9(\bbr^d, \bbr)$,
$B_k$ are the standard $N$-dimensional real valued Brownian motions
on a stochastic basis $(\Omega, \mathscr{F}, \{\mathscr{F}_t\}, \bbp)$,
$1\leq k\leq N <\9$,
and $\mu= \frac 12 \sum_{k=1}^N  \phi_k ^2$.
The last term $XdW(t)$ in \eqref{equa-X-rough} is taken in the sense of controlled rough path,
see Definition \ref{def-X-rough} below.
In particular,
the rough integration coincides with the usual It\^o integration
if the corresponding processes are $\{\mathscr{F}_t\}$-adapted
(see \cite[Chapter 5]{FH14}).

The physical significance of \eqref{equa-X-rough} is well known.
In the conservative case considered here  (i.e., ${\rm Re}W=0$),
$\|X(t)\|^2_{L^2}$ is pathwisely conserved.
Hence,
with the normalization $\|X_0\|_{L^2}=1$,
the quantum system evolves on the unit ball of $L^2$
and so verifies the conservation of probability.
Moreover,
the noise in \eqref{equa-X-rough} can be viewed as a random potential
acting in the system.
In crystals it corresponds to scattering of excitons by phonons,
because of thermal vibrations of molecules,
see \cite{BCIR94, BCIRG95} in the one dimensional case,
and also \cite{RGBC95} for the two dimensional case.
Another important application can be found in the quantum measurement in open quantum systems (\cite{BG09}).
We  also refer to \cite{SS99} for more physical applications in the deterministic case,
such as nonlinear optics,
Bose-Einstein condensation and
the Gross-Pitaevskii equation.

The blow-up phenomena
are extensively studied in the literature for the
focusing mass-critical
nonlinear Schr\"{o}dinger equation (NLS):
\begin{align}
 &i\partial_tu+\Delta u+|u|^{\frac 4d} u=0,  \label{equa-nls}  \\
 &u(0)=u_{0} \in H^1(\bbr^d). \nonumber
\end{align}
As a matter of fact, equation \eqref{equa-nls} admits a number of symmetries and conservation laws.
Precisely,
it is invariant under the translation, scaling, phase rotation and Galilean transform,
i.e.,
if $u$ solves (\ref{equa-nls}),
then so does
\be\label{symmerty}
v(t,x)=\lambda_0^{-\frac{d}{2}} u\lf(\frac{t-t_0}{\lambda_0^2},\frac{x-x_0 }{\lambda_0}- \frac{\beta_0(t-t_0)}{\lbb_0}\rt)e^{i\frac{\beta_0}{2}\cdot(x-x_0)-i\frac{|\beta_0|^2}{4}(t-t_0)+i\t_0},
\ee
with $v(t_0,x)=\lambda_0^{-\frac{d}{2}} u_0\lf(\frac{x-x_0}{\lambda_0}\rt)e^{i\frac{\beta_0}{2}\cdot(x-x_0)+i\t_0}$,
where
$(\lambda_0, \beta_0, \theta_0) \in \bbr^+ \times \bbr^d \times \bbr$,
$x_0\in\R^d$, $t_0\in\R$.
In particular,
the $L^2$-norm of solutions is  preserved under the symmetries above,
and thus
\eqref{equa-nls} is  called the mass-critical equation.
Another important symmetry is related to the
pseudo-conformal transformation in the pseudo-conformal space
$\Sigma:=\{u\in H^1(\bbr^d), \||x|u\|_{L^2(\bbr^d)}<\9 \}$,
\be\label{pseudo}
 \frac{1}{(-t)^{\frac{d}{2}}}u\left(\frac{1}{-t},\frac{x}{-t}\right)e^{-i\frac{|x|^{2}}{4t}}, \quad\ \ t\not =0.
\ee
The conservation laws related to \eqref{equa-nls} contain
\begin{align}
   & {\rm Mass}:\ \ M(u)(t):=\int_{\R^d}|u(t)|^2dx=M(u_0), \label{mass}  \\
   & {\rm Energy}:\ \ E(u)(t) := \half\int_{\R^d}|\nabla u(t)|^2dx-\frac{d}{2d+4}\int_{\R^d}|u(t)|^{2+\frac{4}{d}}dx=E(u_0), \label{energy} \\
   & {\rm Momentum}: \ \ Mom(u):={\rm Im} \int_{\R^d}\nabla u\bar{u}dx
       ={\rm Im} \int_{\R^d}\nabla u_{0}\bar{u}_{0}dx. \label{momentum}
\end{align}

It is  well known that  equation \eqref{equa-nls} is locally well posed
and the solutions exist globally if
the mass of initial data is below that of the ground state $Q$,
which is the unique positive spherically symmetric solution to the soliton equation
\begin{align} \label{equa-Q}
    \Delta Q - Q + Q^{1+\frac{4}{d}} =0.
\end{align}

However,
the situation becomes much more delicate
if $\|u_0\|_{L^2}=\|Q\|_{L^2}$.
In this case, equation (\ref{equa-nls}) has two special solutions:
the solitary wave solution $u(t,x)=Q(x)e^{it}$ that exists globally,
and the so-called pseudo-conformal blow-up solution $S_T(t,x)$, $T\in\bbr$,
obtained from the pseudo-conformal transformation (\ref{pseudo})
and the solitary wave solution:
\be \label{S}
S_T(t,x) : =\frac{1}{(T-t)^{\frac{d}{2}}}Q\left(\frac{x}{T-t}\right)e^{\frac{i}{T-t}-i\frac{|x|^{2}}{4(T-t)}}, \ \ x\in \bbr^d, \ t<T.
\ee
Note that, $\|S_T\|_{L^2}=\|Q\|_{L^2}$,
and $S_T$ blows up at time $T$ with
the blow-up speed $(T-t)^{-1}$.

In particular,
the mass of ground state
$\|Q\|_{L^2}$ characterizes the threshold for the global well-posedness
and blow-up of solutions to NLS,
and $S_T$ is exactly the minimal mass blow-up solution.

It should be mentioned that,
minimal mass blow-up solutions are of significant importance  in the
study of blow-up phenomena.
A remarkable result proved by Merle (\cite{Mu}) is that,
up to the symmetries \eqref{symmerty},
the pseudo-conformal blow-up solution is the unique
minimal mass blow-up solution to NLS.
Moreover,
for the inhomogeneous NLS,
a sufficient condition for the nonexistence of minimal mass blow-up solutions
was proved in \cite{M96}.
Under suitable degenerate conditions,
Banica, Carles and Duyckaerts \cite{BCD} constructed the minimal mass blow-up solutions.
The sharp non-degenerate conditions
for the existence of minimal mass blow-up solutions,
as well as the strong rigidity theorem of uniqueness,
have been
proved by Rapha\"el and Szeftel \cite{R-S}.

For minimal mass blow-up solutions of other dispersive equations,
see, e.g., \cite{K-L-R} for the mass-critical fractional NLS,
\cite{CS13} for the inhomogeneous Hartree equation,
\cite{S17,SG19} for the nonlinear Schr\"{o}dinger system,
and \cite{MM02} for the mass-citical gKdV equation.
We also refer to
\cite{B-W, F17, M-R, P}
for the case of supercritical mass (i.e., $\|u_0\|_{L^2} >\|Q\|_{L^2})$,
where a new type of blow-up solutions
with log-log blow-up rate has been investigated.

For  stochastic Schr\"odinger equations (with It\^o's integration),
there are many works devoted to the  well-posedness results.
See, e.g. \cite{BRZ16,BD03,H18,BM13} for the subcritical case,
and also the recent works \cite{FX18.1,FX18.2,Z18} for the defocusing mass- and energy-critical cases.
As regards
the blow-up phenomena,
de Bouard and Debussche (\cite{BD02,BD05}) first studied the noise effect on blow-up.
Quite surprisingly,
the conservative noise
has the effect to accelerate blow-up immediately with positive probability
in the focusing mass-supercritical case,
namely, the nonlinearity in \eqref{equa-X-rough}
has the exponent $\a\in (1+\frac 4d, 1+\frac {4}{d-2})$ (\cite{BD02,BD05}).
The proof is based on the fact that,
the non-degenerate noise
is able to push the solution
to a blow-up regime,
in which the energy is sufficiently negative
to formalize singularity.

Furthermore,
in the focusing mass-critical case,
the numerical results in \cite{BDM01,DM02,DM02.2} suggest that
the colored conservative noise has a tendency to delay blow-up,
and the white noise may even prevent blow-up.
We also refer to \cite{BRZ17} for the damped effect on blow-up
of the non-conservative noise
(i.e., ${\rm Re} W \not = 0$)
in the mass-(super)critical case.

However, to the best of our knowledge,
quite few results are known
for the quantitative description of blow-up dynamics
in the stochastic setting.
It is unclear whether $\|Q\|_{L^2}$ still serves as
the threshold for the global well-posedness and blow-up
of solutions in the stochastic situation.
One main difficulty here
lies in the loss of the symmetries \eqref{symmerty}
and  the conservation law of energy \eqref{energy},
which are, actually, completely destroyed because of the presence of noise.

Here  we
obtain
the global well-posedness of equation \eqref{equa-X-rough}
if the mass of initial data is below that of the ground state.
Moreover, in both dimensions one and two,
we construct the
minimal mass blow-up solutions to equation \eqref{equa-X-rough},
which indeed evolve asymptotically like the
pseudo-conformal blow-up solutions near the blow-up time.
In particular,
these results yield
that the mass of the ground state is exactly the threshold for the
global well-posedness and blow-up of solutions to \eqref{equa-X-rough}.

Furthermore,
similar results are also obtained
for a class of nonlinear Schr\"odinger equations
with lower order perturbations,
which
has not been well studied
because of the
loss of symmetries \eqref{symmerty} and the conservation law of energy.

The idea of proof is mainly based on the rescaling approach as in \cite{BRZ14, BRZ16,BRZ18}
and on the modulation method  developed in \cite{R-S},
which includes
a bootstrap device and the backward propagation from the singularity.
The latter enables to reduce the analysis of blow-up solutions of \eqref{equa-X-rough}
to that of the dynamics of finite dimensional geometrical parameters
and a small remainder corresponding to
a random nonlinear Schr\"odinger equation with lower order perturbations
(see \eqref{equa-u-RNLS} below).
In particular,
it provides  quantitative descriptions of blow-up dynamics
in the absence of symmetries and conservation of energy.

It should be mentioned that,
because of the backward propagation procedure,
the blow-up solution constructed is no longer
adapted,
and thus the last term in equation \eqref{equa-X-rough}
should not be taken in the sense of It\^o.
It turns out that
it can be appropriately interpreted in the sense of controlled rough path.
Similar situation arises also in the random three dimensional vorticity equation in \cite{BR17,RZZ19}.

A key role is then played by the
equivalence between solutions to the rough equation \eqref{equa-X-rough}
and the random equation \eqref{equa-u-RNLS}.
Quite differently from \cite{BRZ14, BRZ16},
the proof here requires finer descriptions of time regularities of solutions
and, actually, the local smoothing space
is also applied to measure the spatial regularity of solutions
in the two dimensional case.  \\

{\bf Notations.}
For any $x=(x_1,\cdots,x_d) \in \bbr^d$
and any multi-index $\nu=(\nu_1,\cdots, \nu_d)$,
let $|\nu|= \sum_{j=1}^d \nu_j$,
$\<x\>=(1+|x|^2)^{1/2}$,
$\partial_x^\nu=\partial_{x_1}^{\nu_1}\cdots \partial_{x_d}^{\nu_d}$,
and
$\<\na\>=(I-\Delta)^{1/2}$.

For $1\le p\le\9$,
$L^p = L^p(\bbr^d)$ is
the space of $p$-integrable (complex-valued) functions
endowed with the norm $\|\cdot\|_{L^p}$,
and $W^{s,p}$ denotes the standard Sobolev space,
$s\in \bbr$.
In particular,
$L^2(\bbr^d)$ is the Hilbert space endowed with the scalar product
$\<v,w\> =\int_{\bbr^d} v(x)\bar w(x)dx$,
and $H^s:= W^{s,2}$.
As usual,
$L^q(0,T;L^p)$ means the space of all integrable $L^p$-valued functions $f:(0,T)\to L^p$ with the norm
$\|\cdot\|_{L^q(0,T;L^p)}$,
and $C([0,T];L^p)$ denotes the space of all $L^p$-valued continuous functions on $[0,T]$ with the sup norm over $t$.
We also use the local smoothing spaces defined by
$L^2(I;H^\a_{\beta})=\{u\in \mathscr{S}': \int_{I} \int \<x\>^{2\beta}|\<\na\>^{\a} u(t,x)|^2  dxdt <\9 \}$,
where $\a, \beta \in \mathbb{R}$.

For any H\"older continuous function $f\in C^\a(I)$, $I\subseteq \bbr^+$,
we write $\delta f_{st} := f(t)-f(s)$, $s,t\in I$,
and $\|f\|_{\a, I} := \sup_{s,t\in I,s\not =t} \frac{|\delta f_{st}|}{|s-t|^\a}$.
Let $C_c^\9$ be the space of all compactly supported smooth functions on $\bbr^d$.
We also set $g_t:= \frac{d}{dt}g$ for any $C^1$ functions.

The symbol $u =
\mathcal{O}(v)$ means that $|u/v|$ stays bounded, and $v_n=o(1)$
means that $|v_n|$ tends to zero as $n\to \9$.
Throughout this paper,
we use $C$ for various constants that may
change from line to line.

\section{Formulation of main results} \label{Sec-Mainresults}

To begin with, we first present the precise definition of solutions to equation \eqref{equa-X-rough}.

\begin{definition} \label{def-X-rough}
We say that $X$ is a solution to \eqref{equa-X-rough} on $[0,\tau^*)$,
where $\tau^*\in (0,\9]$ is a random variable,
if $\bbp$-a.s. for any $\vf\in C_c^\9$,
$t \mapsto \<X(t), \vf\>$ is continuous on $[0,\tau^*)$
and for any $0<s<t<\tau^*$,
\begin{align*}
   \<X(t)-X(s), \vf\>
   - \int_s^t \<i X, \Delta\vf\>  + \<i|X|^{\frac 4d} X, \vf\>  - \< \mu X, \vf\> dr
   = \sum\limits_{k=1}^N \int_s^t \<i\phi_k X, \vf\> dB_k(r).
\end{align*}
Here the integral $\int_s^t \<i\phi_k X, \vf\> d B_k(r)$
is taken in the sense of controlled rough path
with respect to the rough paths $(B, \mathbb{B})$,
where $\mathbb{B}=(\mathbb{B}_{jk})$, $\mathbb{B}_{jk,st}:= \int_s^t \delta B_{j,sr} dB_k(r)$
with the integration taken in the sense of It\^o.
That is,
$\<i\phi_k X, \vf\> \in C^\a([s,t])$,
\begin{align} \label{phikX-st}
   \delta (\<i\phi_k X, \vf\>)_{st}
   = -\sum\limits_{j=1}^N \<\phi_j\phi_k X(s), \vf\> \delta B_{j,st}
     + \delta R_{k,st},
\end{align}
and
$\|\<\phi_j\phi_k X, \vf\> \|_{\a, [s,t]} <\9, \ \
   \|R_k\|_{2\a, [s,t]} <\9.
$
\end{definition}

We assume that
\begin{enumerate}
   \item[(A0)] ({\it Asymptotical flatness})
  For $1\leq k\leq N$, $\phi_k$ satisfies that
   for any multi-index $\nu \not = 0$,
\begin{align} \label{decay}
   \lim_{|x|\to \9} \<x\>^2 |\partial_x^\nu \phi_k(x)| =0.
\end{align}

  \item[(A1)] ({\it Flatness at the origin})
  For  $1\leq k\leq N$ and for any multi-index $0\leq|\nu|\leq 5$,
\begin{align} \label{degeneracy}
   \partial_x^\nu \phi_k(0)=0.
\end{align}
\end{enumerate}

\begin{remark}
Assumption $(A0)$ is  related to the asymptotical behavior at infinity
of the spatial functions of noise
and guarantees the well-posedness of equation \eqref{equa-X-rough}.
While,
Assumption $(A1)$
characterizes the local behavior near the origin
and is mainly imposed for the construction of minimal mass blow-up solutions.
\end{remark}

Theorem \ref{Thm-X-LWP} summarizes the local well-posedness and blow-up alternative results.

\begin{theorem} \label{Thm-X-LWP}
Let $d\geq 1$.
Assume $(A0)$.
Then, for each $X_0\in H^1$, $\bbp-a.s.$
there exists a unique solution $X$ to \eqref{equa-X-rough} on $[0,\tau^*)$
in the sense of
Definition \ref{def-X-rough},
where $\tau^* \in (0,\9]$ is a random variable,
such that $\bbp$-a.s. for any $T<\tau^*$,
\begin{equation} \label{X-stri-LS-bdd}
    X|_{[0,T]}\in C([0,T];H^1)\cap L^\g(0,T;W^{1,\rho})
                   \cap L^2(0,T; H^\frac 32_{-1}),
\end{equation}
where $(\rho,\g)$ is any Strichartz pair,
i.e., $\frac{2}{\g} = d(\frac 12- \frac{1}{\rho})$,
$(\rho, \g,d) \not= (\9, 2,2)$.

Moreover, we have that for $\bbp$-a.e. $\omega$,
either $\tau^*(\omega)=\9$, or
\begin{equation} \label{subcritical-a}
    \lim\limits_{t\to \tau^*(\omega)} \|X(\omega, t)\|_{H^1}=\9.
\end{equation}
\end{theorem}

\begin{remark}
The proof of Theorem \ref{Thm-X-LWP} is similar to that in \cite{BRZ16},
based on the Strichartz estimates for Schr\"odinger equations with lower order perturbations
and the equivalence between solutions to equations \eqref{equa-X-rough} and \eqref{equa-u-RNLS} in Theorem \ref{Thm-Equiv-X-u} below.
\end{remark}

The next result is concerned with the global well-posedness of \eqref{equa-X-rough}
if the mass of initial data is below that of the ground state.

\begin{theorem} \label{Thm-X-GWP}
Let $d\geq 1$.
Assume $(A0)$.
If $\|X_0\|_{L^2} < \|Q\|_{L^2}$,
then the corresponding solution $X$ to \eqref{equa-X-rough} exists globally almost surely.
\end{theorem}

The main result concerning the existence of minimal mass blow-up solutions is formulated below.
\begin{theorem} \label{Thm-X-blowup}
Let $d=1,2$.
Assume Assumptions $(A0)$ and $(A1)$ to hold.
Then,
for $\bbp$-a.e. $\omega$
there exists $\tau^*(\omega) \in (0,\9)$,
such that for any $T\in (0,\tau^*(\omega)]$
there exists $X_0(\omega) \in H^1$
satisfying that
$\|X_0(\omega)\|_{L^2} = \|Q\|_{L^2}$
and  the corresponding solution $X(\omega)$ to \eqref{equa-X-rough}
blows up at time $T$.
Moreover,
there exist $\delta, C(\omega, T) >0$ such that for $t$ close to $T$,
\begin{align} \label{asym-X-ST}
     \|X(\omega, t) - e^{W(\omega, t)} S_T(t) \|_{H^1} \leq C(\omega, T)(T-t)^\delta.
\end{align}
\end{theorem}

\begin{remark}
$(i)$ Theorems \ref{Thm-X-GWP} and \ref{Thm-X-blowup}
show that, in dimensions one and two,
the mass of the ground state $\|Q\|_{L^2}$ is indeed
the threshold for the global well-posedness and blow-up of solutions to
the rough nonlinear Schr\"odinger equation \eqref{equa-X-rough}.

$(ii)$
The asymptotic behavior \eqref{asym-X-ST} yields that,
up to a phase shift related to the noise,
the pseudo-conformal blow-up solutions are exactly the main blow-up profile
near the blow-up time.
We also would like to mention that,
it is possible to construct  solutions that blow up
at finitely many points,
which will be done in the forthcoming work.

$(iii)$
It is interesting to see whether uniqueness holds for the minimal mass blow-up solutions
to \eqref{equa-X-rough}
which, however, is still unclear.
One difficulty arises from the loss of
the conservation law of energy
due to the presence of noise.
\end{remark}

The first step of the proof is   based on the rescaling approach.
More precisely,
equation \eqref{equa-X-rough} can be formally transformed,
via $u=e^{-W}X$,
to the random equation
\begin{align}  \label{equa-u-RNLS}
\ &i\partial_tu+ e^{-W}\Delta(e^{W}u)+|u|^{\frac{4}{d}}u =0, \;\;(t,x)\in\R\times\R^d,  \\
\ &u(0,x)=u_{0}(x)   \nonumber
\end{align}
with $u_0 = X_0$ and $e^{-W}\Delta(e^{W}u) = \Delta u+ b \cdot \nabla u+c u$,
where
\begin{align}
 b(t,x)&=2\na W(t,x) = 2i \sum\limits_{k=1}^N \na \phi_k(x) B_k(t),  \label{b} \\
 c(t,x)&=\sum\limits_{j=1}^d (\partial_j W(t,x))^2 + \Delta W(t,x)  \nonumber \\
       &= - \sum\limits_{j=1}^d (\sum\limits_{k=1}^N \partial_j \phi_k(x) B_k(t))^2
         + i\sum\limits_{k=1}^N \Delta \phi_k(x) B_k(t).    \label{c}
\end{align}

The transformation $u=e^{-W}X$
is known as the Doss-Sussman transformation in finite dimensional case.
It turns out to be quite robust also in the study of stochastic   equations
in infinite dimensional spaces.
One advantage of the rescaling approach is,
that it enables one to reduce
the problem of stochastic equations to that
of equations with random coefficients,
to which one is able to obtain sharp estimates in a pathwise manner.
See, e.g., \cite{BRZ17.0} for stochastic logarithmic Schr\"odinger equations,
\cite{HRZ18} for the scattering behavior and \cite{BRZ18,Z19} for optimal control problems.
See also \cite{BR17,RZZ19} for random three dimensional vorticity equations.

The solutions to the random equation \eqref{equa-u-RNLS}
can be defined similarly as in \cite{BRZ16}.

\begin{definition} \label{def-u-RNLS}
We say that $u$ is a solution to \eqref{equa-u-RNLS} on $[0,\tau^*)$,
where $\tau^* \in (0,\9]$ is a random variable,
if $\bbp-a.s.$  $u\in C([0,\tau^*); H^1)$,
$|u|^{\frac{4}{d}}u\in L^1(0,\tau^*; H^{-1})$,
and $u$ satisfies
\begin{align} \label{equa-u-RNLS-def}
     u(t)= u(0) + \int_0^t i e^{-W(s)}\Delta(e^{W(s)}u(s))+i|u(s)|^{\frac{4}{d}}u(s) ds,  \ \
      t\in [0,\tau^*),
\end{align}
as an equation in $H^{-1}$.
\end{definition}

\begin{remark} \label{Rem-u-weak-mild}
Under Assumption $(A0)$,
the solution to \eqref{equa-u-RNLS} in the sense of Definition \ref{def-u-RNLS}
is equivalent
to that taken in the mild sense below
\begin{align*}
   u(t) = e^{it\Delta} u(0) + \int_0^t e^{i(t-s)\Delta} (i|u(s)|^{\frac{4}{d}}u(s) + i(b(s)\cdot \na + c(s))u(s)) ds, \ \
   \forall t\in [0,\tau^*),
\end{align*}
as an equation in $H^{-1}$.
Moreover,
similarly to Theorem \ref{Thm-X-LWP},
by \cite{BRZ16} one also has the  local well-posedness, blow-up alternative results
and that for any $0<T<\tau^*$ and for any Strichartz pair $(\rho, \g)$,
\begin{align} \label{u-Stri-ls-bdd}
   \|u\|_{C([0,T]; H^1)} +\|u\|_{L^\g(0,T; W^{1,\rho})} +\|u\|_{L^2(0,T; H^{\frac 32}_{-1})} <\9, \ \ \bbp-a.s..
\end{align}
\end{remark}

An important role here
is played by the equivalence between solutions to the rough equation \eqref{equa-X-rough}
and the random equation \eqref{equa-u-RNLS}.

\begin{theorem} \label{Thm-Equiv-X-u}
$(i)$. Let $u$ be the solution to \eqref{equa-u-RNLS} on $[0,\tau^*)$
with $u(0) = u_0\in H^1$
in the sense of Definition \ref{def-u-RNLS},
where $\tau^* \in (0,\9]$ is a random variable.
Then, $\bbp$-a.s.,
$X := e^{W} u$ is the solution to equation \eqref{equa-X-rough} on $[0,\tau^*)$
with $X(0)= u_0$
in the sense of Definition \ref{def-X-rough}.

$(ii)$. Let $X$ be the solution to equation \eqref{equa-X-rough} on $[0,\tau^*)$ with $X(0)=X_0\in H^1$
in the sense of Definition \ref{def-X-rough},
satisfying that $\bbp$-a.s.
$\|X\|_{C([0,T]; H^1)} + \|X\|_{L^2(0,T; H^{\frac 32}_{-1})} <\9$,
$T\in (0,\tau^*)$,
and
$$\|e^{-it\Delta}e^{-W(t)}X(t) - e^{-is\Delta}e^{-W(s)}X(s)\|_{L^2} \leq C(t)(t-s),\ \ \forall 0\leq s<t<\tau^*.$$
Then,
$u:=e^{-W} X$ solves equation \eqref{equa-u-RNLS} on $[0,\tau^*)$ with $u(0)=X_0$
in the sense of Definition \ref{def-u-RNLS}.
\end{theorem}

Similar results were proved in \cite{BRZ16}
where It\^o's integration was considered.
Quite differently,
the proof of Theorem \ref{Thm-Equiv-X-u}
requires  finer descriptions of the time regularity of solutions,
and the local smoothing space will be also applied
to measure the spatial regularity of solutions in
the two dimensional case.
The proof is contained in Section \ref{Sec-Proof-Main}.

By virtue of Theorem \ref{Thm-Equiv-X-u},
we can now reduce the proof of Theorems \ref{Thm-X-GWP} and \ref{Thm-X-blowup}
to that of the following two results corresponding to the random equation \eqref{equa-u-RNLS}.

\begin{theorem} \label{Thm-u-GWP}
Let $d\geq 1$. Assume $(A0)$.
If $\|u_0\|_{L^2} < \|Q\|_{L^2}$,
then the corresponding solution $u$ to \eqref{equa-u-RNLS} exists globally almost surely.
\end{theorem}

\begin{theorem}\label{Thm-u-blowup}
Let $d=1,2$.
Assume Assumptions $(A0)$ and $(A1)$ to hold.
Then, for $\bbp$-a.e. $\omega$
there exists $\tau^*(\omega) \in (0,\9)$
such that for any $T\in (0,\tau^*(\omega)]$,
there exists $u_0(\omega)\in H^1$
satisfying
$\|u_0(\omega)\|_{L^2} = \|Q\|_{L^2}$
and the corresponding  solution $u(\omega)$ to \eqref{equa-u-RNLS}
blows up at time $T$.
Moreover,  there exist $\delta, C(\omega, T)>0$ such that for $t$ close to $T$,
\begin{align}
     \|u(\omega, t) - S_T(t) \|_{H^1} \leq C(\omega, T)(T-t)^\delta.
\end{align}
\end{theorem}

Theorem \ref{Thm-u-GWP} is proved by the analysis of the evolution of energy
and the sharp Gagliardo-Nirenberg inequality obtained in \cite{Wenn}.

Moreover,
the strategy of the proof of Theorem \ref{Thm-u-blowup} is mainly
inspired by the modulation method in the recent work \cite{R-S}
and the
compactness arguments in \cite{Mu}.
It should be mentioned that,
the main challenge here
lies in the absence of the symmetries and the conservation of energy.

As a matter of fact,
the proof in the original work \cite{Mu}
relies crucially on the pseudo-conformal symmetry,
which,
however, is lost in equation \eqref{equa-u-RNLS} because of the lower order perturbations.
Moreover, quite different from the exponential decay setting in \cite{Mu},
a polynomial decay problem arises from the error term
if one linearizes equation \eqref{equa-u-RNLS} around
the pseudo-conformal blow-up solution.

Here, we apply the framework recently developed by Rapha\"{e}l and Szeftel \cite{R-S}
in the setting of inhomogeneous nonlinear Schr\"odinger equations.
This framework is, actually, quite robust
and applies to more general situations,
in which the symmetries may be lost.
However,
unlike in \cite{R-S},
the conservation law of energy  is also destroyed in equation \eqref{equa-u-RNLS}.
This leads us to a finer control of the energy
and to the spatial flatness condition \eqref{degeneracy}.
It turns out that,
the flatness at the origin of the spatial functions of noise
is sufficient to
decouple the interactions between
the lower order perturbations and the main blow-up profile.
This also reflects the nature of spatial localization in the construction of minimal mass blow-up solutions.

Below
let us give a brief outline of the main steps in the proof of Theorem \ref{Thm-u-blowup}.

We first construct a sequence of approximating solutions
$u_n$
with   $u(t_n)= S_T(t_n)$,
where $t_n$ is any sequence  converging  to $T(>0)$.
The key point here is
that each $u_n$ admits a geometrical decomposition
\be\label{un-dec*}
    u_n(t,x)
    =\lbb_n^{-\frac{d}{2}}(Q_{\calp_n} + \ve_n)(t,\frac{x-\a_n}{\lbb_n}) e^{i\theta_n}
    (=:w_n(t,x) + R_n(t,x)),
\ee
where
$Q_{\calp_n}(t,y) := Q(y) e^{i\beta_n(t) \cdot y - i \frac{\g_n(t)}{4} |y|^2}$
is a two parameters deformation of the ground state $Q$,
and the geometrical parameters $\calp_n:=(\lambda_n(t), \alpha_n(t), \beta_n(t), \gamma_n(t), \theta_n(t))$
and the remainder $R_n$ satisfy the appropriate orthogonality conditions in \eqref{ortho-cond-Rn-wn} below,
which are closely related to the generalized kernels of the linearized operators
$L_{+}=-\Delta + I -(1+{\frac{4}{d}})Q^{\frac{4}{d}}$ and
$L_{-}= -\Delta +I -Q^{\frac{4}{d}}$.

Hence,
the analysis of blow-up dynamics
is now reduced to that
of the finite dimensional geometrical parameters and the remainder term.
Plugging \eqref{un-dec*} into \eqref{equa-u-RNLS}
we then obtain the  equation for the remainder
\begin{align} \label{equa-R}
   i\partial_t R_n+\Delta R_n+(|w_n+R_n|^{\frac{4}{d}}(w_n+R_n)-|w_n|^{\frac{4}{d}}w_n)+b \cdot \nabla R_n+c  R_n=-\eta_n,
\end{align}
where
\be \label{etan-Rn}
    \eta_n = i\partial_t w_n+\Delta w_n+|w_n|^{\frac{4}{d}}w_n+b \cdot \nabla w_n+c w_n.
\ee
(See also the equation of $\ve_n$ in \eqref{equa-ve} below.)
Then,
taking the inner product of equation \eqref{equa-R}
with the directions spanning the generalized null space of $L_{\pm}$,
one is able to obtain the crucial estimates of modulation equations
driving the geometrical parameters.
See Section \ref{Sec-Approx-Modulat} below.

The next step in Section \ref{Sec-Unif-Esti}
is devoted to the key uniform estimates of approximating solutions,
to which a bootstrap device and
the  backward propagation  from the singularity
have been applied.
It should be mentioned that,
the crucial quadratic terms of $R_n$
are controlled
by the energy  \eqref{energy}
and a Lyapounov functional  \eqref{def-I}
involving a local Morawetz type term.
Moreover,
the orthogonality conditions enable us to obtain the negligible errors
of the linear terms of $R_n$,
while the terms of order higher than two
can be easily controlled by using
Gagliardo-Nirenberg's inequality.

Finally, in Section \ref{Sec-Proof-Main},
using compactness arguments we are able to extract a
limit $u_0$ at the initial time $0$,
which then yields a solution
by the  well-posedness result.
In view of the fact that
the pseudo-conformal blow-up solutions
are the main profile in the construction of approximating solutions,
we thus obtain that the resulting solution is
exactly a minimal mass blow-up solution.

It should be mentioned that,
because of the backward propagation procedure above,
the existence of the limit $u_0$
requires the information of the whole sample paths of Brownian motions
on $[0,T)$,
and thus the resulting blow-up solution  is no longer $\{\mathscr{F}_t\}$-adapted.
Instead of  It\^o's integration,
the last term in equation \eqref{equa-X-rough}
can be reinterpreted in the sense of controlled rough path.
Actually, these two notions of stochastic integrations
coincide with each other if the adaptedness of processes is fulfilled.
The proof of the key equivalence between solutions
to equations \eqref{equa-X-rough} and \eqref{equa-u-RNLS}
is  also contained in Section \ref{Sec-Proof-Main}.

The remainder of this paper is structured as follows.
Section \ref{Sec-Pre} contains some preliminaries,
including the coercivity of operators $L_{\pm}$,
the expansion formulas
and basic notions of controlled rough path.
Section \ref{Sec-Approx-Modulat} and \ref{Sec-Unif-Esti}
are devoted to the approximating solutions
and constitute the most technical part of this paper.
In Section \ref{Sec-Approx-Modulat},
we obtain the geometrical decomposition of approximating solutions
and also the estimates of modulation equations.
Then, in Section \ref{Sec-Unif-Esti},
we derive the key uniform estimates
of approximating solutions,
including the analysis of the energy and the Lyapounov functional
involving a local Morawetz type term.
The proof of main results are contained in Section \ref{Sec-Proof-Main}.
Finally,
some technical proofs are postponed to the Appendix,
i.e., Section \ref{Sec-Appendix} for simplicity.

\section{Preliminaries}  \label{Sec-Pre}

This section contains the preliminaries needed in the proof,
including the
coercivity of linearized operators,
the expansion formulas
and basic notions of controlled rough path.

\subsection{Coercivity of linearized operators} \label{Subsec-Coer}

We denote  $Q$ the ground state that solves  the soliton equation \eqref{equa-Q}.
It follows from \cite[Theorem 8.1.1]{C} that
$Q$ is smooth and decays at infinity exponentially fast,
i.e., there exist $C, \delta>0$ such that for any multi-index $|\nu|\leq 2$,
\be\label{Q-decay}
|\partial_x^\nu Q(x)|\leq C e^{-\delta |x|}, \ \ x\in \bbr^d.
\ee

Let $L=(L_+,L_-)$ be the linearized operator around the ground state state,
defined by
\begin{align} \label{L+-L-}
     L_{+}:= -\Delta + I -(1+{\frac{4}{d}})Q^{\frac{4}{d}}, \ \
    L_{-}:= -\Delta +I -Q^{\frac{4}{d}}.
\end{align}
The generalized null space of operator $L$ is
spanned by $\{Q, xQ, |x|^2 Q, \na Q, \Lambda Q, \rho\}$,
where
$\Lambda := \frac{d}{2}I + x\cdot \na$,
and $\rho$ is the unique $H^1$ spherically symmetric solution to the equation
\begin{align} \label{def-rho}
L_{+}\rho= - |x|^2Q,
\end{align}
which satisfies the exponential decay property (see, e.g., \cite{K-M-R, MP18}),
i.e., for some $C,\delta>0$,
\ben
|\rho(x)|+|\nabla \rho(x)|
\leq Ce^{-\delta|x|}.
\enn
Moreover, we have (see, e.g., \cite[(B.1), (B.10), (B.15)]{Wenm})
\be \ba \label{Q-kernel}
&L_+ \na Q =0,\ \ L_+ \Lambda Q = -2 Q,\ \ L_+ \rho = -|x|^2 Q, \\
&L_{-} Q =0,\ \ L_{-} xQ = -2 \na Q,\ \ L_{-} |x|^2 Q = - 4 \Lambda Q.
\ea\ee

For any complex valued $H^1$ function $f = f_1 + i f_2$
in terms of the real and imaginary parts,
we set
\be
(Lf,f) :=\int f_1L_+f_1dx+\int f_2L_-f_2dx.
\ee
Let $\mathcal{K}$ denote the
set of all complex valued $H^1$ functions $f=f_1+if_2$
satisfying the orthogonality conditions below
\be\ba\label{ortho-cond}
&\int Qf_1dx=0,\;\; \int xQf_1dx=0,\;\;\int |x|^2Qf_1dx=0,\\
&\int \nabla Qf_2dx=0,\;\;\int \Lambda Qf_2dx=0,\;\;
\int \rho f_2dx=0.
\ea\ee

The coercivity property below is crucial in the proof of Theorem \ref{Thm-u-blowup}.
\begin{lemma} {(\cite[Theorem 2.5]{Wenm})}  \label{Lem-coerc-L}
We have that for some $\nu>0$,
\be
(Lf,f)\geq \nu\|f\|_{H^1}^2,\ \ \forall f\in \mathcal{K}.
\ee
\end{lemma}

As a consequence we have
\begin{corollary}\label{Cor-coer-f-H1}
There exist positive constants $\nu_1, \nu_2>0$,  such that
\begin{align} \label{coer-f-H1}
(Lf,f)\geq&  \nu_1\|f\|_{H^1}^2
             -\nu_2\big(\<f_1,Q\>^2+\<f_1,xQ\>^2+\<f_1,|x|^2Q\>^2 \nonumber  \\
          &\qquad \qquad \qquad +\<f_2,\nabla Q\>^2+\<f_2,\Lambda Q\>^2+\<f_2,\rho\>^2\big),\ \ \forall f \in H^1,
\end{align}
where $f_1$ and $f_2$ are the real and imaginary parts of $f$, respectively.
\end{corollary}

The following localized version of the coercivity property
will be also useful.
\begin{corollary}(Localized coercivity)\label{Cor-coer-f-local}
Let $\Phi$ be a positive smooth radial function on $\R^d$,
such that
$\Phi(x) = 1$ for $|x|\leq 1$,
$\Phi(x) = e^{-|x|}$ for $|x|\geq 2$,
$0<\Phi \leq 1$,
and $\lf|\frac{\nabla\Phi}{\Phi}\rt|\leq C$ for some $C>0$.
Set $\Phi_A(x) :=\Phi\lf(\frac{x}{A}\rt)$, $A>0$.
Then,
for $A$ large enough we have
\be\ba \label{coer-f-local}
\int|\nabla f|^2\Phi_A +|f|^2-(1+\frac 4d)Q^{\frac4d}f_1^2-Q^{\frac4d}f_2^2dx\geq \nu \int(|\nabla f|^2+|f|^2)\Phi_A dx,\
\forall f\in \mathcal{K},
\ea\ee
where $\nu>0$, and $f_1, f_2$ are the real and imaginary parts of $f$, respectively.
\end{corollary}
The proofs of Corollaries \ref{Cor-coer-f-H1} and \ref{Cor-coer-f-local}
are postponed to the Appendix for simplicity.

\subsection{Expansion of the nonlinearity} \label{Subsec-Expan}
Let $d=1,2$.
Let
$f(z):= |z|^{\frac 4d}z$, $z\in \mathbb{C}$,
and for $v,R\in \mathbb{C}$
set
\begin{align}
   f'(v)\cdot R :=&  \partial_z f(v) R+ \partial_{\ol{z}} f(v) \ol{R}
                 = (1+\frac 2d) |v|^{\frac 4d} R + \frac 2d |v|^{\frac 4d-2} v^2 \ol{R},   \label{f-linear} \\
   N_{f,2}(v,R) :=& \frac 12 \partial_{zz}f(v)R^2 + \partial_{z\ol{z}} f(v) |R|^2 + \frac 12 \partial_{\ol{z}\ol{z}} f(v) \ol{R}^2 \nonumber \\
                 =& \frac 1d(1+\frac 2d)|v|^{\frac 4d -2} \ol{v} R^2
                    + \frac 2d (1+\frac 2d) |v|^{\frac 4d -2} v |R|^2
                    +   \frac 1d(\frac 2d-1) |v|^{\frac 4d -4} v^3 \ol{R}^2.  \label{f-quadratic} \\
   N_f(v,R):=& f(v+R)-f(v) - f'(v)\cdot R. \label{Nf-def}
\end{align}
Then, expanding $f(v+R)$ around $v$ we have
\begin{align}
   f(v+R) =&  f(v) + f'(v)\cdot R
         +  N_{f,2}(v,R)
         + \calo(\sum\limits_{k=3}^{1+\frac 4d} |v|^{1+\frac 4d-k} |R|^k). \label{f-Taylor} \\
   N_f(v,R) =& \calo(\sum\limits_{k=2}^{1+\frac 4d} |v|^{1+\frac 4d-k} |R|^k). \label{Nf-Taylor}
\end{align}

Similarly, for $F(z):= \frac{d}{2d+4} |z|^{2+ \frac 4d}$, $z\in \mathbb{C}$,
we have the expansion
\begin{align} \label{F-Taylor*}
  F(u)=& F(v)
         + \frac 12 |v|^{\frac 4d}\ol{v} R
         + \frac 12  |v|^{\frac 4d} v \ol{R} \nonumber \\
       & + \frac{1}{2d} |v|^{\frac 4d -2}\ol{v}^2 R^2
         + \frac 12 (1+\frac 2d) |v|^{\frac 4d} |R|^2
         + \frac{1}{2d} |v|^{\frac 4d -2} v^2 \ol{R}^2
         + \calo(\sum\limits_{k=3}^{2+\frac 4d} |v|^{2+\frac 4d -k} |R|^k).
\end{align}

The following lemma contains the well-known Gagliardo-Nirenberg  inequality,
which is useful to control the remainder terms in the expansion formulas.
\begin{lemma} \label{Lem-GN}
(\cite[Theorem 1.3.7]{C})
Let $d\geq 1$ and $2\leq p< \infty$.
Then, we have
\begin{align}  \label{G-N}
\|f\|_{L^p}\leq C \|f\|_{L^2}^{1 - d (\frac 12-\frac 1p)} \|\nabla f\|_{L^2}^{d(\frac 12-\frac 1p)}, \ \ \forall f\in H^1,
\end{align}
where $C>0$. Moreover,
in the case where $d=1$,
\begin{align}
     \|f\|_{L^\9} \leq C \|f\|_{H^1}, \ \ \forall f\in H^1.
\end{align}
\end{lemma}

\subsection{Controlled rough path} \label{Subsec-Roughpath}
Below we briefly recall the basic notions of controlled rough path
that are needed in the proof.
For more details of the theory of (controlled) rough path we refer to \cite{G04, MP18} and the references therein.

Given a path $X \in C^\a([0,T]; \bbr^N)$, $0<T<\9$,
we say that $Y \in C^\a([0,T]; \bbr^N)$ is controlled by $X$ if
there exists $Y' \in C^\a([0,T]; \bbr^{N\times N})$ such that
the remainder term $R^Y$ implicitly given by
\begin{align*}
   \delta Y_{j,st} = \sum\limits_{k=1}^N Y'_{jk}(s) \delta X_{k,st} + \delta R^Y_{j,st}
\end{align*}
satisfies $\|R_j^Y\|_{2\a, [s,t]} <\9$,
$1\leq j\leq N$.
This defines the controlled rough path $(Y,Y')\in \mathscr{D}_X^{2\a}([0,T]; \bbr^N)$,
and $Y'$ is the so called Gubinelli's derivative.

Note that, the $N$-dimensional Brownian motions
$B=(B_j)_{j=1}^N$ can be enhanced to a rough path
${\bf B} = (B, \mathbb{B})$,
where $\mathbb{B}_{jk,st}:= \int_s^t \delta B_{j,sr} dB_k(r)$
with the integration taken in the sense of It\^o.
It is known (\cite[Section 3.2]{FH14}) that $\|B\|_{\a,[0,T]} <\9$, $\|\mathbb{B}\|_{2\a, [s,t]}<\9$, $\bbp$-a.s.,
where $\frac 1 3 <\a<\frac 12$.

Given a path $Y$ controlled by the $N$-dimensional Brownian motion,
i.e., $Y\in \mathscr{D}_B^{2\a}([S,T]; \bbr^N)$,
$0<S<T<\9$,
we can define the rough integration of $Y$ against ${\bf B}=(B,\mathbb{B})$ as follows
(see \cite[Theorem 4.10]{FH14}),
for each $1\leq k\leq N$,
\begin{align} \label{def-rp}
   \int_S^T Y_k(r) dB_k(r)
   := \lim\limits_{|\mathscr{P}|\to 0} \sum\limits_{i=0}^{n-1}
       (Y_k(t_i) \delta B_{k,t_it_{i+1}}
         + \sum\limits_{j=1}^N Y'_{kj}(t_i) \mathbb{B}_{jk,t_it_{i+1}}),
\end{align}
where $\mathscr{P}:= \{t_0, t_i,\cdots, t_n\}$ is a partition of $[S,T]$
so that $t_0 =S$, $t_n=T$,
and $|\mathscr{P}|:= \max_{0\leq i\leq n-1} |t_{i+1} - t_i|$.

\section{Approximating solutions and modulation parameters} \label{Sec-Approx-Modulat}

\subsection{Approximating solutions}
Below we fix $T\in (0,\9)$,
which may be taken sufficiently small and will be determined later.
Let $t_n$, $n\geq 1$, be
an increasing sequence such that $\lim_{n\to \infty}t_n=T$
and $u_n$ be the corresponding solution to the equation
\begin{align}\label{equa-un-tn}
 &i\partial_tu_n+\Delta u_n+|u_n|^{\frac 4d}u_n+b \cdot \nabla u_n+c u_n=0,   \\
 &u_n(t_n)=S_T(t_n), \nonumber
\end{align}
where $b,c$ are given by \eqref{b} and \eqref{c} respectively.

Proposition \ref{Prop-dec-un} below contains the key geometrical decomposition of $u_n$,
which enables us to reduce the analysis of blow-up dynamics
to that of finite dimensional modulation parameters and a small remainder term.

\begin{proposition}[Geometrical decomposition]\label{Prop-dec-un}
For any $t_n$  sufficiently close to $T$,
there exist  $t^*_n<T$ and unique modulation parameters
$\mathcal{P}_n:= (\lambda_n, \alpha_n, \beta_n, \gamma_n, \theta_n)
\in C^1((t_n^*, t_n); \bbr^{2d+3})$,
such that $u_n\in C([t^*_n,t_n],H^1)$,
$u_n$ has the geometrical decomposition
\be\label{un-dec}
    u_n(t,x)
    =\lbb_n^{-\frac d2} (Q_{\calp_n} + \ve_n)(t,\frac{x-\a_n}{\lbb_n}) e^{i\theta_n}
    (=:w_n(t,x) + R_n(t,x)),
\ee
with
\be \label{Rn-wtQn}
   Q_{\calp_n}(t,y) := Q(y) e^{i\beta_n(t) \cdot y - i\frac{\g_n(t) }{4}|y|^2},\ \
   t\in [t_n^*, t_n],\ y\in \bbr^d,
\ee
and the following orthogonality conditions hold on $[t^*_n,t_n] $:
\be\ba\label{ortho-cond-Rn-wn}
&{\rm Re}\int (x-\a) w_n(t)\ol{R_n}(t)dx=0,\ \
{\rm Re} \int |x-\a|^2 w_n(t) \ol{R_n}(t)dx=0,\\
&{\rm Im}\int \Lambda w_n(t) \ol{R_n}(t)dx=0,\ \
{\rm Im}\int \nabla w_n(t) \ol{R_n}(t)dx=0,\ \
{\rm Im}\int \wt \rho_n(t) \ol{R_n}(t)dx=0,
\ea\ee
where
\be \label{rhon}
 \wt \rho_n(t,x)=  \lbb^{-\frac d2}_n(t) \rho_{\calp_n}(t,\frac{x-\a_n(t)}{\lbb_n(t)}) e^{i\theta_n(t)} ,\ \ with\
 \rho_{\calp_n}(t,y) = \rho(y)^{i\beta_n(t)\cdot y - i \frac{ \g_n(t)}{4}|y|^2},
\ee
and $\rho$ is given by \eqref{def-rho}.
\end{proposition}

The proof of Proposition \ref{Prop-dec-un} is based on Lemma \ref{Lem-imp-thm} below.
\begin{lemma}\label{Lem-imp-thm}
Let $ \calu_{r}(u_0) :=\{u\in H^1; \|u-u_0\|_{H^1}\leq r\}$,
$ \calu_{r}(\calp_0) :=\{\calp\in \bbr^{2d+3}; |\calp-\calp_0|\leq r\}$,
where $\calp_0:= (\lbb_0, \a_0, \b_0, \g_0, \t_0)$ and
$r>0$.
Assume that $u_0 \in H^1$ has the decomposition
\begin{align} \label{w0-Ur}
u_0(x)=\lambda_0^{-\frac d2}Q (\frac{x-\alpha_0}
{\lambda_0} )e^{i (\beta_0 \cdot \frac{x-\alpha_0}
{\lambda_0}-\frac{\gamma_0}{4}|\frac{x-\alpha_0}
{\lambda_0}|^2+\theta_0 )}
+R_0\
(=: w_0 + R_0),\ \ x\in \R^d,
\end{align}
$\|R_0\|_{L^2} \leq C \lbb^2_0 T$
and the orthogonality conditions in \eqref{ortho-cond-Rn-wn} hold
with $\calp_0$,
$w_0$ and $R_0$ replacing
$\calp_n$, $w_n$ and $R_n$, respectively.
Then,
there exist $T_0, r_0>0$
and a unique $C^1$ map $\Psi: \calu_{r_0}(u_0) \mapsto \calu_{r_0}(\calp_0)$
such that
for any $0<T\leq T_0$ and for any $u\in \calu_{r_0}(\calp_0)$,
\be \label{dec-u-Ur}
 u = w+ R,
\ee
where $w$ and $R$ are as in Proposition \ref{Prop-dec-un}
with $\Psi(u)$ replacing $\mathcal{P}_n$,
and the following orthogonality conditions hold:
\be\ba\label{oc2}
&{\rm Re}\int (x-\alpha) w \ol{R}dx=0,\;\;
{\rm Re}\int |x-\alpha|^2 w \ol{R}dx=0,\\
&{\rm Im}\int (\frac d2 w+(x-\alpha)\cdot \nabla w) \ol{R}dx=0,\ \
{\rm Im}\int \nabla w \ol{R}dx=0,\ \
{\rm Im}\int \wt \rho \ol{R}dx=0.
\ea\ee
\end{lemma}

The proof of Lemma \ref{Lem-imp-thm} is postponed to the Appendix for simplicity.

{\bf Proof of Proposition \ref{Prop-dec-un}.}
Let $u_n$ be the solution to \eqref{equa-un-tn},
$\mathcal{P}_0 :=(\lambda_0,\alpha_0,\beta_0,\gamma_0,\theta_0) = (T-t_n,0,0,T-t_n,\frac{1}{T-t_n})$
and $u_0 =w_0:=S_T(t_n)$.
Let $\calu_r(\calp_0)$ be as in Lemma \ref{Lem-imp-thm}.

Then, $R_0 = 0$,
and thus Lemma \ref{Lem-imp-thm} yields that
there exists $r >0$, depending on $S_T(t_n)$,
and a unique $C^1$ map $\Psi$ from $\calu_{r}(u_0)$ to
$\calu_{r}(\mathcal{P}_0)$
such that \eqref{dec-u-Ur} and \eqref{oc2} hold.
Moreover,
the $H^1$-continuity of $u_n$ implies that
for $t_n^*$ close to $t_n$,
$u_n(t) \in \calu_{r}(u_0)$ for any $t\in [t_n^*, t_n]$.
Thus,
letting $w_n(t)$ be defined as in \eqref{w0-Ur}
with $\Psi(u_n(t))$ replacing $\calp_0$
and $R_n(t):= u_n(t) - w_n(t)$,
we obtain that for any $t\in [t_n^*, t_n]$,
$\calp_n(t):= \Psi(u_n(t)) \in \calu_r(\calp_0)$
and \eqref{dec-u-Ur} and \eqref{oc2} hold
with $u_n(t)$, $w_n(t)$, $R_n(t)$ and $\calp_n(t)$
replacing $u$, $w$, $R$ and $\Psi(u)$, respectively,
which immediately yields \eqref{un-dec} and \eqref{ortho-cond-Rn-wn}.

As regards the $C^1$-regularity of modulation parameters,
in view of the $C^1$-regularity of $\Psi$ and equation \eqref{equa-un-tn},
it suffices to show that $u_n\in C([t^*_n,t_n]; H^3)$.

For this purpose,
we set $\calx(I) :=C(I; L^2) \cap L^{2+\frac 4d}(I\times \bbr^d) \cap L^2(I; H^{\frac 12}_{-1})$
and let $\calx'(I)$ be its dual space.
Then, similarly to \eqref{u-Stri-ls-bdd},
we have $\|u\|_{\calx(t^*_n,t_n)} + \|\na u\|_{\calx(t^*_n,t_n)} <\9$.
In particular,
we take a finite partition $\{I_j\}_{j=1}^L$ of $[t^*_n,t_n]$
such that
$\|u\|_{L^{2+\frac 4d}(I_j \times \bbr^d)} \leq \ve$,
where $I_j:=[t_n-j\Delta t, t_n-(j-1)\Delta t]$, $1\leq j\leq L-1$,
$I_{L}:=[t^*_n, t_{n}-(L-1) \Delta t]$
and $L<\9$.

For any multi-index $|\nu|=2$,
we obtain from \eqref{equa-un-tn} that
\begin{align*}
   i\partial_t \pa_x^\nu u + \Delta \pa_x^\nu u + (b\cdot \na + c)  \pa_x^\nu u
   + (\frac 2d +1) |u|^{\frac 4d} \pa_x^\nu u
   + \frac 2d |u|^{\frac 4d -2} u^{2} \pa_x^\nu \ol{u}
    + G(u) =0,
\end{align*}
where
$$G(u):= \sum\limits_{|\mu_1+\mu_2|=2,|\mu_2|\leq 1}(\pa_x^{\mu_1}b \cdot \na + \pa_x^{\mu_1} c) \pa_x^{\mu_2} u
         + \sum\limits_{|\mu_1+\cdots+\mu_{1+\frac 4d}|=2, |\mu_j|\leq 1} \prod\limits_{j=1}^{1+\frac 2d} \pa_x^{\mu_j} u \prod\limits_{j=1}^{\frac 2d} \pa_x^{\mu_{1+\frac 2d +j}} \ol{u}.  $$
Then, using the Strichartz and local smoothing estimates (see, e.g., \cite[Theorem 3.3]{Z18})
we obtain that
for each $0\leq j\leq L-1$,
\begin{align*}
  \|\pa_x^\nu u\|_{\calx(I_{j+1})}
  \leq& C( \|\pa_x^\nu u (t_n - j\Delta t)\|_2
       + \| |u|^{\frac 4d} \pa_x^\nu u \|_{L^{\frac{2d+4}{d+4}}(I_{j+1}\times \bbr^d)}
       + \|G(u)\|_{\calx'(I_{j+1})}),
\end{align*}
which along with the asymptotical flatness of $\pa_x^{\mu_1} b$, $\pa_x^{\mu_1} c$
and the H\"older inequality yields
\begin{align*}
   \|\pa_x^\nu u\|_{\calx(I_{j+1})}
   \leq&  C\big(\|\pa_x^\nu u (t_n - j\Delta t)\|_2
       + \|u\|^{\frac 4d}_{L^{2+\frac 4d}(I_{j+1}\times \bbr^d)} \| \pa_x^\nu u \|_{\calx(I_{j+1})} \\
       & \qquad + \|u\|_{L^2(I_{j+1}; H^{\frac 32}_{-1})}
       + \|u\|^{1+\frac 4d}_{L^{2+\frac 4d}(I_{j+1}; W^{1,2+\frac 4d})} \big)  \\
   \leq& C\big(\|\pa_x^\nu u (t_n - j\Delta t)\|_2
       + \ve^{\frac 4d}  \| \pa_x^\nu u \|_{\calx(I_{j+1})}
       +1).
\end{align*}
Hence, for $\ve$ small enough we obtain that for any $|\nu|=2$,
\begin{align} \label{pu-bdd}
   \|\pa_x^\nu u\|_{\calx(I_{j+1})}
   \leq C (\|\pa_x^\nu u (t_n - j\Delta t)\|_2 + 1), \ \ 0\leq j\leq L-1.
\end{align}

Thus,
taking into account $u(t_n) =S_T(t_n) \in H^m$ for any $m\geq 0$,
we can use \eqref{pu-bdd} and inductive arguments
to obtain that
$\|\pa_x^\nu u\|_{\calx(t_n^*,t_n)} <\9$ for any $|\nu|=2$,
which, in particular,
yields that $u\in C([t^*_n,t_n]; H^2)$.
One can also use similar arguments as above to show that
$u\in C([t^*_n,t_n]; H^3)$
and so obtain the $C^1$-regularity of modulation parameters.

Therefore, the proof is complete.
\hfill $\square$

\subsection{Modulation parameters}
This subsection contains the analysis of the dynamics of modulation parameters and the remainder.
Below
we set $g_t:= \frac{d}{dt}g$ for any $C^1$ function $g$ for simplicity.
Set
\be P_n(t):=   |\lbb_n(t)|+|\a_n(t)| + |\beta_n(t)| + |\g_n(t)|,
\ee
and denote \emph{the vector of modulation equations} by
\begin{align} \label{Mod-def}
   Mod_n(t) :=&|\lambda_n\lambda_{n,t}+\gamma_n|+|\lambda_n^2\gamma_{n,t}+\gamma_n^2|+|\lambda_n\alpha_{n,t}-2\beta_n|
                +|\lambda_n^2\beta_{n,t}+\gamma_n\beta_n| \nonumber \\
              & +|\lambda_n^2\theta_{n,t}-1-|\beta_n|^2|,
\end{align}
which can be also reformulated in the  rescaled time $\frac{ds}{dt}:= \frac{1}{\lbb_n^{2}(t)}$ as follows
\begin{align*}
   Mod_n(s)
   = &| \frac{\lambda_{n,s}}{\lbb_n}+\gamma_n|+|\gamma_{n,s}+\gamma_n^2|
                +|\frac{\alpha_{n,s}}{\lbb_n}-2\beta_n|
                +|\beta_{n,s}+\gamma_n\beta_n|
      +| \theta_{n,s}-1-|\beta_n|^2|.
\end{align*}

Note that,
since $S_T$ is the main profile for the approximating solutions near the blow-up time,
one has $Mod_n(t) \approx 0$ for $t$ close to $T$,
i.e.,
the leading order ODE system driving modulation parameters takes the form
\begin{align*}
  \lambda\lambda_t+\gamma=0,\
\lambda^2\gamma_t+\gamma^2=0,\
 \lambda\alpha_t-2\beta =0,\
  \lambda^2\beta_t+\beta\gamma=0,\
  \lambda^2\t_t-1-|\beta|^2=0.
\end{align*}

As a matter of fact,
the modulation equations arise from
the equation of remainder $\ve_n$
formulated in \eqref{equa-ve} below,
which is obtained by linearizing the equation of \eqref{equa-u-RNLS}
around the main profile $w_n$
that is close to $S_T$ near the blow-up time.
Precisely,
plugging the decomposition (\ref{un-dec}) into (\ref{equa-un-tn})
and using the identities
\begin{align} \label{Qcalp-fQcalp}
   \Delta Q_{\calp_n} - Q_{\calp_n} + f(Q_{\calp_n})
   = |\beta_n - \frac{\g_n}{2} y|^2 Q_{\calp_n}
     - i \g_n \Lambda Q_{\calp_n}
     + 2i \beta_n \cdot \na Q_{\calp_n},
\end{align}
where $f(z) = |z|^\frac 4d z$ for any $z\in \mathbb{C}$,
and
\begin{align} \label{ds-Qcalp}
   \partial_s Q_{\calp_n} = i(\beta_{n,s} \cdot y - \frac{\g_{n,s}}{4}|y|^2) Q_{\calp_n},
\end{align}
we obtain the equation of $\ve_n$ below
\begin{align} \label{equa-ve}
&i \partial_s \ve_n -M(\varepsilon_n)
-i\frac{\lambda_{n,s}}{\lambda_n}\Lambda\varepsilon_n
-i\frac{\alpha_{n,s}}{\lambda_n}\cdot
\nabla\varepsilon_n
 -(\t_{n,s}-1)\varepsilon_n  + N_f(Q_{\calp_n}, \ve_n) \nonumber  \\
=& - (\lbb_n \wt b_n \cdot \na  + \lbb_n^2 \wt c_n) (Q_{\calp_n} + \ve_n)  \nonumber  \\
 & +(\beta_{n,s}+\beta_n\gamma_n)\cdot y Q_{\calp_n}
   +i (\frac{\lambda_{n,s}}{\lambda_n}+\gamma_n )\Lambda Q_{\calp_n}
   +i (\frac{\alpha_{n,s}}{\lambda_n}-2\beta_n )\cdot
      \nabla Q_{\calp_n}    \nonumber   \\
 & +(\t_{n,s}-1-|\beta_n|^2) Q_{\calp_n}
-\frac{\gamma_n^2+\gamma_{n,s}}{4}|y|^2 Q_{\calp_n} .
\end{align}
Here,  $M(\varepsilon_n)$ denotes the linearized operator around $Q_{\calp_n}$ with respect to $\varepsilon_n$
\begin{align}
M(\varepsilon_n)
:=-\Delta \varepsilon_n+\varepsilon_n- (1+\frac 2d) |Q_{\calp_n}|^{\frac 4d}\varepsilon_n
                 -\frac 2d|Q_{\calp_n}|^{\frac 4d -2}{Q^2_{\calp_n}}\ol{\varepsilon_n},
\end{align}
the term $N_f(Q_{\calp_n}, \ve_n)$ contains the nonlinear terms of $\varepsilon_n$, i.e.,
\begin{align}
&N_f(Q_{\calp_n}, \ve_n)
:= f(Q_{\calp_n} + \ve_n) - f(Q_{\calp_n}) - f'(Q_{\calp_n}) \cdot \ve_n \nonumber \\
&=|Q_{\calp_n}+\varepsilon_n|^{\frac 4d}(Q_{\calp_n}+\varepsilon_n)
-|Q_{\calp_n}|^{\frac 4d}Q_{\calp_n}
-(1+\frac 2d)|Q_{\calp_n}|^{\frac 4d}\varepsilon_n
-\frac 2d|Q_{\calp_n}|^{\frac {4}{d}-2}{Q^2_{\calp_n}}\ol{\varepsilon_n},
\end{align}
and the coefficients of lower order perturbations
$$\wt b_n(y) := b(\lbb_n y + \a_n), \
\wt c_n(y) := c(\lbb_n y + \a_n), \ \ y\in \bbr^d. $$

Writing
$Q_{\calp_n} = \Sigma_n+i \Theta_n$,
$M(\ve_n)= M_1(\ve_n) + i M_2(\ve_n)$
and $\ve_n = \ve_{1,n} + \ve_{2,n}$
in terms of the real and imaginary parts.
We have that for some $\delta>0$,
\begin{align}
  \Sigma_n =Q + \mathcal{O}(P_n^2e^{-\delta|y|}),&\
   \Theta_n =(\beta_n \cdot y-\frac{\gamma_n }{4}|y|^2)Q(y)+ \mathcal{O}(P_n^2e^{- \delta|y|}), \label{vf-psi-Q} \\
 |\Sigma_{n,s} |+|\Theta_{n,s} |
   =& \mathcal{O} ((Mod_n + P_n^2)e^{- \delta|y|}) , \label{vfs-psis}
\end{align}
which implies that
\begin{align*}
& M_1(\ve_n)  = -\Delta \varepsilon_{n,1}+\varepsilon_{n,1}
                 -(1+\frac 2d)Q^{\frac 4d}\varepsilon_{n,1}
                 -\frac 2dQ^{\frac {4}{d}-2}\Sigma_n^2\varepsilon_{n,1}
                 -\frac 4dQ^{\frac {4}{d}-2}\Sigma_n\Theta_n\varepsilon_{n,2}
                 +\calo(P_n^2 e^{-\delta|y|}|\varepsilon_n| ), \nonumber \\
& M_2(\ve_n)  =  -\Delta \varepsilon_{n,2}+\varepsilon_{n,2}
                 -(1+\frac {2}{d})Q^{\frac 4d}\varepsilon_{n,2}
                  +\frac 2dQ^{\frac {4}{d}-2}\Sigma_n^2\varepsilon_{n,2}
                  -\frac 4dQ^{\frac {4}{d}-2}\Sigma_n\Theta_n\varepsilon_{n,1}
                  +\calo(P_n^2e^{-\delta|y|}|\varepsilon_n|).
\end{align*}
In particular, this yields that
$M_1(\ve_n)$ and $M_2(\ve_n)$ are, actually,
small deformations of the linearized operator $(L_+,L_-)$,
namely,
for some $\delta >0$,
\begin{align}
  & M_1(\ve_n) = L_{+} \ve_{n,1} - \frac 4d \Sigma_n^{\frac 4d -1} \Theta_n \ve_{n,2} + \mathcal{O}(P_n^2e^{-\delta |y|}|\ve_n|), \label{M1-L+} \\
  & M_2(\ve_n) = L_{-} \ve_{n,2} - \frac 4d \Sigma_n^{\frac 4d -1} \Theta_n \ve_{n,1} + \mathcal{O}(P_n^2e^{-\delta |y|}|\ve_n|). \label{M2-L-}
\end{align}

We also note that,
by the orthogonality conditions (\ref{ortho-cond-Rn-wn}),
\be \ba \label{ortho-cond-B}
  & {\rm Re} \int y Q_{\calp_n}(t)  \ol{\ve_n}(t) dy =0,\  \ {\rm Re}  \int |y|^2 Q_{\calp_n}  (t) \ol{\ve_n}(t) dy =0, \\
  & {\rm Im} \int \Lambda Q_{\calp_n} (t) \ol{\ve_n}(t) dy =0, \
    {\rm Im} \int \na Q_{\calp_n} (t) \ol{\ve_n}(t) dy =0, \
    {\rm Im}  \int  \rho_{\calp_n}(t) \ol{\ve_n}(t) dy =0.
\ea\ee

Proposition \ref{Prop-Mod} below
contains the key estimate for the modulation equations.

\begin{proposition} \label{Prop-Mod}
Let $t_n^*$ be close to $t_n$
such that
$\calp + \|\ve\|_{L^2}$
is sufficiently small on $[t_n^*, t_n]$.
Then,
there exists $C>0$, independent of $n$, such that
for any $t\in[t_n^*,t_n]$,
\be \label{Mod-bdd}
Mod_n(t)\leq
C(P_n^2(t)\|\varepsilon_n(t)\|_{L^2}+\|\varepsilon_n(t)\|_{L^2}^2
+ \sum\limits_{k=3}^{1+\frac 4d}\|\varepsilon_n(t)\|_{H^1}^k+P_n^6(t)).
\ee
\end{proposition}

{\bf Proof of Proposition \ref{Prop-Mod}.}
We mainly prove in detail the estimate for $|\frac{\lambda_{n,s}}{\lambda_n}+\gamma_n |$ below,
as the other four modulation equations  can be estimated similarly.
For simplicity,
the dependence on $n$ is omitted below.

In the estimates below
we shall use the estimates
\eqref{vf-psi-Q} and \eqref{vfs-psis} and the orthogonality conditions \eqref{ortho-cond-B} above.
We will also frequently use the fact that
\begin{align} \label{Q-y-integ}
  (1+|y|^n)  (Q+\na Q)  \in L^\9 \cap L^2, \ \ n\geq 1.
\end{align}

Now,
taking the inner product of (\ref{equa-ve}) with $|y|^2  {Q}_{\calp}$
and then taking the imaginary part
we obtain that,
if $\ve:= \ve_1 + i\ve_2$ in terms of real and imaginary parts,
\begin{align} \label{equa-ve1ypsi-ve2yvf}
\partial_s&\{\<\varepsilon_1,|y|^2\Sigma)
  +\<\varepsilon_2,|y|^2\Theta\>\}
  -\<\varepsilon_1,|y|^2\Sigma_s\>
-\<\varepsilon_2,|y|^2 \Theta_s\>\nonumber \\
& +\<M_1(\ve),|y|^2\Theta\>
  - \<M_2(\ve),|y|^2\Sigma\>
 -2\beta\{\<\nabla\varepsilon_1,|y|^2 \Sigma\>+\<\nabla\varepsilon_2,|y|^2 \Theta\>\} \nonumber \\
&+\gamma\{\<\Lambda\varepsilon_1,|y|^2 \Sigma\>+\<\Lambda\varepsilon_2,|y|^2 \Theta\>\}
 +|\beta|^2\{\<\varepsilon_1,|y|^2 \Theta\>-\<\varepsilon_2,|y|^2\Sigma\>\}   \nonumber\\
&+{\rm Im}\<N_f(Q_{\calp}, \ve),|y|^2 {Q}_\calp\> \nonumber \\
=& - {\rm Im}\<(\lbb \wt b \cdot \na + \lbb^2 \wt c)(\ve + Q_\calp), |y|^2 {Q}_\calp\>  \nonumber \\
 &+(\frac{\lambda_s}{\lambda}+\gamma )\{\<\Lambda\Sigma,|y|^2\Sigma)+\<\Lambda\Theta,|y|^2\Theta\>
     +\<\Lambda\varepsilon_1,|y|^2\Sigma\> + \<\Lambda\varepsilon_2,|y|^2\Theta)\} \nonumber \\
&+ (\frac{\alpha_s}{\lambda}-2\beta )\{\<\nabla\Theta,|y|^2\Theta\> + \<\nabla\Sigma,|y|^2\Sigma\>
   + \<\nabla\varepsilon_1,|y|^2\Sigma\> + \<\nabla\varepsilon_2,|y|^2\Theta\>\}  \nonumber \\
&-(\t_s-1-|\beta|^2)\{\<\ve_1,|y|^2\Theta\>-\<\ve_2,|y|^2\Sigma\>\},
\end{align}
Here,
the left-hand side above
is arranged according to the orders of $\ve_1$ and $\ve_2$,
while the right-hand side contains the lower order perturbations and
the modulation  equations.

The key point is, that
the orthogonality condition \eqref{ortho-cond-B}
enables us to gain the negligible smallness $\mathcal{O}(P^2\|\ve\|_{L^2})$
for the linear terms of $\ve_1$ and $\ve_2$,
while the nonlinear terms can be estimated easily by
using \eqref{vf-psi-Q} and \eqref{vfs-psis}.
The flatness condition \eqref{degeneracy}
will be also used  to control the lower order perturbations.
Moreover,
the coefficients of modulation equations
on the right-hand side of \eqref{equa-ve1ypsi-ve2yvf}
are all negligible,
except the one $\<\Lambda \Sigma, |y|^2\Sigma\>$
which gives the non-small coefficient $-\|yQ\|_{L^2}^2$
and thus yields the bound of $|\frac{\lbb_s}{\lbb} +\g|$.

To be precise,
using the orthogonality condition \eqref{ortho-cond-B} we have
\begin{align} \label{psv1-v2-0}
   \partial_s\{\<\varepsilon_1,|y|^2\Sigma\>+\<\varepsilon_2,|y|^2\Theta\>\}
   = \partial_s {\rm Re} \int |y|^2 Q_\calp \ol \ve dy
    =0,
\end{align}
and by \eqref{vfs-psis},
\be\label{52}
|\<\varepsilon_1,|y|^2\Sigma_s\>
+\<\varepsilon_2,|y|^2\Theta_s\>|
\leq C \left(Mod +P^2\right)\|\varepsilon\|_{L^2},
\ee

Moreover,   by \eqref{vf-psi-Q}, \eqref{M1-L+} and \eqref{M2-L-},
\begin{align*}
   &\<M_1(\ve),|y|^2\Theta\> - \<M_2(\ve),|y|^2\Sigma\> \\
  =& \< L_{+} \ve_1 + \frac 4d Q^\frac 4d \ve_1, |y|^2 \Theta\>
     -\< L_{-} \ve_2, |y|^2 \Sigma\>
     + \mathcal{O}(P^2\|\ve\|_{L^2}) \\
  =& \<L_{-}\ve_1, (\beta\cdot y-\frac{\gamma}{4}|y|^2)|y|^2Q\>
      -\<L_{-} \ve_2, |y|^2 Q\>
      + \mathcal{O}(P^2\|\ve\|_{L^2}).
\end{align*}
This yields that the second and third lines on the left-hand side of \eqref{equa-ve1ypsi-ve2yvf}
equal to
\begin{align}\label{lin1}
&|\<L_-\varepsilon_1,(\beta\cdot y-\frac{\gamma}{4}|y|^2)|y|^2Q\>
    - 2\beta\<\nabla\varepsilon_1,|y|^2Q\>
        + \gamma\<\Lambda\varepsilon_1,|y|^2Q\>
     -\<L_-\varepsilon_2,  |y|^2Q\> | \nonumber \\
& +\mathcal{O}(P^2\|\varepsilon\|_{L^2})\nonumber \\
=&|\<\varepsilon_1,(\beta\cdot y-\frac{\gamma}{4}|y|^2)L_-|y|^2Q\>
   - \<\varepsilon_2,L_-|y|^2Q \>|
  +\mathcal{O}(P^2\|\varepsilon\|_{L^2}),  \nonumber\\
=&|4\<\varepsilon_1,(\beta\cdot y-\frac{\gamma}{4}|y|^2)\Lambda Q\>
    -4\<\varepsilon_2, \Lambda Q\>|
   +\mathcal{O}(P^2\|\varepsilon\|_{L^2}),
\end{align}
where the last step is due to the identity
$L_{-} |y|^2 Q = - 4 \Lambda Q$ in  \eqref{Q-kernel}.

We claim that
\begin{align} \label{lin1.1-P2}
  4\<\varepsilon_1,(\beta\cdot y-\frac{\gamma}{4}|y|^2)\Lambda Q\>
  -4\<\varepsilon_2, \Lambda Q\>
  = \mathcal{O} (P^2 \|\ve\|_{L^2}).
\end{align}
To this end,
since $Q_{\calp} =\Sigma+i\Theta$,
using \eqref{vf-psi-Q}
we have
\begin{align*}
&{\rm Re} \Lambda Q_{\calp}
= \Lambda Q+\mathcal{O}(P^2e^{-\delta|y|}),  \\
&{\rm Im}  \Lambda Q_\calp
= (\beta\cdot y-\frac{\gamma}{4}|y|^2)\Lambda Q
   + (\beta\cdot y-\frac{\gamma}{2}|y|^2)Q
   + \mathcal{O}(P^2e^{-\delta|y|}).
\end{align*}
Taking into account \eqref{vf-psi-Q} again
we get
\begin{align*}
&4\<\varepsilon_1,(\beta\cdot y-\frac{\gamma}{4}|y|^2)\Lambda Q\>
   -4\<\varepsilon_2, \Lambda Q\> \nonumber \\
=& 4[\<\ve_1, {\rm Im} \Lambda Q_{\calp}\>
      -  \<\ve_2, {\rm Re} \Lambda Q_{\calp} \>]
   - 4 \beta \<\ve_1,  {\rm Re} y Q_{\calp} \>
    + 2 \g \<\ve_1, {\rm Re} |y|^2 Q_{\calp} \>
   + \mathcal{O}(P^2\|\varepsilon\|_{2}) \nonumber \\
=&4 {\rm Im} \int \Lambda Q_{\calp}  \ol{\varepsilon} dy
   -4\beta {\rm Re}\int   y Q_{\calp} \ol{\varepsilon} dy
   + 2\gamma {\rm Re} \int  |y|^2 Q_{\calp} \ol{\varepsilon}  dy
   + \mathcal{O}(P^2\|\varepsilon\|_{2}).
\end{align*}
Then, applying the orthogonality condition \eqref{ortho-cond-B}
we obtain \eqref{lin1.1-P2}, as claimed.

Hence,
putting together \eqref{psv1-v2-0}, \eqref{52},
\eqref{lin1} and \eqref{lin1.1-P2}
we obtain that
the linear terms  on the left-hand side
of \eqref{equa-ve1ypsi-ve2yvf} is bounded by
\begin{align} \label{lin1-bdd}
 C(Mod(t) + P^2)\|\varepsilon\|_{L^2}.
\end{align}

The nonlinear term involving $N_f(Q_\calp, \ve)$
can be bounded  directly
by using Lemma \ref{Lem-GN} and \eqref{Q-y-integ}.
Actually,
using \eqref{Nf-Taylor}, \eqref{Q-y-integ} and
Gagliardo-Nirenberg's inequality
we get
\begin{align} \label{quadra-bdd}
  |{\rm Im}\<N_f(Q_\calp, \ve),|y|^2  {Q}_\calp\>|
  \leq C \sum\limits_{k=2}^{1+\frac 4d} \int |\ve|^k dy
  \leq C ( \|\ve\|^2_{L^2}+  \sum\limits_{k=3}^{1+\frac 4d} \|\ve\|_{H^1}^k).
\end{align}

Next we treat the right-hand side of \eqref{equa-ve1ypsi-ve2yvf}.
Set $ \partial_y^\nu\wt  \phi_k  (y):= (\partial_y^\nu \phi_k) (\lbb y + \a)$.
We have
\begin{align} \label{ve-bc.0}
   & {\rm Im}\<(\lbb \wt b \cdot \na + \lbb^2 \wt c)(\ve + Q_\calp), |y|^2 {Q_\calp}\>     \nonumber  \\
  =&2\lambda \sum\limits_{k=1}^N B_k
        \{\<\nabla \wt  \phi_k\cdot(\nabla\Sigma+\nabla\varepsilon_1),|y|^2\Sigma\>
      + \<\nabla \wt \phi_k\cdot(\nabla\Theta+\nabla\varepsilon_2),|y|^2 \Theta\> \}  \nonumber  \\
&+\lambda^2 \sum\limits_{k=1}^N B_k
     \{\<\Delta \wt  \phi_k(\Sigma+\varepsilon_1),|y|^2\Sigma\>
        +  \<\Delta  \wt \phi_k(\Theta+\varepsilon_2),|y|^2\Theta\> \}  \nonumber  \\
&-\lambda^2
     \{\< \sum\limits_{j=1}^d(\sum\limits_{k=1}^N \partial_j  \wt \phi_k B_k)^2(\Theta+\varepsilon_2),|y|^2\Sigma\>
       - \<\sum\limits_{j=1}^d(\sum\limits_{k=1}^N \partial_j  \wt \phi_k B_k)^2(\Sigma+\varepsilon_1),|y|^2\Theta\>\}.
\end{align}
Using Taylor's expansion, the flatness condition (\ref{degeneracy})
and that $\pa_y^\nu \phi_k \in L^\9$ for any multi-index $\nu$,
we have
\be\label{Taylor}
|(\partial_y^\nu \wt  \phi_k) (y)|
\leq C(\lambda y+\alpha)^{6-|\nu|}
\leq C P^{6-|\nu|} (1+|y|^6), \ \ 0\leq |\nu|\leq 5.
\ee
Then,
using the integration by parts formula and \eqref{Q-y-integ} again we have
\begin{align}
  |2\lambda B_k\<\nabla \wt \phi_k\cdot(\nabla \Sigma +\nabla\varepsilon_1),|y|^2 \Sigma\>|
  \leq C P^6 (1+\|\ve\|_{L^2}).
\end{align}
We can estimate the other terms on the right-hand side of \eqref{ve-bc.0} similarly
and so obtain
\begin{align} \label{inhom-bdd}
    |{\rm Im}\<(\lbb \wt b \cdot \na + \lbb^2 \wt c)(\ve + Q_\calp), |y|^2  {Q_\calp}\>|
    \leq  C P^6(1+\|\ve\|_{L^2}).
\end{align}

It remains to treat the coefficients of  modulation equations
on the right-hand side of \eqref{equa-ve1ypsi-ve2yvf}.
Straightforward computations yield that
\begin{align}
   \<\Lambda\Sigma,|y|^2\Sigma\>
   = - \int |y|^2 |\Sigma|^2 dy
   = - \|yQ\|^2_{L^2} + \mathcal{O}(P^2),
\end{align}
where in the last step we used \eqref{vf-psi-Q}.
Moreover,
since $\na Q$ is odd while $Q$ is even,
we have $\<\na Q, |y|^2 Q\> =0$,
which along with \eqref{vf-psi-Q} implies that
\begin{align*}
    \<\na \Sigma, |y|^2 \Sigma \>
    = (\na Q, |y|^2 Q) + \mathcal{O}(P^2)
    = \mathcal{O}(P^2).
\end{align*}
The remaining coefficients of $(\frac{\lbb_s}{\lbb} + \g)$ can be
easily bounded
by $C(\|\ve\|_{L^2} + P^2)$,
by  using \eqref{vf-psi-Q} and \eqref{Q-y-integ}.
Thus,  the right-hand side of \eqref{equa-ve1ypsi-ve2yvf}
is bounded by
\begin{align} \label{RHS-bdd}
   -\|yQ\|_{L^2}^2 (\frac{\lambda_s}{\lambda}+\gamma)
    + CMod(\|\ve\|_{L^2} + P^2)
    +C P^6 (1+\|\ve\|_{L^2}).
\end{align}

Therefore,
plugging
\eqref{lin1-bdd}, \eqref{quadra-bdd} and \eqref{RHS-bdd}
into equation \eqref{equa-ve1ypsi-ve2yvf}
we obtain
\begin{align} \label{esti-lbbslbb-gamma}
\|yQ\|_{L^2}^2|\frac{\lambda_s}{\lambda}+\gamma|
\leq& C\big( (P +\|\varepsilon\|_{L^2}) Mod+P^2\|\varepsilon\|_{L^2}+\|\varepsilon\|_{L^2}^2 \nonumber \\
     & \qquad + \sum\limits_{k=3}^{1+\frac 4d}\|\varepsilon\|_{H^1}^k + P^6(1+\|\ve\|_{L^2}) \big).
\end{align}

Similar arguments apply also to the other four modulation equations.
Actually, taking the inner products of   equation (\ref{equa-ve})
with $y {Q_{\calp} }$, $i\Lambda {Q_{\calp} }$, $i\nabla {Q_{\calp} }$, $i{\rho_{\calp}}$,
respectively,
then taking the imaginary parts
and using analogous arguments as above,
we can obtain the same bounds for
$\frac d2 \|Q\|^2_{L^2}|\frac{\a_s}{\lbb}-2\beta|$,
$\frac 14 \|yQ\|_{L^2}^2|\g_s+\g^2|$,
$\frac d2 \|Q\|_{L^2}^2 |\beta_s + \beta \g|$
and $ \frac 12 \|yQ\|_{L^2}^2 |\theta_s - 1-|\beta|^2|$,
respectively.

Therefore,
taking $P + \|\ve\|_{L^2}$ small enough
such that $C(P+ \|\ve\|_{L^2}) \leq \frac 12$
we obtain \eqref{Mod-bdd}
and finish the proof.
\hfill $\square$

\section{Uniform estimates of approximating solutions}   \label{Sec-Unif-Esti}

\subsection{Main results}   \label{Subsec-main-unif-esti}
This section is devoted to the key
uniform estimates of approximating solutions.
The main result is formulated in Theorem \ref{Thm-dec-un-t0} below.

\begin{theorem}[Uniform estimates]\label{Thm-dec-un-t0}
Let $0<\delta< \frac{1}{12}$.
Then,
there exists small $\tau^*>0$,
such that  for  any $T\in (0,\tau^*]$
and  for $n$ large enough,
$u_n$ admits the unique
geometrical decomposition $u_n=\omega_n+R_n$ on $[0, t_n]$ as in \eqref{un-dec},
and the following estimates hold:

$(i)$ For the reminder term,
\be \label{wn-Tt}
\|\nabla R_n(t)\|_2\leq (T-t)^2,\quad\|R_n(t)\|_2\leq (T-t)^3.
\ee

$(ii)$ For the modulation parameters,
\begin{align}
&\lf|\la_n(t) - (T-t) \rt| + \ \lf|\gamma_n(t)  - (T-t) \rt|\leq (T-t)^{3+\delta},   \label{lbbn-Tt} \\
&|\al_n(t)|+|\beta_n(t)|\leq (T-t)^{2+\delta},  \label{anbn-Tt}\\
&|\theta_n(t) - \frac{1}{T- t}|\leq (T-t)^{1+\delta}. \label{thetan-Tt}
\end{align}
\end{theorem}

The proof of
Theorem \ref{Thm-dec-un-t0} is based on continuity arguments
and Proposition \ref{Prop-dec-un-t0-boot} below.

\begin{proposition}[Bootstrap] \label{Prop-dec-un-t0-boot}
Fix $0<\delta< \frac{1}{12}$ and $n$ large enough.
Suppose that there exists $t^*\in(0,t_n)$ such that
$u_n$ admits the unique geometrical decomposition \eqref{un-dec} on $[t^*,t_n]$
and the estimates \eqref{wn-Tt}-\eqref{thetan-Tt} hold.

Then,
there exists $t_*\in [0, t^*)$ such that \eqref{un-dec} holds
on $[t_*, t_n]$
and the  bounds in estimates \eqref{wn-Tt}-\eqref{thetan-Tt} can be refined to $1/2$,
i.e., for any $t\in [t_*, t_n]$,
\begin{align}
&\|\nabla R_n(t)\|_{L^2}\leq \frac{1}{2}(T-t)^2,\quad\|R_n(t)\|_{L^2}\leq \frac{1}{2}(T-t)^3,  \label{wn-Tt-boot-2} \\
&\lf|\la_n(t) - (T-t) \rt|+ \lf|\gamma_n(t) - (T-t) \rt|\leq \frac{1}{2}(T-t)^{3+\delta}, \label{lbbn-Tt-boot-2} \\
&|\al_n(t)|+ |\beta_n(t)|\leq \frac{1}{2}(T-t)^{2+\delta},  \label{anbn-Tt-boot-2} \\
&  |\theta_n(t) - \frac{1}{T-t} |\leq \frac{1}{2}(T-t)^{1+\delta}. \label{thetan-Tt-boot-2}
\end{align}
\end{proposition}

In order to prove Proposition \ref{Prop-dec-un-t0-boot},
we apply Lemma \ref{Lem-imp-thm}
with $u_n(t^*)$ and $\calp_n(t^*)$ replacing $u_0$
and $\calp_0$
to obtain that there exists $t_*\in [0,t^*)$  such that
$u_n$ has the geometrical decomposition
$u_n=w_n + R_n$ on $[t_*, t_n]$ as in \eqref{un-dec}.

Moreover, using the local well-posedness theory and the $C^1$-regularity of modulation parameters,
we can take $t_*$ sufficiently close to $t^*$ such that for any $t\in [t_*,t_n]$,
\begin{align}
&\|\nabla R_n(t)\|_{L^2}\leq 2(T-t)^2, \ \ \|R_n(t)\|_{L^2}\leq 2(T-t)^3, \label{wn-Tt2} \\
&\lf|\la_n(t)-(T-t)\rt| + \lf|\gamma_n(t) - (T-t)\rt|\leq 2(T-t)^{3+\delta}, \label{lbbn-Tt3delta}   \\
&|\al_n(t)| + |\beta_n(t)|\leq 2(T-t)^{2+\delta}, \label{anbn-Tt2delta}  \\
&   |\theta_n - \frac{1}{T-t} |\leq 2(T-t)^{1+\delta}. \label{thetan-Tt1delta}
\end{align}

In particular, estimate \eqref{lbbn-Tt3delta} yields that
if $T^{2+\delta} \leq \frac 14$,
\be \label{lbbn-Tt*}
   \frac 12 (T-t)\leq\la_n(t), \g_n(t) \leq 3(T-t),\ \ \forall t\in [t_n^*, t_n].
\ee

By virtue of estimates \eqref{wn-Tt2}-\eqref{lbbn-Tt*},
we have the refined estimates below.

\begin{lemma}   \label{pre-est}
Assume   estimates \eqref{wn-Tt2}-\eqref{thetan-Tt1delta} to hold.
Then, there exists $C>0$ independent of $n$
such that  for all $t\in [t_*, t_n]$,
\begin{align}
   & P_n(t)\leq C\lambda_n(t), \label{Pn-ven-lbbn} \\
    &\|\varepsilon_n(t)\|_{L^2}=\|R_n(t)\|_{L^2} \leq C \lbb_n^3(t), \ \
     {\lambda^{-1}_n}\|\na \ve_n(t)\|_{L^2}=\|\na R_n(t)\|_{L^2} \leq C \lbb_n^2(t), \label{wn-lbbn}  \\
   &\|w_n\|_{L^2}= \|Q\|_{L^2}, \ \ \|\nabla w_n\|_{L^2} + \|\frac{x-\a_n}{\lbb_n} \cdot \na w_n\|_{L^2} \leq  C \lbb_n^{-1}. \label{R-naR}
\end{align}
\end{lemma}

We also note that the estimate of modulation equations in \eqref{Mod-bdd} still holds on $[t_*,t_n]$.
Then,
using Lemma \ref{pre-est}
we have the refined estimate of modulation equations below.

\begin{lemma} \label{Lem-Mod-w-lbb}
Assume  estimates \eqref{wn-Tt2}-\eqref{thetan-Tt1delta} to hold
with $T$ sufficiently small.
Then,
we have that for some $C>0$ independent of $n$,
\be \label{Mod-w-lbb}
Mod_n(t) \leq C (\lbb_n(t))^5, \ \ \forall t\in[t_*,t_n].
\ee
\end{lemma}

Lemma \ref{Lem-eeta} below contains the estimates of $\eta_n$ given by \eqref{etan-Rn},
which follow straightforwardly from
the explicit expression of $\eta_n$ and the estimate \eqref{Taylor}.

\begin{lemma}  \label{Lem-eeta}
Assume estimates \eqref{wn-Tt2}-\eqref{thetan-Tt1delta} to hold with $T$ sufficiently small.
There exists a constant $C>0$ independent of $n$,
such that for all $t\in[t_*,t_n]$,
\begin{align}
\|\partial^\nu \eta_n(t)\|_{L^2} \leq  \frac{C}{\lbb_n^{2+|\nu|}(t)} (Mod_n(t) + (\lbb_n(t))^6)
                               \leq C \lbb^{3-|\nu|}_n(t),  \ \
     \forall\ 0\leq |\nu|\leq 2.   \label{eta-L2}
\end{align}
\end{lemma}

Below Subsections \ref{Subsec-Energy} and \ref{Subsec-Mono}
contain the key analysis of the energy \eqref{energy}
and the monotonicity formula of the Lyapounov functional \eqref{def-I}, respectively.
We shall assume that  estimates  \eqref{wn-Tt2}-\eqref{thetan-Tt1delta}
hold on $[t_*, t_n]$
with $T$ sufficiently small throughout Subsections \ref{Subsec-Energy} and \ref{Subsec-Mono}.

\subsection{Energy estimate} \label{Subsec-Energy}

Because of the presence of  lower order perturbations,
the energy \eqref{energy} of solutions to equation \eqref{equa-u-RNLS}
is no longer conserved.
Here, under Assumption $(A1)$,
we are able to control the variation of energy
up to a suitable order,
which is a key ingredient to
obtain the refined estimate \eqref{anbn-Tt-boot-2}
for the modulation parameters $\a_n$ and $\beta_n$.

\begin{theorem}  \label{Thm-energy}
There exists $C>0$, independent of $n$,
such that  for any $t\in[t_*,t_n]$,
\be \label{esti-Eut-Eutn}
      |E(u_n)(t)-E(u_n)(t_n)| \leq C(T-t)^{3}.
\ee
\end{theorem}

{\bf Proof.}
Using (\ref{equa-un-tn}) and the integration by parts formula we have
\begin{align} \label{dtE}
\frac{d}{dt}E(u_n)
=&-{\rm Im}\int(\ol{b}\cdot\nabla \ol{u_n}+ \ol{c} \ol{u_n})(\Delta u_n+|u_n|^\frac 4d u_n) dx \nonumber \\
=&-2\sum\limits_{k=1}^N B_k {\rm Re}\int \nabla^2 \phi_k(\nabla u_n,\nabla \ol{u_n})dx
+\frac{1}{2}\sum\limits_{k=1}^N B_k \int \Delta^2 \phi_k|u_n|^2dx  \nonumber \\
& +\frac{2}{d+2}\sum\limits_{k=1}^N  B_k \int\Delta \phi_k|u_n|^{2+\frac{4}{d}}dx
 -\sum\limits_{j=1}^d {\rm Im}\int \na(\sum\limits_{k=1}^N \partial_j \phi_k B_k)^2 \cdot \na u_n \ol{u_n} dx.
\end{align}
Then, by \eqref{un-dec},
\begin{align} \label{e}
 |\frac{d}{dt}E(u_n)|
 \leq&\frac{C}{\lambda_n^2} \sum\limits_{k=1}^N
        \int (|\na^2 \phi_k| + |\Delta \phi_k| + |\partial_j\phi_k \na \partial_j \phi_k|)
               (\lambda_n y+\alpha_n)(|\nabla Q_\calp|^2+Q^{2+\frac{4}{d}})(y) dy  \nonumber \\
   &+C \sum\limits_{k=1}^N \sum\limits_{j=1}^d \int (|\Delta^2 \phi_k| + |\partial_j \phi_k \na \partial_j \phi_k|)
              (\lambda_n y+\alpha_n) Q^2(y) dy \nonumber \\
 & +C(\|R\|_{H^1}^2 + \|R\|_{L^{2+\frac 4d}}^{2+\frac{4}{d}})
=: K_1+K_2 + K_3.
\end{align}
Using  \eqref{Q-y-integ} and (\ref{Taylor})
we have
\be\label{e1}
    K_1 + K_2 \leq C \lambda_n^2.
\ee
Moreover,
applying the Gagliardo-Nirenberg inequality (\ref{G-N})
and using \eqref{wn-lbbn} we get
\be\label{e2}
K_3  \leq C(\|R\|_{H^1}^2 + \|R\|_{H^1}^{2+\frac 4d})
     \leq C\lambda_n^4.
\ee
Thus, plugging (\ref{e1}) and (\ref{e2}) into (\ref{e}) and using \eqref{lbbn-Tt*}
we obtain
\be
 |\frac{d}{dt}E(u_n)| \leq C(T-t)^{2},
\ee
which yields \eqref{esti-Eut-Eutn} and finishes the proof.
\hfill $\square$

As a consequence of Theorem \ref{Thm-energy}
and the coercivity of operator $L$,
we have
\begin{lemma} \label{Lem-Pn-ven-lbbn}
There exists $C>0$ independent of $n$ such that
for all  $t\in[t_*,t_n]$,
\be \label{b-g-ve-lbb}
|\beta_n(t)|^2\leq C|\lambda_n^2(t)-\g_n^2(t)|+ C(T-t)^{5}.
\ee
\end{lemma}

{\bf Proof.}
Let $F(z):= \frac{d}{2d+4}|z|^{2+\frac 4d}$, $z\in \mathbb{C}$.
We suppress the $n$ dependence below.

First note that, since
$ M(u)(t)=M(u)(t_n)=\|Q\|_{L^2}^2$,
\begin{align}  \label{mcl}
2 {\rm Re} \int Q_{\calp}  \bar{\varepsilon}dy + \|\varepsilon(t)\|_{L^2}^2 = 0.
\end{align}
The identity
\eqref{mcl} allows to gain an extra factor $\|\ve\|_{L^2}$
when estimating ${\rm Re} \int Q_{\calp}  \bar{\varepsilon}dy$.

Now,
rewriting $E(u)$ in the renormalized variables
and using \eqref{mcl}
we have
\begin{align}  \label{Eu.0}
E(u)
=&\frac{1}{2\lambda^2}\int|\nabla Q_{\calp} +\nabla\varepsilon|^2 dy
   - \frac{1}{\lbb^2} \int F(Q_\calp + \ve) dy \nonumber \\
 & +\frac{1}{\lambda^2}{\rm Re}\int Q_{\calp} \ol{\varepsilon}dy+\frac{1}{2\lambda^2}\int|\varepsilon|^2dy.
\end{align}
Then, using the identities $E(Q)=0$,
$|\na Q_\calp|^2
   = |\na Q|^2 + |\beta - \frac{\g}{2} y|^2 |Q|^2$
and the  expansion \eqref{F-Taylor*}
and then rearranging each term according to the order of $\ve$
we come to
\begin{align}  \label{Eu.1}
E(u)=&\frac{1}{2\lambda^2}
       \int |\beta - \frac \g 2 y|^2 Q^2dy
-\frac{1}{\lambda^2}{\rm Re}\int(\Delta Q_{\calp}
- Q_{\calp}+|Q_{\calp}|^{\frac 4d}Q_{\calp})\bar{\varepsilon}dy  \nonumber \\
&+\frac{1}{2\lambda^2}{\rm Re}\int|\nabla\varepsilon|^2+|\varepsilon|^2
  - (1+\frac {2}{d})|Q_{\calp}|^{\frac 4d}|\varepsilon|^{2}
-\frac{2}{d}|Q_{\calp}|^{\frac {4}{d}-2}Q_{\calp}^2\bar{\varepsilon}^2dy  \nonumber \\
&   +{\calo (  \frac{1}{\lbb^{2}}  \sum\limits_{k=3}^{2+\frac 4d}\|\varepsilon\|_{H^1}^k)}.
\end{align}
In order to estimate the linear term of $\ve$ in \eqref{Eu.1} above,
using \eqref{Qcalp-fQcalp} and \eqref{ortho-cond-B} we get
\begin{align} \label{esti-E12}
    & {\rm Re}\int(\Delta Q_{\calp}
- Q_{\calp}+|Q_{\calp}|^{\frac 4d}Q_{\calp}) \ol{\varepsilon}dy \nonumber \\
  = &  {\rm Im}\int(\gamma\Lambda Q_{\calp} -2\beta\cdot
\nabla Q_{\calp}) \ol{\varepsilon}dy
    + {\rm Re} \int |\beta - \frac \g 2 y|^2 Q_{\calp} \ol{\varepsilon}dy
   = \mathcal{O}(P^2\|\ve\|_{L^2}).
\end{align}

Moreover,
the quadratic interactions of $\ve$
can be bounded from below by using the coercivity of operator $L$.
Precisely,
using \eqref{vf-psi-Q}
we get that if $\ve := \ve_1 + i \ve_2$,
\begin{align} \label{ve-quadra-Eu}
  & {\rm Re}\int|\nabla\varepsilon|^2+|\varepsilon|^2-(1+\frac {2}{d})|Q_{\calp}|^{\frac 4d}|\varepsilon|^{2}
-\frac{2}{d}|Q_{\calp}|^{\frac {4}{d}-2}Q_{\calp}^2\bar{\varepsilon}^2dy  \nonumber \\
   =& \<L_{+}\ve_1, \ve_1\> + \<L_{-}\ve_2, \ve_2\> + \mathcal{O}(P\|\ve\|^2_{L^2})
   \geq  \nu \|\ve\|_{H^1}^2 + \mathcal{O}(P\|\ve\|^2_{L^2}) ,
\end{align}
where the last step is due to Corollary \ref{Cor-coer-f-H1}
and the estimates
\begin{align} \label{Q-ve1-ve2-inner}
 \< Q,\varepsilon_1 \>
= -\frac 12 \|\varepsilon\|_{L^2}^2+ \calo(P) \|\varepsilon\|_{L^2},
\<yQ,\varepsilon_1\>
=&  \calo(P)\|\varepsilon\|_{L^2},
\<|y|^2Q,\varepsilon_1\>
=  \calo(P)\|\varepsilon\|_{L^2}, \nonumber   \\
 \<\nabla Q,\varepsilon_2\> =  \calo(P)\|\varepsilon\|_{L^2},
\<\Lambda Q,\varepsilon_2\>=&  \calo(P)\|\varepsilon\|_{L^2},
\<\rho,\varepsilon_2\> =  \calo(P) \|\varepsilon\|_{L^2},
\end{align}
which follow  from \eqref{vf-psi-Q}, \eqref{ortho-cond-B} and (\ref{mcl}).

Hence, plugging \eqref{esti-E12} and \eqref{ve-quadra-Eu} into \eqref{Eu.0}
we obtain that,
for  $\|\ve\|_{L^2}$ and $P$ small enough
(or equivalently, $t_*$ close to $t_n$),
\begin{align*}
\lambda^2 E(u)
 \geq  \frac{1}{2}\|Q\|_{L^2}^2|\beta|^2+\frac{1}{8}\|yQ\|_{L^2}^2\gamma^2
+ \frac \nu 2 \|\varepsilon\|_{H^1}^2
  - C(P^2\|\ve\|_{L^2}+P\|\ve\|^2_{L^2} + \sum\limits_{k=3}^{2+\frac 4d} \|\ve\|^k_{H^1}).
\end{align*}

Therefore,
taking into account
$E(u(t_n)) =E(S_T(t_n))=\frac{1}{8}\|yQ\|_{L^2}^2$
we arrive at
\begin{align} \label{esti-Eu*}
   \frac{1}{2}\|Q\|_{2}^2|\beta(t)|^2
+ \frac \nu 2 \|\varepsilon(t)\|_{H^1}^2
\leq& \frac{1}{8} \|yQ\|_{L^2}^2 |\lambda^2(t) - \g^2(t)|
+\lambda^2(t) |E(u)(t)-E(u)(t_n)|   \nonumber \\
& + C(P^2\|\ve\|_{L^2}+P\|\ve\|^2_{L^2} + \sum\limits_{k=3}^{2+\frac 4d} \|\ve\|^k_{H^1}),
\end{align}
which yields \eqref{b-g-ve-lbb},
due to \eqref{lbbn-Tt*}, \eqref{wn-lbbn} and \eqref{esti-Eut-Eutn}.
The proof is complete.
\hfill $\square$

\subsection{Monotonicity of generalized energy}  \label{Subsec-Mono}

Let $\chi(x)=\psi(|x|)$ be a smooth radial function on $\R^d$,
where $\psi$ satisfies
$\psi'(r) = r$ if $r\leq 1$,
$\psi'(r) = 2- e^{-r}$ if $r\geq2$,
and
\be\label{chi}
\lf|\frac{\psi^{'''}(r)}{\psi^{''}(r)}\rt|\leq C,
\ \ \frac{\psi'(r)}{r}-\psi^{''}(r) \geq0.
\ee
Let $\chi_A(y) :=A^2\chi(\frac{y}{A})$, $A>0$, $y \in \bbr^d$,
and $F(u)$, $f(u)$, $f'(v)\cdot R$ and $N_f(v,R)$
be  as in Subsection \ref{Subsec-Expan}.
As in \cite{R-S},
we define the generalized energy by
\begin{align} \label{def-I}
I_n(t) := &\frac{1}{2}\int |\nabla R_n|^2+\frac{|R_n|^2}{\lambda_n^2}dx
           -{\rm Re}\int F(u_n)-F(w_n)-f(w_n)\ol{R_n}dx \nonumber \\
&+\frac{\gamma_n}{2\lambda_n}{\rm Im} \int (\nabla\chi_A) (\frac{x-\alpha_n}{\lambda_n})\cdot\nabla R_n\ol{R_n}dx.
\end{align}

The main result of this subsection is formulated below.
\begin{theorem}  \label{Thm-I-mono}
For $n$ large enough and for any $t\in[t_*,t_n]$, we have
\be \label{dIt-mono}
\frac{d I_n}{dt}
\geq \frac{C_1}{\lambda_n} \int (|\nabla R_n|^2+\frac{|R_n|^2}{\lambda_n^2}) e^{-\frac{|x-\alpha_n|}{A\lambda_n}}dx
    -C_2(A)\lambda_n^4,
\ee
where $C_1, C_2(A)>0$ are independent of $n$.
\end{theorem}

\begin{remark}
Theorem \ref{Thm-I-mono} enables us to obtain
the order $\mathcal{O}(\lbb_n^5)$ of the generalized energy,
which in turn gives us the important refined estimate \eqref{wn-Tt-boot-2}
for the remainder $R_n$.
See Subsection \ref{Subsec-proof-bootstrap} below.
\end{remark}

In order to prove Theorem \ref{Thm-I-mono},
we write $I_n$ into two parts $I_n= I_n^{(1)}+ I_n^{(2)}$,
where
\begin{align}
I_n^{(1)}&:=\frac{1}{2}\int |\nabla R_n|^2+\frac{|R_n|^2}{\lambda_n^2}dx
           -{\rm Re}\int F(u_n)-F(w_n)-f(w_n)\ol{R_n}dx, \label{I1}\\
I_n^{(2)}&:=\frac{\gamma_n}{2\lambda_n}{\rm Im} \int \nabla\chi_A (\frac{x-\alpha_n}{\lambda_n} )\cdot\nabla R_n \ol{R}_n dx.   \label{I2}
\end{align}
The evolution formulas of  $I^{(i)}_{n,t}:= \frac{d}{dt} I_n^{(i)}$, $i=1,2$,
are contained in Lemmas \ref{Lem-I1t} and  \ref{Lem-I2t} below, respectively.
See also \cite{R-S} for similar computations
for the inhomogeneous nonlinear Schr\"odinger equations in dimension two.

In the sequel,
we omit the dependence on $n$  for simplicity.

\begin{lemma}   \label{Lem-I1t}
For every $t\in[t_*,t_n]$, we have
\begin{align} \label{I1t}
I^{(1)}_{t}
   =&-\frac{\gamma}{\lambda^2} {\rm Re}\int ((1+\frac 2d)|w|^\frac 4d|R|^2
      +\frac 12 (1+\frac 2d) |w|^{\frac 4d -2} \ol{w}^2R^2
      +\frac 12 (\frac 2d -1) |w|^{\frac 4d -2} w^2 \overline{R}^2)dx\nonumber\\
&-\frac{\gamma}{\lambda} {\rm Re}
      \int(\frac{x-\alpha}{\lambda})\cdot\nabla \ol{w}
      \bigg\{ \frac 2d (1+\frac 2d) |w|^{\frac {4}{d}-2}w|R|^2
        + \frac 1d (1+ \frac 2d)|w|^{\frac {4}{d}-2}\ol{w}R^2  \nonumber \\
 & \qquad \qquad \qquad \qquad \qquad \ \   + \frac 1d (\frac 2d -1) |w|^{\frac 4d -4} w^3 \overline{R}^2 \bigg\}dx
   +\frac{\gamma}{\lambda^4} \|R\|_{L^2}^2
  +\mathcal{O}( \lambda^4).
\end{align}
\end{lemma}

{\bf Proof.}
Since
$ \partial_t F(u)
= {\rm Re} \lf(f(u) \partial_t \ol{u}\rt)$,
$\partial_t f(u)
   = (1+\frac 2d) |u|^\frac 4d \partial_t u
     + \frac 2d |u|^{\frac 4d -2} u^2 \partial_t\ol{u}$,
we have
\begin{align}
  I^{(1)}_t
  =& {\rm Re} \<-\Delta R + \frac{1}{\lbb^2} R - f(u) + f(w), \partial_t R\>
    - {\rm Re} \< N_f(w,R), \partial_t w\>
    - \frac {\lbb_t}{\lbb^3} \|R\|^2_{L^2}.
\end{align}
Then,
using equation \eqref{equa-R}
we obtain
\begin{align} \label{equa-I1t}
I^{(1)}_t
=&-{\rm Im}\int (\Delta R-\frac{1}{\lambda^2}R + f'(w)\cdot R) \ol{\eta} dx
    -{\rm Re}\int N_f(w,R) \pa_t \ol{w} dx \nonumber\\
  &    -\frac{1}{\lambda^2}{\rm Im}\int \frac 2d|w|^{\frac{4}{d}-2}w^2 \ol{R}^2dx
      -{\rm Im}\int  N_f(w,R) \ol{\eta} dx
   -\frac{1}{\lambda^2}{\rm Im}\int N_f(w,R)  \ol{R}dx  \nonumber \\
   & -\frac{\lbb_t}{\lambda^3} \|R\|_{L^2}^2
    -{\rm Im} \int
(\Delta R-\frac{1}{\lambda^2} R +f(u)-f(w)) (\ol{b}\cdot \nabla \ol{R}+\ol{c}\ol{R}) dx  \nonumber \\
=:&  \sum\limits_{j=1}^7 I^{(1)}_{t,j}.
\end{align}
Below we estimate each term on the right-hand side above separately.

$(i)$ {\it Estimate of the first order term $I^{(1)}_{t,1}$.}
Using H\"older's inequality we have
\begin{align*}
     |I^{(1)}_{t,1}|
    \leq  \|\na \eta\|_{L^2} \|\na R\|_{L^2}
        + \frac{1}{\lbb^2} \|\eta\|_{L^2} \|R\|_{L^2}
        + (1+\frac 4d)\int |\eta| |w|^{\frac 4d} |R| dx.
\end{align*}
Note that,
by the pointwise bound $|w|\leq C \lbb^{-\frac d2}$
and H\"older's inequality,
\begin{align}
   \int |\eta| |w|^{\frac 4d} |R| dx
   \leq \frac{1}{\lbb^2}C  \|\eta\|_{L^2} \|R\|_{L^2}.
\end{align}
Thus,
taking into account  \eqref{wn-lbbn} and \eqref{eta-L2} we obtain
\begin{align} \label{esti-I1t-i}
    |I^{(1)}_{t,1}| \leq C( \lbb^2 \|\na R\|_{L^2}
                     + \lbb \|R\|_{L^2} )\leq C\lbb^4.
\end{align}

(ii) {\it Estimates of the nonlinear terms $I^{(1)}_{t,j}$, $2\leq j\leq 6$.}
Let $N_{f,2}(w,R)$ be defined as in \eqref{f-quadratic}
with $w$ replacing $v$
and write
\begin{align} \label{I1_2-I1_21-I1_22}
  I^{(1)}_{t,2} &= -{\rm Re} \int  N_{f,2}(w,R) \partial_t \ol{w} dx
              -{\rm Re} \int(N_f(w,R) - N_{f,2}(w,R) )   \partial_t \ol{w} dx \nonumber \\
          & =:   I^{(1)}_{t,21} +  I^{(1)}_{t,22}.
\end{align}

Letting $y=\frac{x-\a}{\lbb}$
and using the identities
\begin{align*}
   \partial_t Q_\calp (t,y)
   =  (i\beta_t \cdot y - i\frac{\g_t}{4} |y|^2) Q_\calp(t,y), \ \
    \na_x w (t,x)
   = (\lbb(t))^{-\frac d2 -1} (\na_y Q_\calp)(t,y) e^{i\theta(t)},
\end{align*}
we have
\begin{align}
  \partial_t w
  = - \frac d2 \frac{\lbb_t}{\lbb} w
    + i \beta_t \cdot y w
    - i \frac{\g_t}{4} |y|^2 w
    - \lbb_t y\cdot \na_x w
    - \a_t \cdot \na_x w
    + i \theta_t w,
\end{align}
which along with (\ref{lbbn-Tt3delta})-(\ref{thetan-Tt1delta}) and Lemma \ref{Lem-Mod-w-lbb}
implies that
\begin{align} \label{pt-R}
 \partial_t w =\frac{i}{\lambda^2}w +\frac{d}{2}\frac{\gamma}{\lambda^2}w
       +\frac{\gamma}{\lambda} y \cdot \nabla w
    - 2 \frac{\beta}{\lbb} \na_x w
    +\mathcal{O}((1+|y|^2)(|w|+|\na w|)).
\end{align}
This yields that for the quadratic terms of $R$,
\begin{align} \label{esti-lin-omega}
   I^{(1)}_{t,21} +  I^{(1)}_{t,3}
   =&- \frac{\gamma}{\lambda^2} {\rm Re}
      \int (1+\frac 2d)|w|^{\frac {4}{d}}|R|^2
     + \frac 12 (1+\frac 2d) |w|^{\frac {4}{d}-2}\ol{w}^2R^2
     + \frac 12 (\frac 2d -1)|w|^{\frac {4}{d}-2}w^2\overline{R}^2dx \nonumber \\
 &- \frac{\gamma}{\lambda} {\rm Re}
    \int (\frac{x-\alpha}{\lambda})\cdot\nabla \ol{w}
     \bigg\{\frac 2d(1+\frac 2d) |w|^{\frac {4}{d}-2}w|R|^2
      + \frac 1 d (1+\frac 2d) |w|^{\frac {4}{d}-2}\ol{w}R^2  \nonumber\\
  & \qquad \qquad  \qquad \qquad \qquad \ \   + \frac 1d (\frac 2d-1) |w|^{\frac 4d-4} w^3 \overline{R}^2\bigg\}dx
      +\calo (\frac{1}{\lambda^2}\|R\|_{L^2}^{2})  .
\end{align}

As regards the remaining higher order terms contained in $ I^{(1)}_{t,22}$,
using the pointwise bounds $|w|\leq C \lbb^{-\frac d2}$,
$|\na w|\leq C \lbb^{-\frac d2 -1}$,
the expansion \eqref{f-Taylor} and equation \eqref{pt-R}
we have
\begin{align}
 | I^{(1)}_{t,22}|
 \leq &\frac{C}{\lbb^2} \sum\limits_{k=3}^{1+\frac 4d} \int |w|^{2+\frac 4d-k}|R|^{k}dx
        +C \sum\limits_{k=3}^{1+\frac 4d} \int (1+ |\frac{x-\alpha}{\lbb}|) |\nabla w||w|^{1+\frac{4}{d}-k}|R|^{k}dx   \nonumber  \\
      & + C\sum\limits_{k=3}^{1+\frac 4d}  \int |w|^{\frac 4d +2-k} |R|^k dx
 \leq  C\sum\limits_{k=1}^{1+\frac 4d}  \lbb^{\frac {d}{2}k-4-d}\|R\|^k_{L^k}.
\end{align}
Note that,
Gagliardo-Nirenberg's inequality and  \eqref{wn-lbbn} yield  that
\begin{align} \label{R-Lk-k}
   \|R\|_{L^k}^k
   \leq C \|R\|_{L^2}^{k+d-\frac d2 k}
          \|\na R\|_{L^2}^{\frac d2 k -d}
   \leq C \lbb^{d+3k -\frac d2 k}.
\end{align}
Hence, we obtain
\begin{align}  \label{esti-I1t-iii.1}
       | I^{(1)}_{t,22}|
 \leq C \sum\limits_{k=3}^{1+\frac 4d} \lbb^{\frac {d}{2}k-4-d} \lbb^{d+3k-\frac d2 k}
 \leq C\sum\limits_{k=3}^{1+\frac 4d} \lbb^{3k-4}\leq C\lbb^4.
\end{align}
Similarly, using  again  \eqref{wn-lbbn}, \eqref{eta-L2}
and Gagliardo-Nirenberg's inequality we get
\begin{align} \label{esti-I1t-iii.2}
 | I^{(1)}_{t,4}|
 &\leq  C\sum\limits_{k=2}^{1+\frac 4d}\lbb^{\frac {d}{2}k-\frac d2-2}\|\eta\|_{L^2}\|R\|^k_{L^{2k}}
 \leq C\sum\limits_{k=2}^{1+\frac 4d}\lbb^{\frac {d}{2}k-\frac d2-2} \lbb^3 \lbb^{\frac d2 +3k-\frac d2 k}
  \leq C\lbb^4,
\end{align}
and
\begin{align} \label{esti-I1t-iii.3}
  | I^{(1)}_{t,5}|
 \leq C \sum\limits_{k=2}^{1+\frac 4d} \lambda^{-2}\int |w|^{1+\frac{4}{d}-k}|R|^{1+k}dx
  \leq C \sum\limits_{k=2}^{1+\frac 4d} \lbb^{3k-1}
  \leq C\lbb^{4}.
\end{align}

Moreover,
by \eqref{wn-lbbn} and \eqref{Mod-w-lbb},
  \begin{align}
  \frac{|\lbb \lbb_t +\g|}{\lbb^4}\|R\|^2_{L^2}
  \leq \frac{Mod(t)}{\lbb^{4}}\|R\|^2_{L^2}
  \leq C \lbb^4,
\end{align}
which implies that
\begin{align} \label{esti-I1t-ii.2}
   I^{(1)}_{t,6}
  = \frac{\g}{\lbb^4}\|R\|^2_{L^2}-\frac{\lbb \lbb_t +\g}{\lbb^4}\|R\|^2_{L^2}
  = \frac{\g}{\lbb^4}\|R\|^2_{L^2}+\calo(\lbb^4).
  \end{align}

Thus, we conclude from the estimates above that
\begin{align} \label{esti-nonl-w}
 \sum\limits_{j=2}^6 I^{(1)}_{t,j}
 = &- \frac{\gamma}{\lambda^2} {\rm Re}
    \int (1+\frac 2d) |w|^{\frac {4}{d}}|R|^2
    +\frac 12 (1+\frac 2d) |w|^{\frac {4}{d}-2}\ol{w}^2R^2
    +\frac 12 (\frac 2d -1) |w|^{\frac {4}{d}-2}w^2\overline{R}^2dx \nonumber \\
 &-\frac{\gamma}{\lambda}
   {\rm Re} \int
    (\frac{x-\alpha}{\lambda})\cdot\nabla \ol{w}
    \bigg\{\frac 2d (1+\frac 2d) |w|^{\frac {4}{d}-2}w|R|^2
     +\frac 1d (1+\frac 2d) |w|^{\frac {4}{d}-2}\ol{w}R^2  \nonumber \\
 &\qquad \qquad \qquad  \qquad  \qquad  \ \   +\frac 1 d(\frac 2d-1) |w|^{\frac 4d-4}w^3\overline{R}^2\bigg\}dx
   + \frac{\g}{\lbb^4} \|R\|_{L^2}^2
   +\calo ( \lbb^4).
\end{align}

(iii) {\it Estimates of the last term $I^{(1)}_{t,7}$.}
By the integration by parts formula and \eqref{Rn-wtQn},
\begin{align} \label{inhom-eta}
 & |{\rm Im}\<\Delta R - \lbb^{-2} R + f(u) - f(w), b\cdot \na R\>| \nonumber \\
 \leq& C \sum\limits_{k=1}^N \big( |{\rm Re}\int B_k \Delta R \nabla \phi_k\cdot \nabla \ol{R} dx |
+\frac{1}{\lambda^2} |{\rm Re}\int B_k R \nabla \phi_k\cdot\nabla \ol{R} dx |  \nonumber  \\
& +\int |B_k |R|^{\frac 4d}R \nabla \phi_k\cdot \nabla \ol{R}  |dx
  +  |\int B_k (f(u)-f(w)-|R|^{\frac{4}{d}}R) \nabla \phi_k\cdot \nabla \ol{R} dx |\big) \nonumber \\
:=& \sum\limits_{j=1}^4 I^{(1)}_{t,bj}.
\end{align}

We estimate each term $I^{(1)}_{t,bj}$, $1\leq j\leq 4$, separately below.
First, using the integration by parts formula to shift the derivatives of $R$
onto  $\phi_k$
we get
\begin{align*}
     &{\rm Re}\int \Delta R \nabla \phi_k\cdot \nabla \ol{R} dx
  = - {\rm Re}
    \int \na^2 \phi_k (\na \ol{R}, \na R)  dx +\frac 12 \int \Delta \phi_k  |\na R|^2dx, \\
     & {\rm Re}\int  R \nabla \phi_k\cdot\nabla \ol{R}  dx
    = -\frac 12 \int  \Delta \phi_k |R|^2 dx.
\end{align*}
Then, taking into account
the boundedness of
$\|\partial_x^\nu \phi_k\|_{L^\9}$
and $\sup_{0\leq t\leq T}|B_k(t)|$,
$1\leq |\nu|\leq 2$,
and using the Gagliardo-Nirenberg inequality we obtain that,
if $p:=1+\frac 4d$,
\begin{align}  \label{esti-I1t-J1J2J5}
 I^{(1)}_{t,b1}+I^{(1)}_{t,b2} + I^{(1)}_{t,b3}
 \leq& C(\|\nabla R\|_{L^2}^2
 + \frac{1}{\lbb^2} \|R\|_{L^2}
 +\|\na R\|_{L^2}\|R\|_{L^{2p}}^{p}) \nonumber \\
\leq& C(\|\nabla R\|_{L^2}^2+\frac{1}{\lambda^2}\|R\|^2_{L^2}
        +\|R\|_{L^2}^{p-dp(\frac 12 - \frac 1p)}\|\nabla R\|_{L^2}^{1+dp(\frac 12 - \frac 1p)})
\leq C\lbb^4,
\end{align}
where we also used \eqref{wn-lbbn} and \eqref{R-Lk-k}.
Moreover, we infer from \eqref{Rn-wtQn} and \eqref{Taylor} that
\begin{align*}
   I^{(1)}_{t,b4}
    \leq& C \sum\limits_{k=1}^{\frac 4 d}   \lbb^d
             \int |\na \phi_k (\lbb y +\a)(\na \ol{R} R^k)(\lbb y +\a) w^{1+\frac{4}{d}-k}(\lbb y +\a)| dy \nonumber \\
    \leq& C  \sum\limits_{k=1}^{\frac 4 d}  \lbb^{5+d}
              \int (1+|y|^6) |\na \ol{R} R^k| (\lbb y +\a) |w|^{1+\frac{4}{d}-k}(\lbb y +\a) dy \nonumber \\
    \leq& C \sum\limits_{k=1}^{\frac 4 d}   \lbb^{5}
              \int (1+|\frac{x-\a}{\lbb}|^6) |w|^{1+\frac{4}{d}-k}(x) |\na \ol{R} R^k|(x)dx.
\end{align*}
Since $\sup_x (1+|\frac{x-\a}{\lbb}|^6)|w|^{\frac{d+4}{d}-k}(x) \leq C \lbb^{\frac{d}{2}k-\frac d2-2}$,
by  \eqref{R-Lk-k},
\begin{align}  \label{esti-I1t-J3}
 I^{(1)}_{t,b4} \leq& C\sum\limits_{k=1}^{\frac 4 d}  \lbb^{3-\frac d2 + \frac d2 k} \|\na R\|_{L^2} \|R\|_{L^{2k}}^k
     \leq \sum\limits_{k=1}^{\frac 4 d}  \lbb^{5 +3k}
      \leq C\lbb^{4}.
\end{align}

The estimate of  the term involving the coefficient $c$
is easier.
Actually, since $|f(u)-f(w)| \leq C \sum_{k=1}^{1+4/d} |w|^{1+4/d-k}|R|^k$,
taking into account \eqref{R-Lk-k}
and the boundedness of $\|\partial^\nu_x \phi_k\|_{L^\9}$
and $\sup_{0\leq t\leq T}|B_k(t)|$, $0\leq |\nu|\leq 2$,
we get
\begin{align*}
    |{\rm Im }\<\Delta R - \lbb^{-2} R + f(u) - f(w), cR\>|
  \leq& C(\|\na R\|^2_{L^2} + \lbb^{-2} \|R\|^2_{L^2}
           + \sum\limits_{k=1}^{1+\frac 4d} \int |w|^{1+\frac 4d -k} |R|^{k+1} dx) \\
   \leq& C (\lbb^4 + \sum\limits_{k=1}^{1+\frac 4d} \lbb^{3k+1})
  \leq C \lbb^4.
\end{align*}

Thus,
we conclude from the estimates above that
\begin{align} \label{esti-I1t-iii.3}
    |I^{(1)}_{t,7}|
  \leq  C\lbb^4.
\end{align}

Therefore, plugging \eqref{esti-I1t-i},
\eqref{esti-nonl-w} and \eqref{esti-I1t-iii.3}  into \eqref{equa-I1t}
we obtain \eqref{I1t}.
\hfill $\square$

\begin{lemma} \label{Lem-I2t}
For all $t\in[t_*,t_n]$, we have
\begin{align} \label{I2t}
I^{(2)}_t \geq &-\frac{\gamma}{4\lambda^4}\int \Delta^2\chi_A(\frac{x-\alpha}{\lambda})|R|^2 dx
           +\frac{\gamma}{\lambda^2} {\rm Re} \int \nabla^2\chi_A(\frac{x-\alpha}{\lambda})(\nabla R,\nabla \ol{R}) dx   \nonumber \\
&+ \frac{\gamma}{\lambda}{\rm Re}\int
   \nabla\chi_A (\frac{x-\alpha}{\lambda} )\cdot\nabla \ol{w}
  \bigg\{\frac 2d(1+\frac 2d)|w|^{\frac {4}{d}-2}w|R|^2
    + \frac 1d (1+\frac 2d) |w|^{\frac {4}{d}-2}\ol{w}R^2   \nonumber\\
 &  \qquad \qquad \qquad \qquad \qquad  \qquad \ \  + \frac 1d (\frac 2d -1) |w|^{\frac 4d-4}w^3 \ol{R}^2 \bigg\}dx
     -C(A) \lbb^4.
\end{align}
\end{lemma}

{\bf Proof.}
Using \eqref{equa-R}, \eqref{etan-Rn}, \eqref{I2}
and integrating by parts formula
we first see that
\begin{align}\label{mf2}
I^{(2)}_t
=&-\frac{\lambda_t\gamma-\lambda\gamma_t}{2\lambda^2} {\rm Im}\int \nabla\chi_A(\frac{x-\alpha}{\lambda})\cdot\nabla R \ol{R}dx
+\frac{\gamma}{2\lambda}{\rm Im}\int \partial_t (\nabla\chi_A(\frac{x-\alpha}{\lambda}))\cdot\nabla R \ol{R}dx \nonumber \\
&+{\rm Im}\int (\frac{\gamma}{2\lambda^2} \Delta\chi_A(\frac{x-\alpha}{\lambda})R
+\frac{\gamma}{\lambda}  \nabla\chi_A(\frac{x-\alpha}{\lambda})\cdot\nabla R) \pa_t \ol{R} dx
=:  \sum\limits_{j=1}^3 I_{t,j}^{(2)}.
\end{align}

Since $\sup_{y} |\na^2 \chi_A(y)(1+|y|)| \leq C(A)$,
by Lemmas \ref{pre-est} and \ref{Lem-Mod-w-lbb},
\begin{align}
|\partial_t\nabla\chi_A(\frac{x-\alpha}{\lambda})|
=&|\nabla^2\chi_A(\frac{x-\alpha}{\lambda})\cdot
((\frac{x-\a}{\lbb})\cdot \frac{\lbb_t \lbb +\g}{\lbb^2}
  -(\frac{x-\a}{\lbb})\cdot \frac{\g}{\lbb^2}
  + \frac{\lbb \a_t - 2\beta}{\lbb^2}
  + \frac{2\beta}{\lbb^2} )|   \nonumber  \\
\leq& \frac{1}{\lambda^2}( Mod+P) C(A)
\leq \frac {1}{\lbb}C(A),
\end{align}
which along with   $\lf|\frac{\lambda_t\gamma-\lambda\gamma_t}{2\lambda^2}\rt|
\leq C\frac{Mod(t)}{\lambda^3} \leq C\lbb^2$
and \eqref{wn-lbbn}
yields  the bound
\be\ba\label{esti-I2.2}
 |I^{(2)}_{t,1} + I^{(2)}_{t,2}|
\leq C(A) \lbb^4.
\ea\ee

Moreover,
using equation \eqref{equa-R}, the expansion \eqref{Nf-def}
and integrating by parts we get
\begin{align} \label{esti-I2.0}
I^{(2)}_{t,3}=&-\frac{\gamma}{4\lambda^4}{\rm Re}\int \Delta^2\chi_A (\frac{x-\alpha}{\lambda} )|R|^2 dx
  +\frac{\gamma}{\lambda^2}{\rm Re}\int \nabla^2\chi_A(\frac{x-\alpha}{\lambda})(\nabla R,\nabla \ol{R}) dx  \nonumber \\
& -{\rm Re} \<\frac{\gamma}{2\lambda^2}\Delta\chi_A(\frac{x-\alpha}{\lambda} ) R
                + \frac{\gamma}{\lambda}\nabla\chi_A (\frac{x-\alpha}{\lambda} )\cdot\nabla R,
                  f'(w) \cdot R \> \nonumber \\
&-{\rm Re} \<\frac{\gamma}{2\lambda^2}\Delta\chi_A(\frac{x-\alpha}{\lambda} ) R
                + \frac{\gamma}{\lambda}\nabla\chi_A (\frac{x-\alpha}{\lambda} )\cdot\nabla R,
                N_f(w,R) \>  \nonumber \\
&-{\rm Re}  \<\frac{\gamma}{2\lambda^2}\Delta\chi_A(\frac{x-\alpha}{\lambda} ) R
                + \frac{\gamma}{\lambda}\nabla\chi_A (\frac{x-\alpha}{\lambda} )\cdot\nabla R,
                (b\cdot \na +c) R\>   \nonumber \\
&-{\rm Re}  \< \frac{\gamma}{2\lambda^2}\Delta\chi_A(\frac{x-\alpha}{\lambda} ) R
                + \frac{\gamma}{\lambda}\nabla\chi_A (\frac{x-\alpha}{\lambda} )\cdot\nabla R,
                 \eta \> \nonumber \\
=:&-\frac{\gamma}{4\lambda^4}{\rm Re}\int \Delta^2\chi_A (\frac{x-\alpha}{\lambda} )|R|^2 dx
  +\frac{\gamma}{\lambda^2}{\rm Re}\int \nabla^2\chi_A(\frac{x-\alpha}{\lambda})(\nabla R,\nabla \ol{R}) dx  \nonumber \\
  & + \sum\limits_{j=1}^4 I^{(2)}_{t,3j}.
\end{align}

Using the integration by parts formula we have
\begin{align*}
      {\rm Re} &\<\na \chi_A(\frac{x-\a}{\lbb}) \cdot  \na R, f'(w)\cdot R\>
   =  -\frac{1}{2\lbb} {\rm Re}
     \int \Delta \chi_A(\frac{x-\alpha}{\lambda} ) \ol{R} f'(w)\cdot R dx  \nonumber \\
    & - \frac 12 {\rm Re}
       \int \na \chi_A (\frac{x-\alpha}{\lambda} )
       \cdot \na w ((\partial_{uu}f(w) + \ol{\partial_{u\ol{u}}f}(w))|R|^2
        + \partial_{u\ol{u}}f(w)\ol{R}^2
        + \ol{\partial_{\ol{u}\ol{u}}f}(w) R^2).
\end{align*}
Then, using the explicit expression \eqref{f-linear} we obtain
\begin{align} \label{esti-I2.31}
 I^{(2)}_{t,31}=& \frac{\gamma}{\lambda}{\rm Re}\int \nabla\chi_A (\frac{x-\alpha}{\lambda} )\cdot\nabla \ol{w}
   \big\{\frac 2d(1+\frac 2d) |w|^{\frac {4}{d}-2}w|R|^2   \nonumber \\
  & \qquad \qquad   +\frac 1d(1+\frac 2d) |w|^{\frac {4}{d}-2}\ol{w}R^2
   +\frac 1d(\frac 2d-1) |w|^{\frac 4d-4}w^3\overline{R}^2\big\} dx.
\end{align}

For the nonlinear term  $I^{(2)}_{t,32}$,
similar arguments as in the estimate (\ref{esti-I1t-iii.1}) lead to
\begin{align} \label{esti-I2.1}
|I^{(2)}_{t,32}|
\leq C(A) \sum\limits_{k=2}^{1+\frac 4d} \int |w|^{1+\frac 4d -k}|R|^k
          (\lbb^{-1} |R| + |\na R|)dx
\leq C(A)\lbb^4.
\end{align}

Similarly, we have
\begin{align}
  |I^{(2)}_{t,33}|\leq & \frac{1}{\lambda}C(A)(\|R\|_{L^2}\|\nabla R\|_{L^2}+\|R\|_{L^2}^2)
+C(A)(\|\nabla R\|_{L^2}^2+\| R\|_{L^2}\|\nabla R\|_{L^2}) \nonumber \\
\leq&  C(A)(\|\nabla R\|_{L^2}^2+\frac{1}{\lambda^2}\|R\|_{L^2}^2)
\leq C(A) \lbb^4.
\end{align}

Moreover, by Lemmas \ref{Lem-Mod-w-lbb} and \ref{Lem-eeta},
\begin{align}  \label{esti-I2.3}
  |I^{(2)}_{t,34}|\leq & \frac{1}{\lambda}C(A)\|\eta\|_{L^2}
\|R\|_{L^2}+C(A)\|\nabla R\|_{L^2}\|\eta\|_{L^2}
\leq  C(A) \lbb^4.
\end{align}

Therefore,
plugging \eqref{esti-I2.31}-\eqref{esti-I2.3} into \eqref{esti-I2.0}
and using \eqref{esti-I2.2}
we obtain \eqref{I2t}.
\hfill $\square$

{\bf Proof of Theorem \ref{Thm-I-mono}.}
We infer from \eqref{I1t} and \eqref{I2t} that
for all $t\in[t_*,t_n]$,
\begin{align*}
&\frac{dI}{dt}
\geq \frac{\gamma}{\lambda^2} {\rm Re}\int \nabla^2\chi_A(\frac{x-\alpha}{\lambda})(\nabla R,\nabla \ol{R}) dx
       + \frac{\gamma}{\lambda^4}\int|R|^2dx   \nonumber \\
&-\frac{\gamma}{\lambda^2} {\rm Re}
  \int ((1+\frac 2d)|w|^{\frac 4d} |R|^2
  +\frac 12 (1+\frac 2d) |w|^{\frac {4}{d}-2} \ol{w}^2R^2
  +\frac 12 (\frac 2d-1) |w|^{\frac 4d-2} w^2 \ol{R}^2)dx\nonumber\\
&+\frac{\gamma}{\lambda} {\rm Re}\int
       (\nabla\chi_A(\frac{x-\alpha}{\lambda}) - (\frac{x-\alpha}{\lambda}))
       \cdot\nabla \ol{w}
       \big\{\frac 2d(1+\frac 2d) |w|^{\frac 4 d- 2} w|R|^2
         +\frac 1d (1+\frac 2d) |w|^{\frac {4}{d}-2} \ol{w}R^2 \big\}dx \nonumber \\
&+  \frac 1d (\frac 2d -1)\frac{\gamma}{\lambda} {\rm Re}
  \int (\nabla\chi_A(\frac{x-\alpha}{\lambda}) - (\frac{x-\alpha}{\lambda})) \cdot \na \ol{w}
   |w|^{\frac 4d-4}w^3 \ol{R}^2 dx  \nonumber \\
&-\frac{\gamma}{4\lambda^4}\int \Delta^2\chi_A(\frac{x-\alpha}{\lambda})|R|^2 dx
 -C(A) \lambda^4
\end{align*}
with $C(A)>0$  independent of $n$,
which, via \eqref{Rn-wtQn}, can be reformulated as follows
\begin{align}
\frac{dI}{dt}
   \geq&
   \frac{\gamma}{\lambda^4}
   \big(\int \nabla^2\chi (\frac{y}{A} )(\nabla \varepsilon,\nabla \ol{\varepsilon}) dy
          + \int |\varepsilon|^2dy
          -\int ((1+\frac 4d)Q^\frac{4}{d}\varepsilon_1^2+Q^\frac{4}{d}\varepsilon_2^2) dy \nonumber \\
&  \qquad -\frac{1}{4A^2}\int \Delta^2\chi (\frac{y}{A} )|\varepsilon|^2 dy\big) \nonumber \\
&+\frac 2d  \frac{\gamma}{\lambda^4}\int (A\nabla\chi (\frac{y}{A} )-y ) \cdot \nabla Q
Q^{\frac 4d -1} ((1+\frac 4d)\varepsilon_1^2+ \varepsilon_2^2)dy
 -C(A)\lambda^4.
\end{align}

We claim that
\be\label{matrix}
\int \nabla^2\chi (\frac{y}{A})(\nabla \varepsilon,\nabla \ol{\varepsilon}) dy
\geq \int \psi^{\prime\prime} (|\frac{y}{A}| )|\nabla \varepsilon |^2 dy.
\ee
To this end,
since the case where $d=1$ is obvious,
we only need to treat the case where $d=2$ below.
Since $\chi(y)=\psi(|y|)$, we have
\begin{align*}
  \partial_{y_{i}y_{j}}\chi(y)
  =\psi^{\prime \prime}(|y|)\frac{y_iy_j}{|y|^2}
    + \psi^{\prime}(|y|)\frac{\delta_{ij}}{|y|}
    -\psi^{\prime}(|y|)\frac{y_iy_j}{|y|^3},\ \ i,j=1,2.
\end{align*}
Then, by the condition (\ref{chi}), we have
\begin{align}
\int \nabla^2\chi(\frac{y}{A}) & (  \nabla \varepsilon,\nabla\bar{\varepsilon})dy
=\int \psi^{\prime \prime}(|\frac{y}{A}|)|\na \ve|^2dy \nonumber \\
&  +\int (\frac{A}{|y|} \psi^\prime  (|\frac{y}{A}|)  - \psi''(\frac y A) )
  (\frac{y_1^2}{|y|^2}|\partial_{2}\varepsilon|^2
+\frac{y_2^2}{|y|^2}|\partial_{1}\varepsilon|^2
-2 \frac{y_1y_2}{|y|^2} {\rm Re} (\partial_{1}\varepsilon\partial_{2}\bar{\varepsilon}) ) dy  \nonumber \\
&\geq \int \psi''(|\frac{y}{A}|)|\na \ve|^2dy,
\end{align}
which implies  \eqref{matrix}, as claimed.

Thus, applying Corollary \ref{Cor-coer-f-local} with $\Phi(x) :=\psi''(|x|)$ we obtain
that for some $\nu>0$
\begin{align} \label{dIdt-esti}
\frac{dI}{dt}
\geq& \nu  \frac{\gamma}{\lambda^4} \int \psi'' (|\frac{y}{A}| )(|\varepsilon|^2+|\nabla \varepsilon|^2) dy
    -\frac{1}{4A^2} \frac{\g}{\lbb^4} \int \Delta^2\chi (\frac{y}{A} )|\varepsilon|^2 dy  \nonumber \\
&+\frac 2d \frac{\gamma}{\lambda^4}\int (A\nabla\chi (\frac{y}{A} )-y ) \cdot \nabla Q
Q^{\frac 4d -1} ((1+\frac 4d) \varepsilon_1^2+ \varepsilon_2^2)dy
 -C(A)\lambda^4.
\end{align}
Note that,
since $\psi'' (\lf|\frac{y}{A}\rt| ) \to 1$
and $\frac{1}{4A^2} \Delta^2\chi (\frac{y}{A} ) \to 0$ as $A\to \9$,
we infer that for $A$ large enough,
$$\frac{1}{4A^2} |\Delta^2\chi (\frac{y}{A} )| \leq \frac {\nu}{4} \psi'' (|\frac{y}{A}| ).$$
Similarly,
since $A \na \chi (\frac yA) = y$ for $|y|\leq A$
and $Q$ is exponentially decay at infinity,
we have that
$|A\nabla\chi (\frac{y}{A} )-y| \leq A(1+|y|)e^{-\delta A}$,
which yields that  for $A$ large enough
$$\frac 2d(2+\frac 4d) |A\nabla\chi (\frac{y}{A} )-y ||\na QQ^{\frac 4d -1}| \leq \frac {\nu}{4} \psi'' (|\frac{y}{A}| ).$$

Therefore,
in view of $\psi''(r) \geq \delta e^{-r}$ for some $\delta>0$,
we obtain that for $A$ large enough
\begin{align}
  \frac{dI}{dt}
\geq& \frac{\nu}{2} \frac{\gamma}{\lambda^4} \int \psi'' (|\frac{y}{A}| )(|\varepsilon|^2+|\nabla \varepsilon|^2) dy
    - C(A) \lbb^4 \nonumber \\
\geq& \frac{\nu \delta}{2} \frac{\gamma}{\lambda^4} \int  e^{-\frac{|y|}{A}} (|\varepsilon|^2+|\nabla \varepsilon|^2) dy
     - C(A) \lbb^4,
\end{align}
which yields \eqref{dIt-mono}.
The proof of Theorem \ref{Thm-I-mono} is complete.  \hfill $\square$

\subsection{Proof of Proposition \ref{Prop-dec-un-t0-boot}}  \label{Subsec-proof-bootstrap}
{\it (i) Estimate of $R_n$.}
By \eqref{un-dec}  and \eqref{def-I},
for $t\in[t_*,t_n]$,
\begin{align} \label{I.0}
   I_n(t)
   =& \frac{1}{2\lbb_n^2} \int |\na \ve_n|^2 + |\ve_n|^2 dy
    - \frac{1}{\lbb_n^2}
     {\rm Re} \int F(Q_{\calp_n} +\ve_n) - F(Q_{\calp_n}) - f(Q_{\calp_n})\ol{\ve_n} dy   \nonumber \\
     & +\mathcal{O}(\|\nabla R_n(t)\|_{L^2}\|R_n(t)\|_{L^2}).
\end{align}
This along with \eqref{F-Taylor*} yields that
\begin{align*}
   &\frac{1}{ \lbb_n^2} {\rm Re} \int F(Q_{\calp_n} +\ve_n) - F(Q_{\calp_n}) - f(Q_{\calp_n})\ol{\ve_n} dy \nonumber \\
   =& \frac{1}{2\lbb_n^2}
      {\rm Re} \int Q^\frac 4d (\frac 2d \ve_n^2 + (1+\frac 2d) |\ve_n|^2) dy
      + o (\frac{1}{\lbb_n^{2}} \|\ve_n\|^2_{H^1}),
\end{align*}
which yields immediately that
\begin{align*}
  I_n(t)  = \frac{1}{2\lbb_n^2} L(\ve_n, \ve_n)
      + \frac{1}{2\lbb_n^2} o(\|\ve_n\|_{H^1}^2)
     +\mathcal{O}(\|\nabla R_n(t)\|_{L^2}\|R_n(t)\|_{L^2}).
\end{align*}
Then,  using Corollary \ref{Cor-coer-f-H1}
and similar arguments as in the proof of \eqref{ve-quadra-Eu}
we obtain that for some $\nu>0$ and for $T$ small enough,
\begin{align} \label{I-lowbdd}
   I_n(t) \geq  \frac{\nu}{2\lbb_n^2} \|\ve_n\|_{H^1}^2
              - C(A) (\|\nabla R_n(t)\|_{L^2}\|R_n(t)\|_{L^2})
        \geq  \frac{\nu}{4}(\|\nabla R_n(t)\|_{L^2}^2+\frac{1}{\la^2_n(t)}\|R_n(t)\|_{L^2}^2).
\end{align}

Moreover,
Theorem \ref{Thm-I-mono} yields that for any $t\in[t_*,t_n]$,
\begin{align}  \label{dIdt-lowbdd}
\frac{dI_n}{dt}
\geq   - C(A)\lbb_n^4.
\end{align}

Thus, we infer from \eqref{I-lowbdd} and \eqref{dIdt-lowbdd} that for any $t\in[t_*,t_n]$,
\begin{align*}
    \frac{\nu}{4}(\|\nabla R_n(t)\|_{L^2}^2+\frac{1}{\la^2_n(t)}\|R_n(t)\|_{L^2}^2)
\leq&  I_n(t)=I_n(t_n)-\int_{t}^{t_n}\frac{dI_n}{dr}(r) dr \nonumber \\
\leq&  I_n(t_n)+C(A)\int_{t}^{t_n}\la_n^4(r) dr.
\end{align*}
Taking into account $I_n(t_n)=0$ and
using \eqref{wn-lbbn}  we obtain that
for $T$ small enough,
\be\ba
\|\nabla R_n(t)\|_{L^2}^2+\frac{1}{\la^2_n(t)}\|R_n(t)\|_{L^2}^2
\leq \frac{4C(A)}{\nu}
\int_{t}^{T} (T-r)^4 dr
\leq \frac{1}{2}(T-t)^4,
\ea\ee
which yields \eqref{wn-Tt-boot-2}.

{\it $(ii)$ Estimates of $\lambda_n$ and $\g_n$.}
First note that, by Lemma \ref{Lem-Mod-w-lbb},
\be
  |(\frac{\gamma_n}{\lambda_n} )_t|
=\frac{|\lambda_n^2\gamma_t-\lambda_n\lambda_{n,t}\gamma_n|}{\lambda_n^3}
\leq 2 \frac{Mod_n(t)}{\lambda_n^3}
\leq C \lambda_n^{2},
\ee
which along with $(\frac {\gamma_n}{ \lambda_n})(t_n)=1$ and \eqref{lbbn-Tt*}
implies that for  $T$ small enough,
\be  \label{gamlbb-1}
 | (\frac{\gamma_n}{\lambda_n} )(t)-1 |
\leq\int_{t}^{t_n}  | (\frac{\gamma_n}{\lambda_n} )_r |dr
\leq C  (T-t)^{3}\leq \frac{1}{2}(T-t)^{2+6\delta}.
\ee
Then, taking into account \eqref{lbbn-Tt*} and \eqref{Mod-w-lbb}
we obtain
\begin{align*}
|(\lambda_n - (T-t))_t|
=|\lambda_{n,t}+\frac{\gamma_n}{\lambda_n}+1-\frac{\gamma_n}{\lambda_n}|
\leq\frac{Mod_n}{\lambda_n}+\frac{1}{2} (T-t)^{2+6\delta}
\leq C  (T-t)^{2+6\delta},
\end{align*}
which implies that for $T$ possibly even smaller,
such that $C T^{\delta} \leq \frac{1}{2}$,
\begin{align} \label{lbb-Tt*}
|\lambda_n - (T-t)|
\leq\int_{t}^{t_n} \lf|(\lambda_n-(T-r))_r\rt|dr
\leq \frac{1}{2}(T-t)^{3+ 5 \delta  }.
\end{align}
Similarly,
by \eqref{lbbn-Tt*}, \eqref{Mod-w-lbb} and  \eqref{gamlbb-1},
\be \label{g-Tt*}
|\gamma_n - (T-t)|\leq \frac{1}{2}(T-t)^{3+ 5 \delta  }.
\ee
Thus, we obtain \eqref{lbbn-Tt-boot-2}.

{\it $(iii)$ Estimates of $\beta_n$ and $\alpha_n$.}
By \eqref{lbbn-Tt*} and \eqref{b-g-ve-lbb},
\be
|\beta_n|^2
\leq C|\lambda_n^2-\gamma_n^2|
+ C(T-t)^{5}
\leq C \lbb_n^2 (1-\frac {\g_n}{ \lbb_n}) + C \lbb_n^{5},
\ee
which along with \eqref{lbbn-Tt*} and \eqref{gamlbb-1}
yields that  for $T$ small enough,
\be \label{b-Tt*}
|\beta_n(t)|
\leq  C\lbb_n |1-\frac {\g_n}{ \lbb_n}|^\frac 12 + C \lbb_n^{\frac{5}{2}}
\leq C\lambda_n^{2+ 3 \delta }
\leq\frac{1}{2}(T-t)^{2+2 \delta }.
\ee

Moreover, using again \eqref{Mod-w-lbb} and \eqref{b-Tt*}
we have
\be
|\alpha_{n,t}|=
|\frac{\lambda_n\alpha_{n,t}-2\beta_n}{\lambda_n}+\frac{2\beta_n}{\lambda_n}|
\leq\frac{Mod_n}{\lambda_n}+\frac{2\beta_n}{\lambda_n}
\leq C\lambda_n^{1+ 2 \delta },
\ee
which yields that for sufficiently small $T$,
\be
\alpha_n(t)\leq\int_{t}^{t_n}|\alpha_{n,r}|dr
\leq \frac{1}{2}(T-t)^{2+  \delta  },
\ee
thereby yielding \eqref{anbn-Tt-boot-2}.

{\it $(iv)$ Estimate of $\theta_n$.}
Using \eqref{lbbn-Tt*}, \eqref{Mod-w-lbb}, \eqref{lbb-Tt*} and \eqref{b-Tt*} we get
\begin{align}  \label{esti-thetan.0*}
(\theta_n-\frac{1}{T-t})_t
=&\theta_{n,t}-\frac{1}{(T-t)^2}
=\frac{\lbb_n^2 \theta_{n,t} -1 -|\beta_n|^2}{\lbb_n^2}
  +\frac{|\beta_n|^2}{\lbb_n^2}
  + \frac{1}{\lbb_n^2}
  - \frac{1}{(T-t)^2}.
\end{align}
Note that, by \eqref{lbbn-Tt*} and \eqref{lbb-Tt*},
\begin{align*}
    |\frac{1}{\lbb_n^2} - \frac{1}{(T-t)^2}|
   \leq \frac{|\lbb_n-(T-t)||\lbb_n+(T-t)|}{\lbb_n^2(T-t)^2 }
   \leq C \lbb_n^{5\delta},
\end{align*}
which along with \eqref{lbbn-Tt*}, \eqref{Mod-w-lbb}, \eqref{b-Tt*}
and \eqref{esti-thetan.0*} yields that
\begin{align}  \label{esti-thetan.0}
   |(\theta_n-\frac{1}{T-t} )_t|
 \leq \frac{Mod_n}{\lbb_n^2}
      + \frac{|\beta_n|^2}{\lbb_n^2}
      + C \lbb_n^{5\delta  }
 \leq C (T-t)^{5 \delta  }.
\end{align}
Thus, integrating both sides above and taking  $T$ very small we obtain
\be
 |\theta_n-\frac{1}{T-t} |\leq \int_{t}^{t_n} | (\theta_n-\frac{1}{T-r} )_r |dr
\leq \frac{1}{2}(T-t)^{1+ 4 \delta },
\ee
which implies \eqref{thetan-Tt-boot-2}.
Therefore, the proof of Proposition \ref{Prop-dec-un-t0-boot} is complete.
\hfill $\square$

\subsection{Proof of Theorem \ref{Thm-dec-un-t0}} \label{Subsec-Proof-Unifesti}
Since the bounds in the proof of Proposition \ref{Prop-dec-un-t0-boot}
are uniform of $n$,
we may take a universal sufficiently small $\tau^*$
such that the estimates in Proposition \ref{Prop-dec-un-t0-boot} hold.
Below we take $T \in (0,\tau^*]$ fixed.

Let
$H(t)$ denote the statement that
the geometrical decomposition \eqref{un-dec} and estimates \eqref{wn-Tt}-\eqref{thetan-Tt} hold on $[t,t_n]$,
and let
$C(t)$ denote the statement that the decomposition \eqref{un-dec} and estimates
\eqref{wn-Tt-boot-2}-\eqref{thetan-Tt-boot-2} hold on $[t,t_n]$,
$0\leq t\leq t_n$.

It is clear that $H(t_n)$ is true.
By Lemma \ref{Lem-imp-thm}
and the continuity of $u_n$ and $\calp_n$,
we also have that
if $C(t)$ is true for some $t\in [0, t_n]$,
then $H(t')$ holds for all $t'$ in a neighborhood of $t$.
Moreover,
Proposition \ref{Prop-dec-un-t0-boot}
yields  that
if $H(t)$ is true for some $t\in [0, t_n]$,
then $C(t)$ is also true.

Furthermore,
we claim that
if $\wt t_m$, $m\geq 1$, is a sequence in $[0, t_n]$
which converges to another $\wt t_* \in [0, t_n]$
and $C(\wt t_m)$ is true for all $\wt t_m$, $m\geq 1$,
then $C(\wt t_*)$ is also true.

To this end,
by virtue of estimates \eqref{lbbn-Tt-boot-2}-\eqref{thetan-Tt-boot-2},
Lemma \ref{Lem-Mod-w-lbb} and \eqref{wn-lbbn},
we infer that
the derivatives of the modulation parameters $\calp'_n(\wt t_m)$
are uniformly bounded on $(\wt t_*,\wt t_m]$, $m\geq 1$,
which yields that there exists a unique $\wt \calp_n^*$
such that $\lim_{m\to \9} \calp_n(\wt t_m) =\wt \calp_n^*$.
Hence,
letting $\calp_n(\wt t_*):= \wt \calp^*_n$
we have that $\calp_n$ is continuous on $[\wt t_*, t_n]$.
Taking into account the decomposition \eqref{un-dec} and estimates \eqref{lbbn-Tt-boot-2}-\eqref{thetan-Tt-boot-2}
we infer that
$w_n(\wt t_m) \to w_n(\wt t_*)$ in $H^1$,
and $|x|^2w_n(\wt t_m)$, $\wt \rho_n(\wt t_m)$
and $R_n(\wt t_m)$ converge to
$|x|^2w_n(\wt t_*)$, $\wt \rho_n(\wt t_*)$
and $R_n(\wt t_*):= u_n(\wt t_*) - w_n(\wt t_*)$, respectively,
in the space $L^2$
(see similar arguments as in the proof  \eqref{4cor} below).
So, the decomposition \eqref{un-dec}
and the orthogonality condition \eqref{ortho-cond-Rn-wn} also hold at time $\wt t_*$.
Moreover,
in view of the continuity of $u_n$ and $\calp_n$,
we also infer that the estimates \eqref{wn-Tt-boot-2}-\eqref{thetan-Tt-boot-2}
hold on $[\wt t_*, t_n]$.
Thus, we conclude that
the statement $C(\wt t_*)$ is also true, as claimed.

Therefore,
by virtue of the abstract bootstrap principle
(see \cite[Proposition 1.21]{T06}),
we prove Theorem \ref{Thm-dec-un-t0}.
\hfill $\square$

\section{Proof of main results}  \label{Sec-Proof-Main}

Let us start with the global well-posedness result in Theorem \ref{Thm-u-GWP}.

{\bf Proof of Theorem \ref{Thm-u-GWP}.}
In view of Remark \ref{Rem-u-weak-mild},
we have that $\bbp$-a.s. there exists a unique solution $u$
to equation \eqref{equa-u-RNLS} on $[0,\tau^*)$ with $u(0)=u_0 \in H^1$,
where $\tau^*\in (0,\9]$ is some positive random variable.
Hence, we only need to prove that $\tau^*=\9$, $\bbp$-a.s..

For this purpose, similarly to \eqref{dtE}, we have
\begin{align} \label{dtE*}
\frac{d}{dt}E(u)
=&-2\sum\limits_{k=1}^N B_k {\rm Re}\int \nabla^2 \phi_k(\nabla u,\nabla \ol{u})dx
+\frac{1}{2}\sum\limits_{k=1}^N B_k \int \Delta^2 \phi_k|u|^2dx  \nonumber \\
& +\frac{2}{d+2}\sum\limits_{k=1}^N  B_k \int\Delta \phi_k|u|^{2+\frac{4}{d}}dx
 -\sum\limits_{j=1}^d {\rm Im}\int \na(\sum\limits_{k=1}^N \partial_j \phi_k B_k)^2 \cdot \na u \ol{u} dx.
\end{align}
Then, for any $T\in(0,\tau^*)$,
since $\phi_k\in C^\9_b$ and $B_k \in C([0,T])$, $\bbp$-a.s., $1\leq k\leq N$,
using H\"older's inequality we obtain that $\bbp$-a.s. for any $t\in [0,T]$,
\begin{align}
   E(u(t) ) \leq E(u_0)
                 + C(\tau^*) \int_0^t (\|u(s) \|_{L^2}^2+ \|\na u(s) \|_{L^2}^2 + \| u(s) \|_{L^{2+\frac 4d}}^{2+\frac 4d}) ds.
\end{align}
Note that,
$\|u(s) \|_{L^2}^2 = \|u_0 \|_{L^2}^2$,
and by the Gagliardo-Nirenberg inequality \eqref{G-N},
\begin{align*}
   \|\na u(s) \|_{L^{2+\frac 4d}}^{2+\frac 4d}
   \leq C \|  u(s) \|_{L^2}^{\frac 4d}  \|\na u(s) \|_{L^2}^{2}
   = C   \|  u_0  \|_{L^2}^{\frac 4d}  \|\na u(s) \|_{L^2}^{2}.
\end{align*}
Thus, we obtain
\begin{align} \label{bdd-E.1}
   E(u(t) ) \leq C(\tau^*)
                 + C(\tau^*) \int_0^t \|\na u(s) \|_{L^{2}}^{2}  ds,\ \ 0<t\leq T.
\end{align}

Moreover,
applying the sharp Gagliardo-Nirenberg  inequality (cf. \cite[(III.5)]{Wenn},
\cite[(4)]{PR07}) we get
\begin{align} \label{bdd-E.2}
   (1-(\frac{\|u_0\|_{L^2}}{\|Q\|_{L^2}})^{\frac{4}{d}})\|\na u(t)\|_{L^2}^2
   =(1-(\frac{\|u(t)\|_{L^2}}{\|Q\|_{L^2}})^{\frac{4}{d}} ) \|\na u(t)\|_{L^2}^2 \leq 2 E(u(t)).
\end{align}
Thus, in view of $\|u_0\|_{L^2} < \|Q\|_{L^2}$,
putting together \eqref{bdd-E.1} and \eqref{bdd-E.2}
and using Gronwall's inequality
we obtain $\bbp$-a.s.
\begin{align}
   \sup\limits_{0\leq t\leq T}\|\na u(t)\|_{L^2}^2 \leq C(\tau^*) <\9,
\end{align}
which yields
the boundedness of $\sup_{0\leq t<\tau^*}\|\na u(t)\|^2_{L^2}<\9$
by letting $T\to \tau^*$.
Therefore,
using the blow-up alternative result
we obtain  $\tau^*=\9$, $\bbp$-a.s.,
and  finish the proof.
\hfill $\square$

{\bf Proof of Theorem \ref{Thm-u-blowup}.}
Let $\tau^*$ be as in Theorem \ref{Thm-dec-un-t0}
and $T\in (0,\tau^*]$ fixed below.
By virtue of Theorem \ref{Thm-dec-un-t0},
we have
\be \label{dec-un-wnRn}
u_n(t,x)
=w_n(t,x)+R_n(t,x)
=\lambda^{-\frac d2}_n(t)Q_{\calp_n}  (t,\frac{x-\alpha_n(t)}
{\lambda_n(t)} )e^{i\theta_n(t)}+R_n(t,x),\ \ \forall t\in[0,t_n],
\ee
with the modulation parameters
\be
    \calp_n(t_n): = \lf(\la_n(t_n),\alpha_n(t_n),\beta_n(t_n),\gamma_n(t_n),\theta_n(t_n)\rt)
    =(T-t_n,0,0,T-t_n,\frac{1}{T-t_n}),
\ee
and so
\be
u_n(t_n) = S_T(t_n),\  \ R_n(t_n)=0,
\ee
where $S_T$ is the pseudo-conformal blow-up solution given by \eqref{S}.

Moreover,
the estimates \eqref{wn-Tt}-\eqref{thetan-Tt} hold for all $t\in[0,t_n]$.
In particular,
\begin{align} \label{lbbn-Tt-equal}
     \frac 12 (T-t) \leq |\lbb_n(t)| \leq 3(T-t), \ \ \forall\ t\in [0, t_n].
\end{align}

Since
\begin{align} \label{Rn-Q}
   \|w_n(0)\|_{L^2} = \|Q\|_{L^2}, \ \
   \|\na w_n(0)\|_{L^2}
   \leq \frac{C}{T} ( \|Q\|_{H^1} + \| |\cdot| Q\|_{L^2}),
\end{align}
where $C$ is independent of $n$,
taking into account the uniform $H^1$-boundedness of $R_n$ in \eqref{wn-Tt}
we infer that $\{u_n(0)\}$ is uniformly bounded in $H^1$.
This yields that  up to a subsequence (still denoted by $\{n\}$),
for some $u_0\in H^1$,
\begin{align}
    u_n(0) \rightharpoonup u_0,\ \ weakly\ in\ H^1,\ as\ n \to \9.
\end{align}

We claim that
\begin{align} \label{unt0-u0-L2}
    u_n(0) \to u_0,\ \  in\ L^2,\ as\ n \to \9.
\end{align}
In particular,
since $\|u_n(0)\|_{L^2} = \|u_n(t_n)\|_{L^2} = \|Q\|_{L^2}$,
we have that $\|u_0\|_{L^2} = \|Q\|_{L^2}$.

For this purpose,
since
$u_n(0)\rightarrow u_0$ in $L^2_{loc}(\R^d)$ by the compactness imbedding,
it suffices to prove the uniform integrability of $\{u_n(0)\}$, i.e.,
\be  \label{un-uninteg-L2}
\lim_{A\rightarrow\infty} \sup\limits_{n\geq 1}\|u_n(0)\|_{L^2(|x|>2A)}=0.
\ee

In order to prove \eqref{un-uninteg-L2},
we take a nonnegative  function $\vf\in C^{\infty}$ such that
$\vf(x)=0$ if $|x|\leq 1$,
$\vf(x) =1$ if $|x|\geq 2$,
$|\na \vf(x)|\leq 2$ for $x\in \R^d$.
Set $\vf_A(\cdot) :=\vf\lf(\frac{\cdot}{A}\rt)$, $A>0$.
Then,
\be \label{un-uninteg.0}
    |\int \vf_A|u_n(0)|^2dx|\leq
\int \vf_A|u_n(t_n)|^2dx
+\int_{0}^{t_n}|\frac{d}{dt}\int \vf_A|u_n(t)|^2dx|dt.
\ee
Since $u_n(t_n)=S_T(t_n)$,
we have that as $A\to \9$,
\begin{align} \label{un-uninteg.1}
    \sup\limits_{n\geq 1} |\int \vf_A|u_n(t_n)|^2dx|
    \leq C\sup\limits_{n\geq 1} \int_{|y|\geq\frac{A}{T-t_n}}|Q(y)|^2dy
    \leq C\int_{|y|\geq \frac{A}{T}} |Q(y)|^2 dy
    \to 0.
\end{align}
Moreover, by \eqref{b}, \eqref{c} and \eqref{equa-un-tn},
\begin{align} \label{dtun-A}
|\frac{d}{dt}\int \vf_A|u_n(t)|^2dx|
=&2|{\rm Im}\int\nabla \vf_A\cdot\nabla u_n(t) \ol{u_n}(t)dx
                        +\sum\limits_{k=1}^N  \int B_k \nabla \vf_A\cdot\nabla \phi_k |u_n|^2(t)dx| \nonumber \\
\leq&\frac{C}{A}(\|u_n(t)\|^2_{L^2(|x|>A)} + \|\nabla u_n(t)\|^2_{L^2(|x|>A)}).
\end{align}
Note that,
by \eqref{wn-Tt}, \eqref{dec-un-wnRn} and \eqref{Rn-Q},
\begin{align} \label{un-L2-A}
   \|u_n(t)\|_{L^2(|x|>A)}
   \leq \|w_n\|_{L^2} + \|R_n\|_{L^2}
   \leq C(\|Q\|_{L^2} + T^3)<\9.
\end{align}
Moreover,
using \eqref{wn-Tt} and  \eqref{dec-un-wnRn} again we have
\begin{align*}
 \|\nabla u_n(t)\|^2_{L^2(|x|>A)}
 \leq&C(\|\nabla w_n(t)\|^2_{L^2(|x|>A)}+\|\nabla R_n(t)\|^2_{L^2(|x|>A)})   \\
\leq&  \frac{ C}{\lbb_n^2(t)}\int_{|y|\geq \frac{A-\a_n(t)}{\lbb_n(t)}} |\na Q(y)|^2 + (1+|y|^2) Q^2 dy  + C.
\end{align*}
Then, taking into account the exponential decay of $Q$
we obtain
\begin{align} \label{naun-L2-A}
    \|\nabla u_n(t)\|_{L^2(|x|>A)}
    \leq \frac{C}{\lbb_n^2(t)}  \int_{|y|\geq \frac{A-1}{\lbb_n(t)}} e^{-\delta |y|} dy  + C
    \leq C(\frac{1}{\lbb_n^2(t)} e^{- \frac{ \delta(A-1)}{2\lbb_n(t)}} +1)
    \leq C,
\end{align}
where $C, \delta>0$ are independent of $n$.

Thus, plugging \eqref{un-L2-A} and \eqref{naun-L2-A} into \eqref{dtun-A}
we obtain
\begin{align} \label{un-uninteg.2}
   \sup\limits_{n\geq 1}| \frac{d}{dt} \int \vf_A |u_n(t)|^2 dx | \leq \frac C A \to 0,\ \ as\ A\to \9,
\end{align}
which along with \eqref{un-uninteg.0} and \eqref{un-uninteg.1} implies \eqref{un-uninteg-L2},
thereby yielding \eqref{unt0-u0-L2}, as claimed.

Now,
since $u_n$ solves the equation \eqref{equa-u-RNLS}
on $[0,t_n]$ with $\lim_{n\rightarrow\infty}t_n=T$,
and
by \eqref{unt0-u0-L2},
$u_n(0)$ converge strongly to $u_0$ in the space $L^2$ as $n\to \9$,
the $L^2$ local well-posedness theory yields that
there exists a unique $L^2$-solution $u$ to \eqref{equa-u-RNLS} on $[0, T)$
satisfying
\begin{align} \label{un-u-0-L2}
\lim_{n\rightarrow \infty}\|u_n(t)-u(t)\|_{L^2}=0,\ \ t\in [0, T).
\end{align}
Moreover,
since $u_0\in H^1$,
the preservation of $H^1$-regularity
implies that $u(t) \in H^1$ for any $t\in [0, T)$.

Below we prove that for some $\delta >0$,
\be \label{u-S-0-H1}
\|u(t)-S_T(t)\|_{H^1} \leq C(T-t)^\delta,
\ee
which yields that $u$ behaves asymptotically like
the pseudo-conformal blow-up solution $S_T$ as $t\to T$,
In particular, $u$ blows up at time $T$.

In order to prove \eqref{u-S-0-H1},
using \eqref{dec-un-wnRn} we have
\be \label{un-St}
\|u_n(t)-S_T(t)\|_{H^1}\leq\|R_n(t)\|_{H^1}+\|w_n(t)-S_T(t)\|_{H^1}, \ \ \forall t\in [0,T).
\ee

The uniform estimate \eqref{wn-Tt} yields immediately that
\be \label{esti-wn}
\|R_n(t)\|_{H^1}\leq C (T-t)^2,\ \ \forall t\in [0, T).
\ee

Moreover,
since
\be
S_T(t)=\lambda^{-\frac d2}_0(t)Q_{\calp_0} (t,\frac{x-\alpha_0(t)}
{\lambda_0(t)})e^{i\theta_0(t)},
\ee
where $Q_{\calp_0} (t,y) := Q(y) e^{i(\beta_0\cdot y - \frac{\g_0}{4} |y|^2)}$
and $\calp_0:= \lf(\la_0(t),\alpha_0(t),\beta_0(t),\gamma_0(t),\theta_0(t)\rt) =(T-t,0,0,$ $T-t,\frac{1}{T-t})$,
taking into account \eqref{dec-un-wnRn} we get
\begin{align} \label{Rn-S-L2.0}
\|w_n(t)-S_T(t)\|_{L^2}
=& \|\lambda_n^{-\frac{d}{2}}(t) Q_{\calp_n} (t, \frac{x-\alpha_n(t)}
{\lambda_n(t)}) e^{i\theta_n(t)}
-\lambda_0^{-\frac{d}{2}}(t) Q_{\calp_0}  (t, \frac{x-\alpha_0(t)}
{\lambda_0(t)}) e^{i\theta_0(t)} \|_{L^2} \nonumber \\
\leq&
\lambda_0^{-\frac{d}{2}}(t)\| Q_{\calp_n} (t, \frac{x-\alpha_n(t)}
{\lambda_n(t)})- Q_{\calp_0} (t, \frac{x-\alpha_n(t)}{\lambda_n(t)})\|_{L^2} \nonumber \\
&+\lambda_0^{-\frac{d}{2}}(t)\|Q_{\calp_0} (t, \frac{x-\alpha_n(t)}
{\lambda_n(t)})
-Q_{\calp_0} (t, \frac{x-\alpha_0(t)}
{\lambda_0(t)})\|_{L^2} \nonumber \\
& +( \lambda_0^{-\frac{d}{2}}(t) {|\lambda_n^{\frac{d}{2}}(t)-\lambda_0^{\frac{d}{2}}(t)|}  +
  |\theta_n(t) - \theta_0(t)|)\|Q\|_{L^2}.
\end{align}
Note that,
by the change of variables,
\begin{align} \label{Rn-S-L2.1}
  &\lambda_0^{-\frac{d}{2}}(t)\|Q_{\calp_n} (t, \frac{x-\alpha_n(t)}
{\lambda_n(t)})-Q_{\calp_0} (t, \frac{x-\alpha_n(t)}{\lambda_n(t)})\|_{L^2}  \nonumber \\
=& (\frac{\lbb_n(t)}{\lbb_0(t)})^{\frac d2} \|Q_{\calp_n}(t)  - Q_{\calp_0}(t) \|_{L^2}
\leq C (\frac{\lbb_n(t)}{\lbb_0(t)})^{\frac d2} (|\beta_n(t) - \beta(t)| + |\g_n(t) - \g_0(t)|).
\end{align}
Moreover, using the change of variables again
and the mean value theorem we get
\begin{align}  \label{Rn-S-L2.2}
   & \lambda_0^{-\frac{d}{2}}(t)\|Q_{\calp_0} (t, \frac{x-\alpha_n(t)}
{\lambda_n(t)})
-Q_{\calp_0} (t, \frac{x-\alpha_0(t)}
{\lambda_0(t)})\|_{L^2} \nonumber \\
    =& \|Q(\frac{\lbb_0(t)}{\lbb_n(t)} y + \frac{\a_0(t)- \a_n(t)}{\lbb_n(t)}) - Q(y)\|_{L^2}
    \leq  C (|\frac{\lbb_0(t)}{\lbb_n(t)} - 1| + |\frac{\a_0(t)-\a_n(t)}{\lbb_n(t)}|).
\end{align}
Hence, plugging \eqref{Rn-S-L2.1} and \eqref{Rn-S-L2.2} into \eqref{Rn-S-L2.0}
we obtain
\begin{align*}
\|w_n(t) - S_T(t)\|_{L^2}
\leq&C \big(|\frac{\lambda_n^{\frac{d}{2}}(t)-\lambda_0^{\frac{d}{2}}(t)}{\lbb_0^{\frac{d}{2}}(t)}|+|\frac{\alpha_0(t)-\alpha_n(t)}{\lbb_n(t)}|
+(\frac{\lbb_n(t)}{\lbb_0(t)})^{\frac d2}|\beta_n(t)-\beta_0(t)| \\
&\qquad +(\frac{\lbb_n(t)}{\lbb_0(t)})^{\frac d2}|\gamma_n(t)-\gamma_0(t)|
   +|\theta_n(t)-\theta_0(t)|
+ |\frac{\lbb_0(t)}{\lbb_n(t)}-1| \big),
\end{align*}
which, via \eqref{lbbn-Tt}-\eqref{thetan-Tt}, yields immediately that
\begin{align} \label{Rn-S-L2}
\|w_n(t) - S_T(t)\|_{L^2}
\leq C (T-t)^{1+ \frac \delta 2},
\end{align}
where $C$ is independent of $n$.

Estimating as above we also have that
\begin{align}  \label{Rn-S-H1}
& \|\nabla w_n(t)-\nabla S_T(t)\|_{L^2} \nonumber \\
\leq&C \big(|\frac{\lbb_0^{1+\frac{d}{2}}(t) - \lbb_n^{1+\frac{d}{2}}(t)}{ \lbb_n(t) \lbb^{1+\frac{d}{2}}_0(t)}|
            + \frac{1}{\lbb_0(t)} \frac{|\a(t)-\a_0(t)|}{ \lbb_n(t) }
            +\frac{(\lbb_n(t))^{\frac d2}}{(\lbb_0(t))^{1+ \frac d2}}|\beta_n(t) - \beta_0(t)|  \nonumber \\
   &\  + \frac{(\lbb_n(t))^{\frac d2}}{(\lbb_0(t))^{1+ \frac d2}} |\g(t) - \g_0(t)|
    + \frac{1}{\lbb_0(t)}
             |\theta_n(t) - \theta_0(t)|
        +  \frac{1}{\lbb_0(t)} |\frac{\lbb_0(t)}{\lbb_n(t)}-1|   \big) \nonumber\\
\leq& C(T-t)^\delta.
\end{align}

Thus, combining \eqref{un-St}, \eqref{esti-wn},
\eqref{Rn-S-L2} and \eqref{Rn-S-H1} altogether we obtain
\be\label{4cor}
\sup\limits_{n\geq 1}\|u_n(t)-S_T(t)\|_{H^1}\leq C(T-t)^\delta.
\ee

In particular,
this yields that
$u_n(t)-S_T(t)$ is uniformlly bounded in $H^1$ for every $t\in [0, T)$,
which along with \eqref{un-u-0-L2} implies that,
up to a subsequence (still denoted by $\{n\}$ which may depend on $t$),
\be
u_n(t)-S_T(t)\rightharpoonup u(t)-S_T(t),\ \ weakly\ in\ H^1,\ as\ n\rightarrow\infty.
\ee

Therefore,
we obtain that for every $t\in [0,T)$,
\be
\|u(t)-S_T(t)\|_{H^1}\leq \liminf_{n}\|u_n(t)-S_T(t)\|_{H^1}\leq C (T-t)^\delta,
\ee
which yields \eqref{u-S-0-H1}.
The proof of Theorem \ref{Thm-u-blowup} is complete.
\hfill $\square$

The remaining part of this section is devoted to the proof of Theorem \ref{Thm-Equiv-X-u}.
Let us first state the time regularity of $u$ below.

\begin{lemma}
let $u$ be the solution to \eqref{equa-u-RNLS} on $[0,\tau^*)$ with $u(0)=u_0 \in H^1$
in the sense of Definition \ref{def-u-RNLS}.
Then, $\bbp$-a.s. for any $0\leq s<t<\tau^*$, we have that
\begin{align} \label{Holder-u}
    \|e^{-it\Delta} u(t) - e^{-is\Delta} u(s)\|_{L^2}
    \leq C(t) |t-s|,
\end{align}
where $C(t)$ depends on $\|u\|_{C([0,t]; H^1)}$
and $\sup_{0\leq s\leq t}|B_k(s)|$,
$1\leq k\leq N$.
\end{lemma}

{\bf Proof.}
We reformulate equation \eqref{equa-u-RNLS} in the mild form
\begin{align*}
   u(t) = e^{it\Delta} u_0 + \int_0^t e^{i(t-s)\Delta}
             (i |u(s)|^{\frac 4d} u(s) + i(b(s)\cdot \na +c(s))u(s) ) ds,
\end{align*}
where $b,c$ are given by \eqref{b} and \eqref{c}, respectively.
This yields that
\begin{align*}
   \|e^{-it\Delta}u(t)-e^{-is\Delta}u(s) \|_{L^2}
   =\| \int_s^t  e^{-ir\Delta} (|u(r)|^\frac 4d u(r)  + (b(r)\cdot \na +c(r))u(r)) dr\|_{L^2}.
\end{align*}
Note that, by the Sobolev embedding $H^1 \hookrightarrow L^{2+\frac 8d}$,
\begin{align*}
   \| \int_s^t  e^{-ir\Delta} |u(r)|^\frac 4d u(r) dr\|_{L^2}
   \leq \int_s^t \|u(r)\|^{1+\frac 4d}_{L^{2+\frac 8d}}  dr
   \leq  C \|u\|^{1+\frac 4d}_{C([0,t]; H^1)}  (t-s).
\end{align*}
Moreover, in view of the boundedness of $\sup_{0\leq s\leq t}|B_k(s)|$,
$1\leq k\leq N$, we also have
\begin{align*}
   \| \int_s^t  e^{-ir\Delta} (b(r)\cdot \na +c(r))u(r) dr\|_{L^2}
   \leq& \int_s^t \|(b(r)\cdot \na +c(r))u(r)\|_{L^2} dr \\
   \leq& C(t) \|u\|_{C([0,t]; H^1)}  (t-s).
\end{align*}
Thus, combining the estimates above
we obtain \eqref{Holder-u}.
\hfill $\square$

{\bf  Proof of Theorem \ref{Thm-Equiv-X-u}.}
We mainly prove the first assertion $(i)$,
as the second assertion $(ii)$ can be proved similarly.
Below, we fix any $T\in (0,\tau^*)$
and recall that $B_k \in C^\nu([0,T])$
for any $\nu\in (\frac 13, \frac 12)$, $1\leq k\leq N$,
$\bbp$-a.s.,
and by  Remark \ref{Rem-u-weak-mild},
$$\|u\|_{C([0,T]; H^1)} + \|u\|_{L^2(0,T; H^{\frac 32}_{-1})}<\9. $$

Note that for any $\vf\in C_c^\9$ and any $0\leq s<t\leq T$,
\begin{align} \label{Delta-X-st}
     \<\delta X _{st}, \vf\>
      = \<(\del e^W)_{st} u(s), \vf\>
     + \<e^{W(s)} \delta u_{st}, \vf\>
     + \<(\delta e^W)_{st} \delta u_{st}, \vf\>.
\end{align}
We shall estimate each term on the right-hand side above separately below.

{\it $(i)$ Estimate of $\<(\del e^W)_{st} u(s), \vf\>$.}
Using Taylor's expansion we have
\begin{align}
  (\del e^W)_{st}
  = e^{W(s)} (\sum\limits_{k=1}^N i\phi_k \delta B_{k,st}
               - \frac 12 \sum\limits_{j,k=1}^N \phi_j\phi_k \delta B_{j,st}\delta B_{k,st})
               +o(t-s).
\end{align}
Then, taking into account (see \cite[Section 3.3]{FH14}, \cite[p.9]{RZZ19})
\begin{align}
   \delta B_{j,st}\delta B_{k,st}
   = \mathbb{B}_{jk,st} + \mathbb{B}_{kj,st}
     + \delta_{jk}(t-s),
\end{align}
we obtain
\begin{align}
   (\del e^W)_{st}
   = e^{W(s)} (-\mu (t-s)
               + \sum\limits_{k=1}^N i\phi_k \delta B_{k,st}
               - \sum\limits_{j,k=1}^N \phi_j\phi_k \mathbb{B}_{jk,st} )
               +o(t-s),
\end{align}
which yields that
\begin{align} \label{deltaeW-st}
   \<(\del e^W)_{st}u(s), \vf\>
   =&   \<- \mu (e^{W(s)} u(s)), \vf\>  (t-s)
      + \sum\limits_{k=1}^N \<i \phi_k  (e^{W(s)} u(s)), \vf\> \delta B_{k,st} \nonumber \\
    & - \sum\limits_{j,k=1}^N  \<\phi_j\phi_k (e^{W(s)} u(s)), \vf\> \mathbb{B}_{jk,st}
      + o(t-s).
\end{align}

{\it $(ii)$ Estimate of $\<e^{W(s)} \delta u_{st}, \vf\>$.}
Let $f(u):= |u|^\frac 4d u$.
We claim that
\begin{align} \label{deltay-st}
  \<e^{W(s)}\delta u_{st}, \vf \>
  =& \< i \Delta (e^{W(s)} u(s)), \vf\> (t-s)
     + \<i f(e^{W(s)}u(s)), \vf\> (t-s)
     + \mathcal{O}((t-s)^{1+\nu}).
\end{align}

In order to prove \eqref{deltay-st},
by \eqref{equa-u-RNLS-def}, we have
\begin{align} \label{K1-K2}
  \<e^{W(s)}\delta u_{st}, \vf \>
  =& \<e^{W(s)} \int_s^t i e^{-W(r)}\Delta (e^{W(r)} u(r)) dr, \vf\>
     + \<e^{W(s)} \int_s^t i f(u(r))dr, \vf\>  \nonumber \\
  =:& K_1 + K_2.
\end{align}
We shall treat $K_1$ and $K_2$ separately below.
For simplicity,
we set $L(r) v := (b(r)\cdot \na + c(r))v$, $v\in H^1$, $r\in [s,t]$.
Then,
since
$e^{-W(r)}\Delta (e^{W(r)} u(r)) - e^{-W(s)}\Delta (e^{W(s)} u(s))
= \Delta (u(r)-u(s)) + (L(r)u(r) - L(s)u(s))$
and $\ol{e^W} = e^{-W}$,
we have
\begin{align} \label{K1-K11-K12}
   K_1 =& \<i\Delta (e^{W(s)}u(s)), \vf\> (t-s)
          +\int_s^t \<u(r)-u(s) , (-i) \Delta  (e^{-W(s)}\vf)\> dr   \nonumber \\
        & + \int_s^t \< L(r)u(r) - L(s)u(s), (-i)e^{-W(s)}\vf\> dr  \nonumber \\
       =:& \<i\Delta (e^{W(s)}u(s)), \vf\> (t-s)
          + K_{11} + K_{12}.
\end{align}
Note that,
integration by parts formula yields that
\begin{align*}
   K_{11}
   =& \int_s^t \<e^{-ir\Delta} u(r) - e^{-is\Delta} u(s) , (-i)e^{-is\Delta} \Delta (e^{-W(s)}\vf ) \> dr \nonumber \\
    &  + \int_s^t \<u(r),(-i) (1-e^{i(r-s)\Delta})  \Delta (e^{-W(s)}\vf ) \> dr \nonumber \\
   \leq& \int_s^t \|e^{-ir\Delta} u(r) - e^{-is\Delta} u(s)\|_{L^2} \| \Delta (e^{-W(s)}\vf) \|_{L^2} dr \nonumber \\
    &  + \int_s^t \|u(r)\|_{L^2} \|(1-e^{i(r-s)\Delta})  \Delta (e^{-W(s)}\vf)\|_{L^2} dr.
\end{align*}
Since $\sup_{s\leq T}\|\partial^\nu_x (e^{-W(s)}\vf)\|_{L^2} \leq C(T)$,
$\forall 0\leq |\nu|\leq 4$,
and $\|(1-e^{i(r-s)\Delta} \psi\|_{L^2} \leq C(r-s) \|\Delta \psi\|_{L^2}$ for any $\psi \in C_c^\9$,
taking into account \eqref{Holder-u}
we obtain
\begin{align} \label{K11-st}
    K_{11} \leq&  C(T) (t-s)^2.
\end{align}

Moreover, using the integration by part formula again we have
\begin{align*}
  K_{12}
  =& \int_s^t \<e^{-ir\Delta} u(r) - e^{-is\Delta} u(s), (-i) e^{-is\Delta} L(s)^* (e^{-W(s)}\vf) \> dr  \nonumber \\
   & + \int_s^t \<u(r), (-i) (1-e^{i(r-s)\Delta})  L(s)^*  (e^{-W(s)}\vf) \>dr \nonumber \\
   & + \int_s^t \<(L(r)-L(s))u(r), (-i) e^{-W(s)} \vf\> dr.
\end{align*}
Since
$\|(L(r)-L(s))u(r)\|_{L^2} \leq C(T) \max_{1\leq k\leq N}|B_k(r)-B_k(s)| \leq C(T)(r-s)^{\nu}$,
similarly to \eqref{K11-st},
we have
\begin{align} \label{K12-st}
   K_{12}
   \leq C(T) \int_s^t   (r-s) + (r-s)^{\nu} dr
   \leq C(T) (t-s)^{1+\nu}.
\end{align}
Thus, plugging \eqref{K11-st} and \eqref{K12-st} into \eqref{K1-K11-K12}
we obtain
\begin{align} \label{K1-st}
   K_1  =  \<i\Delta (e^{W(s)}u(s)), \vf\> (t-s)
          + \mathcal{O}((t-s)^{1+\nu}).
\end{align}

Next we treat the delicate term $K_2$ in \eqref{K1-K2}.
We see that
\begin{align} \label{K2.0}
   K_2 = \<i f(e^{W(s)}u(s)), \vf\> (t-s)
         + \int_s^t \<f(u(r))- f(u(s)), (-i) e^{-W(s)} \vf\> dr.
\end{align}
Setting $u_\tau(r):= \tau u(r) + (1-\tau)u(s)$, $\tau\in [0,1]$,
we have
\begin{align*}
   f(u(r))- f(u(s))
   =& \int_0^1  \partial_zf(u_\tau(r)) d\tau  \delta u_{sr}
     +  \int_0^1 \partial_{\ol{z}} f(u_\tau(r)) d\tau \ol{\delta u_{sr}} \nonumber \\
   =:& h_1(u_\tau(r)) \delta u_{sr}  + h_2(u_\tau(r)) \ol{\delta u_{sr}},
\end{align*}
where $\partial_zf(z)=(1+\frac 2d) |z|^{\frac 4d}$,
$\partial_{\ol{z}} f(z)= \frac 2d |z|^{\frac 4d-2} z^2$,
$z\in \mathbb{C}$.
Then, by \eqref{K2.0},
\begin{align} \label{K2.1}
   K_2 =& \<i f(e^{W(s)}u(s)), \vf\> (t-s)
         + \int_s^t \<u(r)- u(s), (-i) \ol{h_1(u_\tau(r))} e^{-W(s)} \vf\> dr  \nonumber \\
        & + \int_s^t \<\ol{u(r)}-  \ol{u(s)}, (-i) \ol{h_2(u_\tau(r))} e^{-W(s)} \vf\> dr.
\end{align}

In order to estimate the right-hand side above,
we take the third term involving $h_2$ for example,
the second term containing $h_1$ can be estimated similarly.
Note that
\begin{align} \label{J1-J2}
   & \int_s^t \<\ol{u(r)}-  \ol{u(s)}, (-i) \ol{h_2(u_\tau(r))} e^{-W(s)} \vf\> dr  \nonumber  \\
  =& \int_s^t \<e^{ir\Delta} \ol{u(r)} - e^{is\Delta} \ol{u(s)}, (-i)e^{is\Delta} \ol{h_2(u_\tau(r))} e^{-W(s)} \vf\> dr  \nonumber \\
   &+ \int_s^t \< ( 1-e^{i(r-s)\Delta})\ol{u(r)}, (-i) \ol{h_2(u_\tau(r))} e^{-W(s)} \vf\> dr  \nonumber \\
  =:& K_{21} +K_{22}.
\end{align}

Since by the H\"older inequality and the Sobolev embedding $H^1 \hookrightarrow L^{\frac 8d}$,
\begin{align*}
   \|\ol{h_2(u_\tau(r))} e^{-W(s)} \vf\|_{L^2}
   \leq C \|\vf\|_{L^\9}
        \int_0^1 \|u_\tau(r)\|^{\frac 4d}_{L^\frac 8d} d\tau
   \leq C \|u\|^{\frac 4d}_{C([0,T];H^1)},
\end{align*}
which along with \eqref{Holder-u}  yields that
\begin{align} \label{J1-st}
   K_{21} \leq \int_s^t \|e^{-ir\Delta} u(r) - e^{-is\Delta} u(s)\|_{L^2}
                     \|\ol{h_2(u_\tau(r))} e^{-W(s)} \vf\|_{L^2} dr
       \leq C(T) (t-s)^2.
\end{align}

As regards $K_{22}$,
we consider the cases where $d=1$ and $d=2$ separately below.
In the case where $d=1$,
we have
\begin{align*}
   K_{22}
   \leq& \int_s^t \|1-e^{-i(r-s)\Delta}u(r)\|_{H^{-1}}
                 \|h_2(u_\tau(r))\|_{H^1} \|e^{-W(s)} \vf\|_{W^{1,\9}} dr \\
   \leq& C(T) \int_s^t (r-s)  \|h_2(u_\tau(r))\|_{H^1} dr.
\end{align*}
Note that, by H\"older's inequality and Sobolev's embedding $H^1 \hookrightarrow L^8 \cap L^\9$,
\begin{align*}
   \|h_2(u_\tau)\|_{H^1}
   \leq&  C\int_0^1(\|u_\tau^4\|_{L^2} + \|\na (u_\tau^3 \ol{u_\tau})\|_{L^2}) d\tau  \\
   \leq&  C\int_0^1(\|u_\tau\|_{L^8}^4 + \|u_\tau\|_{L^\9}^3 \|\na u_\tau\|_{L^2})d\tau
   \leq  C \|u\|^4_{C([0,T];H^1)}.
\end{align*}
This yields that, in the case where $d=1$,
\begin{align} \label{J2-d=1}
   K_{22} \leq C(T) \int_s^t (r-s) dr
   \leq C(T) (t-s)^2.
\end{align}

Moreover, in the case where $d=2$
we shall measure the spatial regularity of $u$ in the local smoothing space.
Precisely,
since $\vf\in C_c^\9$,
setting $D:=supp (\vf)$ we have
\begin{align*}
   \ol{K_{22}}
   =& \int_s^t \< \<\na\>^{-\frac 12}(( 1-e^{-i(r-s)\Delta})u(r),
                   i\<\na\>^\frac 12({h_2(u_\tau(r))} e^{W(s)} \ol{\vf}) \> dr \nonumber \\
     \leq& \int_s^t \|\chi_{D}\<\na\>^{-\frac 12}(( 1-e^{-i(r-s)\Delta} {u(r)}\|_{L^2}
                 \|{h_2(u_\tau(r))}\|_{H^{\frac 12}}
                 \|e^{W(s)}\ol{\vf}\|_{W^{1, \9}} dr \nonumber \\
     \leq& C(T) \int_s^t (r-s)\|{u(r)}\|_{H^{\frac 32}_{-1}}
                 \|{h_2(u_\tau(r))}\|_{H^{\frac 12}} dr,  \nonumber
\end{align*}
where $\chi_D$ is the characteristic function of $D$,
and  we also used the
inequalities
$\|\chi_D \<\na\>^{\frac 32} \psi\|_{L^2} \leq \|\psi\|_{H^{\frac 32}_{-1}}$
and
$\|\psi_1\psi_2\|_{H^\frac 12 } \leq C \|\psi_1\|_{W^{1,\9}} \|\psi_2\|_{H^{\frac 12}}$.
Note that,
by the  Leibniz law for fractional derivatives and
the Sobolev embedding $H^1 \hookrightarrow W^{\frac 12, 4}$,
\begin{align*}
   \|{h_2(u_\tau(r))}\|_{H^{\frac 12}}
   \leq \int_0^1 \|\<\na\>^\frac 12 ({u_\tau(r)})^2 \|_{L^2} d\tau
   \leq  \int_0^1  \|\<\na\>^\frac 12 {u_\tau(r)}  \|^2_{L^4} d \tau
   \leq   C \|u\|^2_{C([0,T];H^1)}.
\end{align*}
This yields that, in the case where $d=2$,
\begin{align} \label{J2-st}
   |K_{22}| \leq& C(T)\|u\|^2_{C([0,T];H^1)} \int_s^t (r-s)\|{u(r)}\|_{H^{\frac 32}_{-1}} dr  \nonumber \\
       \leq& C(T)  \|u\|^2_{C([0,T];H^1)} \|{u}\|_{L^2(0,T; H^{\frac 32}_{-1})} (t-s)^\frac 32.
\end{align}
Thus, plugging \eqref{J1-st}, \eqref{J2-d=1} and \eqref{J2-st} into \eqref{J1-J2}
we obtain
\begin{align} \label{h2-st}
    \int_s^t \<\ol{u(r)}- \ol{u(s)}, (-i) \ol{h_2(u_\tau(r))} e^{-W(s)} \vf\> dr
    \leq C(T) (t-s)^\frac 32.
\end{align}
The term involving $h_1(u_\tau(r))$
on the right-hand side of \eqref{K2.1} can be estimated similarly.
Thus,
we conclude that
\begin{align} \label{K2-st}
   K_2 = \<i f(e^{W(s)}u(s)), \vf\> (t-s)
          + \mathcal{O}((t-s)^\frac 32).
\end{align}

Therefore,
plugging \eqref{K1-st} and \eqref{K2-st} into \eqref{K1-K2}
we obtain \eqref{deltay-st}, as claimed.

{\it $(iii)$ Estimate of $\<(\delta e^W)_{st} \delta u_{st},\vf\>$.}
Using the integration by parts formula and H\"older's inequality we have
\begin{align}
     \<(\delta e^W)_{st} \delta u_{st}, \vf\>
   =& \int_s^t \<u(r), (-i)e^{-W(r)}\Delta(e^{W(r)} \ol{(\delta e^W)_{st}} \vf)\>
      +  \<f(u(r)), (-i) \ol{(\delta e^W)_{st}} \vf\> dr  \nonumber \\
   \leq& \int_s^t \|u(r)\|_{L^2} \|e^{W(r)} \ol{(\delta e^W)_{st}}\|_{H^2} \|\vf\|_{W^{2,\9}}
     + \|u(r)\|^{1+\frac 4d}_{L^{2+\frac 8d}} \|\ol{(\delta e^W)_{st}}\|_{L^2} \|\vf\|_{L^\9} dr.
\end{align}
Since for any multi-index $|\nu|\geq 0$,
\begin{align*}
   \|\partial_x^\nu (\delta e^W)_{st} \|_{L^\9} \leq C(T) (t-s)^{\nu},
\end{align*}
taking into account Sobolev's embedding $H^1 \hookrightarrow L^{2+\frac 8d}$,
we obtain
\begin{align}  \label{deltaeW-deltay-st}
       \<(\delta e^W)_{st} \delta u_{st}, \vf\>
  \leq C(T) (t-s)^{1+\nu}.
\end{align}

Now, plugging \eqref{deltaeW-st}, \eqref{deltay-st} and \eqref{deltaeW-deltay-st} into \eqref{Delta-X-st}
and using $X=e^W u$
we obtain
\begin{align}
   \<\delta X_{st}, \vf\>
   =& \<i \Delta X(s) + i f(X(s)) - \mu X(s), \vf\>(t-s)  \nonumber\\
    &  + \sum\limits_{k=1}^N \<i\phi_k X(s), \vf\> \delta B_{k,st}
      -  \sum\limits_{j,k=1}^N  \<\phi_j\phi_k X(s), \vf\> \mathbb{B}_{jk,st}
      + o(t-s).
\end{align}
This yields that for any $\vf \in C_c^\9$,
$\<X, \vf\> \in \mathscr{D}_B^{2\a}([0,T]; \bbr)$
with the Gubinelli derivative $Y'_k = \<i\phi_k X,\vf\>$, $1\leq k\leq N$.
In particular,
we infer that
$X:=e^W u$ satisfies equation \eqref{equa-X-rough} in the sense of Definition \ref{def-X-rough}
and \eqref{phikX-st} holds.
Therefore, the proof is complete.
\hfill $\square$

Now, by virtue of Theorem \ref{Thm-Equiv-X-u},
we obtain Theorems \ref{Thm-X-GWP} and \ref{Thm-X-blowup}
from Theorems \ref{Thm-u-GWP} and \ref{Thm-u-blowup}, respectively.

\section{Appendix} \label{Sec-Appendix}

{\bf Proof of Corollary \ref{Cor-coer-f-H1}.}
For any $f_1,f_2\in H^1$,
let
\begin{align} \label{wtf1-wtf2}
  \tilde{f}_1=f_1+a_1\cdot\nabla Q+b_1\Lambda Q+c_1\rho, \ \
  \tilde{f}_2=f_2+a_2Q+b_2\cdot xQ+c_2|x|^2Q,
\end{align}
where $a_1=(a_{1,j})$ and $b_2=(b_{2,j})$ are vectors.

Let $\wt f:= \wt f_1 + i \wt f_2$
be such that $\widetilde{f}\in \mathcal{K}$.
This yields that,
for $1\leq j\leq d$,
\begin{align*}
& a_{1,j}=\frac{2}{\|Q\|_{L^2}^2}\<f_1,xQ\>,\;
b_1=\frac{1}{\|xQ\|_{L^2}^2}\<f_1,|x|^2Q\>
    + 2\frac{\<\rho,|x|^2Q\>}{\|xQ\|_{L^2}^4}\<f_1,Q\>,\;
c_1=\frac{ 2}{\|xQ\|_{L^2}^2}\<f_1,Q\>, \\
& a_2= \frac{2}{\|xQ\|_{L^2}^2}\<f_2,\rho\>
      + 2\frac{\<\rho,|x|^2Q\>}{\|xQ\|_2^4}\<f_2,\Lambda Q\>,\;
b_{2,j}=\frac{2}{\|Q\|_{L^2}^2}\<f_2,\partial_{x_j} Q\>,\;
c_2=\frac{1}{\|xQ\|_{L^2}^2}\<f_2,\Lambda Q\>.
\end{align*}

By direct calculation we have
\begin{align*}
 \<L\tilde{f},\tilde{f}\>=&\<Lf,f\> -4b_1\<f_1,Q\>-2c_1\<f_1,|x|^2Q\> + 2b_1c_1\|xQ\|_{L^2}^2 -c_1^2\<\rho,|x|^2Q\>
 \nonumber \\
&-4\<f_2,b_2 \cdot \nabla Q\> -8c_2\<f_2,\Lambda Q\>+ |b_2|^2\|Q\|_{L^2}^2+4c_2^2\|xQ\|_{L^2}^2,
\end{align*}
which along with the inequality $ab\leq \ve a^2 + \ve^{-1} b^2$
and  Lemma \ref{Lem-coerc-L} yields
\begin{align} \label{Lwtf-wtf}
   \<Lf,f\>
    \geq  \nu \|\wt f\|^2_{H^1}  - \ve \|f\|_{H^1}^2 - C(|a_1|^2+b_1^2+c_1^2+a_2^2+|b_2|^2+c_2^2).
\end{align}
Moreover, expanding $\|\wt f\|_{H^1}^2$ by \eqref{wtf1-wtf2}
and then using the inequality
$|\<f_j, g\>_{H^1}| \leq \ve \|f_j\|^2_{H^1} + \ve^{-1} \|g\|^2_{H^1}$ for any $g\in H^1$, $j=1,2$,
we obtain that for some $C_1, C_2>0$,
\begin{align} \label{wtf-H1}
\|\tilde{f}\|_{H^1}^2
     \geq& C_1 \|f\|_{H^1}^2
           - C_2(|a_1|^2+b_1^2+c_1^2+a_2^2+|b_2|^2+c_2^2).
\end{align}

Therefore,
combining \eqref{Lwtf-wtf} and \eqref{wtf-H1}
and taking $\ve$ small enough
we obtain \eqref{coer-f-H1}.
\hfill $\square$

{\bf Proof of Corollary \ref{Cor-coer-f-local}.}
Let $\tilde{f} :=f\Phi_A^{\frac{1}{2}}$.
Since
$\nabla f\Phi_A^{\frac{1}{2}}=\nabla\tilde{f}-\frac{\nabla\Phi_A}{2\Phi_A}\tilde{f}$,
we have
\begin{align}\label{a4}
&\int|\nabla f|^2\Phi_A +|f|^2-(1+\frac{4}{d})Q^{\frac4d}f_1^2-Q^{\frac4d}f_2^2dx \nonumber \\
=&\int|\nabla \tilde{f}|^2+|\tilde{f}|^2-(1+\frac{4}{d})Q^{\frac4d}\tilde{f}_1^2-Q^{\frac4d}\tilde{f}_2^2dx \nonumber \\
&+ \int(1-\Phi_A)(|f|^2-(1+\frac{4}{d})Q^{\frac4d}f_1^2-Q^{\frac4d}f_2^2)dx \nonumber \\
&+\frac{1}{4}\int |\frac{\nabla\Phi_A}{\Phi_A} |^2|\tilde{f}|^2dx
-{\rm Re} \int\frac{\nabla\Phi_A}{\Phi_A}\cdot\nabla\tilde{f} \ol{\tilde{f}}dx
=: \sum\limits_{i=1}^4 K_i.
\end{align}

Since
$\Phi_A(x)=1$ for $|x|\leq A$ and by \eqref{ortho-cond},
$\int Qf_1dx=0$,
we infer that
\be\ba
\int Q\tilde{f}_1dx
&=\int_{|x|>A} Q(\tilde{f}_1-f_1)dx=\int_{|x|>A} \tilde{f}_1(Q-Q\Phi_A^{-\half})dx.
\ea\ee
Note that,
for $A\leq |x|\leq 2A$, $\Phi^{-\frac 12}_A(x) \leq C$
and so, by \eqref{Q-decay},
$|Q(x)-Q\Phi_A^{-\half}(x)| \leq C e^{-\delta |x|}$ for some $\delta>0$.
Moreover,
for $|x|\geq 2|A|$,
since $\Phi^{-\frac 12}_A(x) =e^{\frac{|x|}{2A}}$,
taking $A$ sufficiently large such that $\delta- \frac{1}{2A} \geq \frac \delta 2>0$,
we obtain that
$Q\Phi^{-\frac 12}_A(x) \leq C e^{-(\delta -\frac{1}{2A})|x|} \leq C e^{-\frac{\delta}{2}|x|}$.
Thus, we conclude that
$|Q(x)-Q\Phi_A^{-\half}(x)| \leq C e^{-\frac \delta 2 |x|}$ for $|x|\geq A$.
This yields that
\begin{align}
  |\int Q\tilde{f}_1dx |
  \leq& \|\wt f_1\|_{L^2} (\int_{|x|>A} |Q-Q\Phi_A^{-\frac 12}|^2 dx)^\frac 12 \nonumber \\
  \leq& C  \|\wt f_1\|_{L^2} (\int_{|x|>A} e^{-\delta |x|} dx)^\frac 12
  =: \delta_A \|\wt f_1\|_{L^2},
\end{align}
where  $\delta_A \to 0$ as $A\to \9$.
In view of the orthogonal conditions
and the exponential decay of $Q$ and $\rho$,
similar arguments as above also yield
that the remaining five inner products in the second part
on the right-hand side of \eqref{coer-f-H1} can be also bounded
by $\delta_A \|\wt f\|_{L^2}$ with $\delta_A \to 0$ as $A\to \9$.
Thus, Corollary \ref{Cor-coer-f-H1} yields that
\be\label{a1}
 K_1 \geq \nu_1\|\tilde{f}\|_{H^1}^2- \nu_2 \delta_A \|\tilde{f}\|_{L^2}^2.
\ee

Moreover, by \eqref{Q-decay},  we have that for $A$ large enough,
\be\label{a2}
 K_2=\int_{|x|>A}(1-\Phi_A)(|f|^2-(1+\frac 4d)Q^2f_1^2-Q^2f^2_2)dx
>0.
\ee

We also note that,
since $|\frac{\nabla\Phi_A}{\Phi_A}|\leq C A^{-1}$,
by  H\"{o}lder's inequality,
\be\label{a3}
 K_3 + K_4 \leq \frac C A \|\tilde{f}\|_{H^1}^2.
\ee

Thus, plugging (\ref{a1}), (\ref{a2}) and (\ref{a3}) into (\ref{a4}),
we obtain  that for $A$ large enough
\begin{align} \label{esti-wtf}
  \int|\nabla f|^2\Phi_A +|f|^2-(1+\frac{4}{d})Q^{\frac4d}f_1^2-Q^{\frac4d}f_2^2dx
  \geq \frac{\nu_1}{2} \|\wt f\|_{H^1}^2.
\end{align}
Using again
$\nabla\tilde{f} = \nabla f\Phi_A^{\frac{1}{2}}+ \frac{\nabla\Phi_A}{2\Phi_A}\tilde{f}$
we see that
\begin{align} \label{esti-wtf.1}
   \frac{\nu_1}{2}\|\wt f\|_{H^1}^2
   =&  \frac{\nu_1}{2} \int |f|^2 \Phi_A dx
     + \frac{\nu_1}{2} \int |\nabla f\Phi_A^{\frac{1}{2}}+ \frac{\nabla\Phi_A}{2\Phi_A}\tilde{f}|^2 dx \nonumber \\
   =&   \frac{\nu_1}{2} \int (|f|^2 + |\na f|^2) \Phi_A dx
      +  \frac{\nu_1}{2} \int {\rm Re} (\na f \Phi^\frac 12_A \frac{\na \Phi_A}{\Phi_A} \ol{\wt f}) dx
      +  \frac{\nu_1}{8} \int |\frac{\na \Phi_A}{\Phi_A} \wt f|^2 dx.
\end{align}
By H\"older's inequality and $|\frac{\nabla\Phi_A}{\Phi_A}|\leq C A^{-1}$,
\begin{align} \label{esti-wtf.1}
   |\frac{\nu_1}{2} \int {\rm Re} (\na f \Phi^\frac 12_A \frac{\na \Phi_A}{\Phi_A} \ol{\wt f}) dx |
   \leq& \frac{\nu_1C}{2A} (\int |\na f|^2 \Phi_A dx)^{\frac 12}
        (\int | f|^2 \Phi_Adx)^\frac 12 \nonumber  \\
   \leq& \frac{\nu_1}{4} \int |\na f|^2 \Phi_A dx
         + \frac{4\nu_1 C^2}{A^2} \int |f|^2 \Phi_A dx.
\end{align}
Moreover, we have
\begin{align} \label{esti-wtf.2}
   \frac{\nu_1}{8} \int |\frac{\na \Phi_A}{\Phi_A} \wt f|^2 dx
  \leq \frac{\nu_1 C^2}{2A^2} \int |f|^2 \Phi_A dx.
\end{align}
Therefore,
putting together \eqref{esti-wtf}-\eqref{esti-wtf.2}
and taking $A$ large enough
we obtain  \eqref{coer-f-local}.
\hfill $\square$

{\bf Proof of Lemma \ref{Lem-imp-thm}.}
The proof is  based on the implicit function theorem.
Let $w$ and $R$ be as in \eqref{dec-u-Ur}.
Set $\mathcal{P} := (\lbb, \a, \b, \g,\t)$
and define the functionals
\begin{align*}
&f_{1,j}(u,\mathcal{P}) ={\rm Re}\int (x_j-\alpha_j)w \ol{R}dx, \ \
f_2(u, \mathcal{P}) ={\rm Re}\int |x-\alpha|^2 w \ol{R}dx, \\
&f_{3}(u, \mathcal{P} ) ={\rm Im}\int (\frac d2 w+(x-\alpha)\cdot \nabla w) \ol{R}dx,\ \
f_{4,j}(u, \mathcal{P}) ={\rm Im}\int \partial_{x_j} w \ol{R}dx, \\
&f_5(u, \mathcal{P}) ={\rm Im}\int \wt \rho \ol{R}dx, \ \ 1\leq j\leq d.
\end{align*}
Note that,
by the definition of $u_0$,
$ f_j(u_0, \mathcal{P}_0)=0$,
$1\leq j\leq 5$.
Moreover,
the orthogonality conditions and
straightforward computations show that for $1\leq j,k\leq d$,
\begin{align*}
  &\pa_{\a_k} f_{1,j} (u_0,\calp_0) = -\frac {\delta_{jk}}{2}  \|Q\|_{L^2}^2 + \mathcal{O}(\|R_0\|_{L^2}),\
   \pa_{\lbb} f_2  (u_0,\calp_0) = -\lbb_0 \|xQ\|_{L^2}^2 + \mathcal{O}(\lbb_0 \|R_0\|_{L^2}), \\
  & \pa_{\a_k} f_3  (u_0,\calp_0) = - \frac{ \beta_{0,k}}{2}\|Q\|_{L^2}^2 + \mathcal{O}(\frac{ \|R_0\|_{L^2}}{\lbb_0}),\
    \pa_{\g} f_3 (u_0,\calp_0) = \frac 14 \|xQ\|_{L^2}^2 + \mathcal{O}(\|R_0\|_{L^2}), \\
  &  \pa_{\lbb} f_{4,j}  (u_0,\calp_0) = \frac {\beta_{0,j}}{2\lbb_0^{2}} \|Q\|^2_{L^2} + \mathcal{O}(\frac{ \|R_0\|_{L^2}}{\lbb_0^{2}}),\
   \pa_{\beta_k} f_{4,j}  (u_0,\calp_0) = - \frac {\delta_{jk}}{ 2 \lbb_0} \|Q\|^2_{L^2} + \mathcal{O}(\frac{\|R_0\|_{L^2}}{\lbb_0}), \\
  & \pa_{\lbb} f_5 (u_0,\calp_0) =  \frac{\g_0}{2 \lbb_0}  \<\rho, |x|^2 Q\> + \mathcal{O}(\frac{\|R_0\|_{L^2}}{\lbb_0}) ,  \
    \pa_{\a_k} f_5 = \frac {\beta_{0,k}}{2 \lbb_0}  \|xQ\|_{L^2}^2 +  \mathcal{O}(\frac{\|R_0\|_{L^2}}{\lbb_0}), \\
  & \pa_{\g} f_5 = -\frac 14\<\rho, |x|^2Q\>  +  \mathcal{O}(\|R_0\|_{L^2}) ,   \
    \pa_{\theta} f_5 = -\frac 12 \|xQ|_{L^2}^2 +  \mathcal{O}(\|R_0\|_{L^2}) .
\end{align*}
and the remaining partial differentials are bounded at most by $C \lbb_0^{-2} \|R_0\|_{L^2} = \mathcal{O}(T)$ for some $C>0$.
This yields that
the Jacob determinant
\be\label{Jacob}
\frac{\partial (f_1, f_2, f_3, f_4, f_5)}{\partial \mathcal{P}} |_{(u_0, \mathcal{P}_0)}
= 2^{-2d-3} \lbb_0^{1-d} \|Q\|_{L^2}^{4d} \|xQ\|_{L^2}^6 + \mathcal{O}(T),
\ee
which is positive for $T$ small enough.

Thus,
the  implicit function theorem
yields the existence of $r_0>0$ and a $C^1$ map $\Psi$ from $\calu_{r_0}(u_0)$
to $\calu_{r_0}(\calp_0)$
and the orthogonality conditions in \eqref{oc2} hold.
\hfill $\square$

\section*{Acknowledgements}
The authors would like to thank Prof. Daomin Cao for valuable comments to improve the manuscript.
Y. Su is supported by NSFC (No. 11601482) and D. Zhang  is supported by NSFC (No. 11871337).


\begin{thebibliography}{99}

\bibitem{BCIR94}
O. Bang, P.L. Christiansen, F. If, K.O. Rasmussen,
Temperature effects in a nonlinear model of monolayer Scheibe aggregates.
{\it Phys. Rev. E} {\bf 49} (1994), 4627--4636.

\bibitem{BCIRG95}
O. Bang, P.L. Christiansen, F. If, K.O. Rasmussen, Y.B. Gaididei,
White noise in the two-dimensional nonlinear Schr\"odinger equation,
{\it Appl. Anal}. {\bf 57} (1995), no. 1-2, 3--15.

\bibitem{BCD}
V. Banica, R. Carles, T. Duyckaerts, Minimal blow-up solutions to the mass-critical inhomogeneous NLS equation. {\it Commun. Partial Differ. Equ.}
{\bf 36} (2011),  no. 3, 487-531.


\bibitem{BR17}
V. Barbu, M. R\"ockner,
Global solutions to random 3D vorticity equations for small initial data.
{\it J. Differential Equations} {\bf 263} (2017), no. 9, 5395--5411.


\bibitem{BRZ14} V. Barbu, M. R\"ockner, D.Zhang,
The stochastic nonlinear Schr\"odinger equation with multiplicative noise:
the rescaling aproach, {\it J. Nonlinear Sciences}, {\bf 24} (2014), 383--409.

\bibitem{BRZ16}
V. Barbu, M. R\"ockner,  D. Zhang,
Stochastic nonlinear Schr\"odinger equations.
{\it Nonlinear Anal}. {\bf 136} (2016), 168--194.

\bibitem{BRZ17.0}
V. Barbu, M. R\"ockner, D. Zhang,
The stochastic logarithmic Schr\"odinger equation.
{\it J. Math. Pures Appl. (9)} {\bf 107} (2017), no. 2, 123--149.

\bibitem{BRZ17}
V. Barbu, M. R\"ockner, D. Zhang,
Stochastic nonlinear Schr\"odinger equations: no blow-up in the non-conservative case.
{\it J. Differential Equations} {\bf 263} (2017), no. 11, 7919--7940.

\bibitem{BRZ18}
V. Barbu, M. R\"ockner, D. Zhang,
Optimal bilinear control of nonlinear stochastic Schr\"odinger equations driven by linear multiplicative noise.
{\it Ann. Probab}. {\bf 46} (2018), no. 4, 1957--1999.


\bibitem{BG09}
A. Barchielli, M. Gregoratti,
Quantum Trajectories and Measurements in Continuous Case.
The Diffusive Case,
{\it Lecture Notes Physics} {\bf 782}, Springer Verlag,
Berlin, 2009.


\bibitem{BM13}
Z. Brze\'{z}niak, A. Millet,
On the stochastic Strichartz estimates and the stochastic nonlinear Schr\"odinger equation on a compact Riemannian manifold.
{\it Potential Anal}. {\bf 41} (2014), no. 2, 269--315.

\bibitem{B-W}
J. Bourgain, W. Wang,
Construction of blowup solutions for the nonlinear Schr\"odinger equation with critical nonlinearity.
Dedicated to Ennio De Giorgi.
{\it Ann. Scuola Norm. Sup. Pisa Cl. Sci. (4)} {\bf 25} (1997), no. 1--2, 197--215


\bibitem{CS13}
D.M. Cao, Y. Su,
Minimal blow-up solutions of mass-critical inhomogeneous Hartree equation. {\it J. Math. Phys.} {\bf 54} (2013), no. 12, 121511, 25pp.


\bibitem{C}
T. Cazenave,  Semilinear Schr\"odinger equations.
{\it Courant Lecture Notes in Mathematics}, 10. New York University,
Courant Institute of Mathematical Sciences,
 New York; American Mathematical Society, Providence, RI, 2003. xiv+323 pp.

\bibitem{BD02}
A. de Bouard, A. Debussche,
On the effect of a noise on the solutions of the focusing supercritical nonlinear Schr\"odinger equation.
{\it Probab. Theory Related Fields} {\bf 123} (2002), no. 1, 76--96.

\bibitem{BD03}
A. de Bouard, A. Debusche,
The stochastic nonlinear Schr\"odinger equation in $H^1$,
{\it Stoch. Anal. Appl.}, {\bf 21} (2003), 97--126.

\bibitem{BD05}
A. de Bouard, A. Debussche,
Blow-up for the stochastic nonlinear Schr\"odinger equation with multiplicative noise.
{\it Ann. Probab}. {\bf 33} (2005), no. 3, 1078--1110.

\bibitem{BDM01}
A. de Bouard, A. Debussche, L. Di Menza,
Theoretical and numerical aspects of stochastic nonlinear Schr\"odinger equations.
{\it Journ\'{e}es "\'{E}quations aux D\'{e}riv\'{e}es Partielles"} (Plestin-les-Gr\`{e}ves, 2001),
Exp. No. III, 13 pp., Univ. Nantes, Nantes, 2001.

\bibitem{DM02}
A. Debussche, L. Di Menza,
Numerical simulation of focusing stochastic nonlinear Schr\"odinger equations.
{\it Phys. D} {\bf 162} (2002), no. 3--4, 131--154.

\bibitem{DM02.2}
A. Debussche, L. Di Menza,
Numerical resolution of stochastic focusing NLS equations.
{\it Appl. Math. Lett}. {\bf 15} (2002), no. 6, 661--669.

\bibitem{F17}
C.J. Fan,
log-log blow up solutions blow up at exactly m points.
{\it Ann. Inst. H. Poincar\'{e} Anal. Non Lin\'{e}aire} {\bf 34} (2017), no. 6, 1429--1482.

\bibitem{FX18.1}
C.J. Fan, W.J. Xu,
Global well-posedness for the defocusing mass-critical stochastic nonlinear Schr\"odinger equation on $\bbr$ at $L^2$ regularity,
\texttt{arXiv:1810.07925}

\bibitem{FX18.2}
C.J. Fan, W.J. Xu,
Subcritical approximations to stochastic defocusing mass-critical nonlinear Schr\"odinger equation on $\bbr$.
{\it J. Differential Equations} {\bf 268} (2019), no. 1, 160--185.

\bibitem{FH14}
P. Friz, M. Hairer,
A course on rough paths.
With an introduction to regularity structures. {\it Universitext}. Springer, Cham, 2014. xiv+251 pp.

\bibitem{G04}
M. Gubinelli, Controlling rough paths. {\it J. Funct. Anal.} {\bf 216} (2004), no. 1, 86--140.

\bibitem{HRZ18}
S. Herr, M. R\"ockner, D. Zhang,
Scattering for stochastic nonlinear Schr\"odinger equations.
{\it Comm. Math. Phys}. {\bf 368} (2019), no. 2, 843--884.


\bibitem{H18}
F. Hornung, The nonlinear stochastic Schr\"oodinger equation via stochastic
Strichartz estimates. {\it J. Evol. Equ.} {\bf 18} (2018), no. 3, 1085--1114.

\bibitem{K-L-R}
J. Krieger, E. Lenzmann, P. Rapha\"{e}l, Nondispersive solutions to the $L^2$-critical halfwave equation, {\it Arch. Ration. Mech. Anal.} {\bf 209} (2013), no. 1,  61-129.

\bibitem{K-M-R}
J. Krieger, Y. Martel, P. Rapha\"{e}l,
Two-soliton solutions to the three-dimensional gravitational Hartree equation.
{\it Comm. Pure Appl. Math.} {\bf 62} (2009), no. 11, 1501--1550.

\bibitem{MM02}
Y. Martel, F. Merle,
Nonexistence of blow-up solution with minimal $L^2$-mass for the critical gKdV equation.
{\it Duke Math. J}. {\bf 115} (2002), no. 2, 385--408.

\bibitem{MP18}
Y. Martel, P. Rapha\"el,
Strongly interacting blow up bubbles for the mass critical nonlinear Schr\"odinger equation.
Ann. Sci. \'{E}c. Norm. Sup\'{e}r. (4) 51 (2018), no. 3, 701-737.


\bibitem{M96}
F. Merle,
Nonexistence of minimal blow-up solutions of equations
$iu_t = -\Delta u-k(x)|u|^{4/N}u$ in $\bbr^N$.
{\it Ann. Inst. H. Poincar\'{e} Phys. Th\'{e}or}. {\bf 64} (1996), no. 1, 33--85.

\bibitem{Mu}
F. Merle,
Determination of blow-up solutions with minimal mass for nonlinear Schr\"odinger equations with critical power.
{\it Duke Math. J.} {\bf 69} (1993), no. 2, 427--454.

\bibitem{M-R}
F. Merle, P. Rapha\"{e}l,
The blow-up dynamic and upper bound on the blow-up rate for critical nonlinear Schr\"odinger equation.
{\it Ann. of Math. (2)} {\bf 161} (2005), no. 1, 157--222.

\bibitem{P}
G. Perelman, On the blow-up phenomenon for the critical nonlinear Schr\"{o}dinger equation in 1D,
{\it Ann. Henri. Poincar\'{e}} {\bf 2} (2001), 605--673.

\bibitem{PR07}
F. Planchon, P. Rapha\"el,
Existence and stability of the log-log blow-up dynamics for the $L^2$-critical nonlinear Schr\"odinger equation in a domain.
{\it Ann. Henri Poincar\'{e}} {\bf 8} (2007), no. 6, 1177--1219.

\bibitem{R-S}
P. Rapha\"{e}l, J. Szeftel,
Existence and uniqueness of minimal blow-up solutions to an inhomogeneous mass critical NLS.
{\it J. Amer. Math. Soc.} {\bf 24} (2011), no. 2, 471--546.


\bibitem{RGBC95}
K.O. Rasmussen, Y.B. Gaididei, O. Bang, P.L. Chrisiansen,
The influence of noise on critical collapse
in the nonlinear Schr\"odinger equation.
{\it Phys. Letters A} {\bf 204} (1995), 121--127.

\bibitem{RZZ19}
M. R\"ockner, R.C. Zhu, X.C. Zhu,
A remark on global solutions to random 3D vorticity equations for small initial data.
{\it Discrete Contin. Dyn. Syst. Ser. B} {\bf 24} (2019), no. 8, 4021--4030.


\bibitem{S17}
Y. Su, Uniqueness of minimal blow-up solutions to nonlinear Schr\"{o}dinger system.
{\it Nonlinear Anal.} {\bf 155} (2017), 186--197.

\bibitem{SG19}
Y. Su, Q. Guo,
Blow-up solutions to nonlinear Schr\"{o}dinger system at multiple points.
{\it Z. Angew. Math. Phys.} {\bf 70} (2019), no. 1, Art. 70, 14pp.


\bibitem{SS99}
C. Sulem, P.L. Sulem,
The nonlinear Schr\"odinger equation: self-focusing
and wave collapse.
{\it Applied Mathematical Sciences} {\bf 139}, Springer, New York,
1999.

\bibitem{T06}
T. Tao, Nonlinear dispersive equations. Local and global analysis.
CBMS Regional Conference Series in Mathematics, 106.
Published for the Conference Board of the Mathematical Sciences, Washington, DC;
by the American Mathematical Society, Providence, RI, 2006. xvi+373 pp.

\bibitem{Wenn}
M. Weinstein,
Nonlinear Schr\"odinger equations and sharp interpolation estimates.
{\it Comm. Math. Phys.} {\bf 87} (1982/83), no. 4, 567--576.

\bibitem{Wenm}
M. Weinstein, Modulational stability of ground states of nonlinear Schr\"{o}dinger equations.
{\it SIAM J. Math. Anal.} {\bf 16} (1985), no. 3, 472--491.

\bibitem{Z18}
D. Zhang, Stochastic nonlinear Schr\"odinger equations in
the defocusing mass and energy critical cases.
\texttt{arXiv:1811.00167v2}

\bibitem{Z19}
D. Zhang, Optimal bilinear control of stochastic
nonlinear Schr\"odinger equations:
mass-(sub)critical case,
arXiv: 1902.03559.



\end{thebibliography}
\end{document}